\def\Box{{$\sqcap$}\llap{$\sqcup$}}

 \magnification=1200 \tolerance=1000 \overfullrule=0pt

 \font\nagytwobf=cmbx10 scaled\magstep2 
\font\nagybf=cmbx10
 scaled\magstep1 
\font\sf=cmss10 
\def\Z{\hbox{\sf Z\kern-.34em Z}}
 \def\R{\hbox{\sf I\kern-.34em R}} 
\def\N{\hbox{\sf I\kern-.34em N}}
\def\Q{\hbox{\sf Q\kern-.7em Q}}

 \def\proofend{\nobreak\Box\par\vskip5pt} 

 \def\l{\ell}

 \def\eps{\varepsilon}

 \def\ph{\varphi}

\centerline{\nagytwobf An attempt to prove an effective Siegel theorem  }
\bigskip
\centerline{\nagytwobf Part One}
\bigskip
\bigskip\bigskip
\centerline{\bf J\'ozsef Beck (Rutgers, USA)}
\bigskip\bigskip\bigskip
We describe a plan how to prove an effective Siegel theorem 
(about the exceptional Dirichlet character). 
The basic idea is quite simple: it is briefly outlined in Section 0.
We completed a proof, but it is ridiculously long; more than 300
pages. Such a long paper probably contains several errors; the
question is whether the errors are substantial or not.
Even if there are errors, is the basic idea still good? 
In Sections 1-5 we give a very detailed
plan for the case of negative discriminants. The missing details
(mostly routine elementary estimations) are in Part Two.
I am happy to send the pdf-file of Part Two to anybody who 
requests it by email.

\bigskip

AMS Math. Subject Classification: 11E41, 11F66, 11N13

Key words and phrases: Dirichlet L-functions, Siegel zero,
class number
\bigskip\bigskip
\medskip
\centerline{\bf 0. A nutshell summary }
\bigskip\bigskip\bigskip\medskip\noindent
We begin with the case of negative discriminants; we want to prove
\medskip
\noindent
{\bf Theorem 1:} {\it For every negative
fundamental discriminant $\Delta<0$ with $|\Delta|>e^{10^{100}}$,
$$L(1,\chi_{\Delta})\ge {1\over (\log |\Delta|)^{146}}.\eqno(0.1)$$}
\medskip\noindent
{\bf Brief outline of the proof.}
Assume that there is a ``bad'' negative fundamental discriminant 
$\Delta =-D<0$ ($D=|\Delta |$ is always positive)
such that
$$L(1,\chi_{-D})< {1\over (\log D)^{146}}\ \ {\rm and} \ \ 
D>e^{10^{100}}.\eqno(0.2)$$
By the class number formula and (0.2), the
class-number $h(-D)$ is relatively small:
$$h(-D)<\sqrt{D}(\log D)^{-146}.$$
Our goal is to derive a contradiction.

First we define a simple ``sieve out'' procedure; we distinguish two cases.
Let $p_0$ denote the smallest prime that is not a divisor of
the discriminant $D$; clearly $p_0<2\log D$.
If $\chi_{-D}(p_0)=1$ then let
$$Z_0=Z'_0=\prod_{\scriptstyle p\le D:\atop\scriptstyle
\chi_{-D}(p)\ne -1}p\eqno(0.3)$$
be the product of all primes $p\le D$
with $\chi_{-D}(p)=1$ or 0.

If $\chi_{-D}(p_0)=-1$ then we have two possible choices
for $Z_0$: either (3) (i.e., $Z_0=Z'_0$),
or we include $p_0$ as an extra prime factor:
$Z_0=Z''_0=p_0Z'_0$.
One of these two choices will lead to a contradiction.

Let $m\ge 1$ be a squarefree integer relatively prime to
$Z_0$, and let
 $\l $ be an integer in $1\le\l <m$ relatively
prime to $m$.
Given $0<A_1<A_2$, let
${\bf v}_{m,\l ;a}(A_1,A_2)$ denote the  vector
$${\bf v}_{m,\l ;a}(A_1,A_2)=$$
$$=(e^{2\pi {\rm i}(j_1D+a)\l /m},
e^{2\pi {\rm i}((j_1+1)D+a)\l /m},
e^{2\pi {\rm i}((j_1+2)D+a)\l /m},
\ldots , e^{2\pi {\rm i}(j_2D+a)\l /m})\eqno(0.4)$$
that plays a key role in the proof;
here $j=j_1$ is the smallest integer with $jD+a\ge A_1$, 
 $j=j_2$ is the largest integer with $jD+a\le A_2$, and
of course ${\rm i}=\sqrt{-1}$.

We need the usual inner product
of complex vectors: if $\psi =(a_1,\ldots ,a_n)$ and
$\xi =(b_1,\ldots ,b_n)$ with $a_j,b_j$  complex numbers, then
$\langle \psi ,\xi\rangle =\sum_{j=1}^na_j\overline{b_j}$
where $\overline{b_j}$ is the complex conjugate of $b_j$.
We write, as usual, 
$$\| \psi\| =\sqrt{\langle\psi ,\psi\rangle },\ \ {\rm so\ \ }
\| \psi\|^2 =\langle\psi ,\psi\rangle $$
is the square of the norm. For arbitrary complex numbers
$C_j$ we have 
$$\left\|\sum_{j=1}^LC_j\psi_j\right\|^2=
{\rm Diagonal}+{\rm OffDiagonal},\eqno(0.5)$$
where
$${\rm Diagonal}=\sum_{j=1}^L|C_j|^2\|\psi_j\|^2
{\rm\ and\ }
{\rm OffDiagonal}=\sum_{j=1}^L\sum_{\scriptstyle 1\le k\le L:\atop
\scriptstyle k\ne j}
C_j\overline{C_k}\langle\psi_j,\psi_k\rangle .$$
We apply the quadratic identity (0.5) with the choice
$$\psi_j={\bf v}_{m,\l ;a}(A_1,A_2){\rm\ \ and\ \ }
C_j=C_{m,\l ;a}={\mu (m)\over\ph (m)}\eqno(0.6)$$
We add up the equalities (0.5) with the choice of (0.6) 
for all $1\le a\le D$ with weight
$\chi_{-D}(a)$.

Later in the proof we apply a routine ``smoothing'' that I
completely skip; here
for simplicity we just choose  $A_1=1$ and $A_2=ND$.

Write
$${\bf V}_a(M;A_1,A_2)
=\sum_{\scriptstyle 1\le m\le M,1\le \l\le m:\atop
\scriptstyle \gcd (m,\l Z_0)=1}
{\mu (m)\over\ph (m)}{\bf v}_{m,\l ;a}(A_1,A_2).$$
By (0.4),
$$\left\| {\bf V}_a(M;A_1,A_2)\right\|^2=
\sum_{j:\ A_1\le jD+a\le A_2}\left|
\sum_{\scriptstyle 1\le m\le M:
\atop\scriptstyle\gcd (Z_0, m)=1}{\mu (m)\over\ph (m)}
\sum_{\scriptstyle 1\le \l\le m:\atop
\scriptstyle \gcd (\l ,m)=1}e^{2\pi {\rm i}(a+jD)\l /m}\right|^2
.\eqno(0.7)$$
Using the well-known Ramanujan's sum 
$$\sum_{\scriptstyle 1\le \l<m\atop\scriptstyle
\gcd (\l ,m)=1}
e^{2\pi {\rm i}n\l /m}={\mu (m')\ph (m)\over\ph (m')}
$$
in (0.7), where $m'=m/\gcd (n,m)$ and $\ph (m)$
is the  Euler's function, we have
$$\sum_{a=1}^D\chi_{-D}(a)\left\| {\bf V}_a(M; A_1,A_2)\right\|^2=$$
$$=\sum_{a=1}^D\chi_{-D}(a)
\sum_{j:\ A_1\le jD+a\le A_2}
\left(\sum_{\scriptstyle d\ge 1:d|a+jD\atop\scriptstyle\gcd 
(Z_0,d)=1}\mu (d)
\sum_{\scriptstyle 1\le k\le M/d:\atop
\scriptstyle \gcd (k,(a+jD)Z_0)=1}
{|\mu (k)|\over\ph (k)}\right)^2=$$
$$=\prod_{p|Z_0}\left( 1-{1\over p}\right)^2\cdot
D\log {M^2\over N}
+{\rm negligible\ error\ term},\eqno(0.8)$$
where the last step requires some routine (mostly elementary)
 estimations, using the following choices of
 the key parameters $M$ and $N$: let 
$$\log N=(\log D)^{15},\ \ M=N\exp\left(-{2\over 3}
\sqrt{\log N}\right) .\eqno(0.9)$$
On the other hand, we have
$$\left\| {\bf V}_a(M;A_1,A_2)\right\|^2={\rm Diagonal}_a(M;A_1,A_2)+{\rm
OffDiagonal}_a(M;A_1,A_2),\eqno(0.10)$$
where
$${\rm OffDiagonal}_a(M;A_1,A_2)=
\sum_{\scriptstyle d\ge 1:\atop\scriptstyle \gcd (d,Z_0)=1}
\sum_{{\scriptstyle (m_1,\l_1)\ne 
(m_2,\l_2):\ 1\le m_1,m_2\le M\atop
\scriptstyle 1\le\l_h\le m_h,\gcd (\l_h,m_h)=1,h=1,2}\atop\scriptstyle
\gcd (m_1,m_2)=d,\gcd (m_1m_2,Z_0)=1}{\mu(m_1)\mu(m_2)
\over\ph (m_1)\ph (m_2)}\cdot $$
$$\qquad\qquad\qquad\cdot
e^{2\pi{\rm i}a({\l_1\over m_1}-{\l_2\over m_2})}
\cdot \sum_{j:\ A_1\le jD+a\le A_2} 
e^{2\pi{\rm i}Dj({\l_1\over m_1}-{\l_2\over m_2})} .\eqno(0.11)$$
Note that
$$\sum_{a=1}^D\chi_{-D}(a)
\left\| {\bf V}_a(M;A_1,A_2)\right\|^2=\sum_{a=1}^D\chi_{-D}(a){\rm
OffDiagonal}_a(M;A_1,A_2),\eqno(0.12)$$
since the diagonal part clearly cancels out.
The critical sum (see the end of (0.11))
$$\sum_{a=1}^D\chi_{-D}(a)
e^{2\pi{\rm i}a({\l_1\over m_1}-{\l_2\over m_2})}
\cdot \sum_{j:\ A_1\le jD+a\le A_2} 
e^{2\pi{\rm i}Dj({\l_1\over m_1}-{\l_2\over m_2})} \eqno(0.13)$$
can be estimated by involving the Gauss sums corresponding to the 
exceptional character $\chi_{-D}$ (and some routine approximation; see
Section 3).
The explicit formula for these Gauss sums brings in the
extra factor $\chi_{-D}(m_1)\chi_{-D}(m_2)$ involving the 
exceptional character $\chi_{-D}$. The crucial fact is that
now in (0.11) we have the product
$$\mu(m_1)\mu(m_2)\cdot \chi_{-D}(m_1)\chi_{-D}(m_2),$$
which is $\mu^2(m_1)\mu^2(m_2)$ for the overwhelming majority
of the pairs $(m_1,m_2)$ showing up in the sum in (0.11)
(indeed, the assumption of small class number $h(-D)$
implies that $\mu(m)=\chi_{-D}(m)$ for the overwhelming majority
of square-free integers $1\le m\le M$ with $\gcd (m,Z_0)=1$).
Since $\mu^2(m_1)\mu^2(m_2)=1$ or 0 is {\it positive}, 
(0.11) becomes a relatively simple Riemann sum 
(apart from a negligible error), and we can approximate (0.12) with the
corresponding definite integral.
To illustrate what is going on here, consider the sum
$$\sum_{j=1}^{D^4}\chi_{-D}(j)\sum_{k=1}^{M/j}{\log k\over k}=
\sum_{j=1}^{D^4}\chi_{-D}(j)\int_{x=1}^{M/j}{\log x\over x}\, dx+
{\rm negligible}=$$
$$=\sum_{j=1}^{D^4}\chi_{-D}(j){1\over 2}\log^2 (M/j)+
{\rm negligible}=$$
$$=\sum_{j=1}^{D^4}\chi_{-D}(j){1\over 2}\left(\log^2 M
-2\log M\log j+\log^2 j\right) +{\rm negligible}=$$
$$=-\log M\sum_{j=1}^{D^4}\chi_{-D}(j)\log j+{\rm negligible}
,\eqno(0.14)$$
where we used the following facts: $\sum_{j=1}^{D^4}\chi_{-D}(j)=0$,
$$\sum_{j=1}^{D^4}\chi_{-D}(j)\log^2 j=O(\sqrt{D}\log^3 D)$$
(P\'olya--Vinogradov plus partial summation),
$\log M$ is much larger than $\log D$ (see (0.9)), and finally
we used (0.16) below.

By using (0.11)-(0.13) and an argument similar to (0.14),
we eventually obtain
$$\sum_{a=1}^D\chi_{-D}(a)
\left\| {\bf V}_a(M;A_1,A_2)\right\|^2=\left(\sum_{j=1}^{D^4}
\chi_{-D}(j)\log j\right)\cdot {\rm RoutineFactor}=$$
$$=\prod_{p|Z_0}\left(1-{1\over p}\right)^2
\prod_{\scriptstyle p|Z_0\atop\scriptstyle \chi_{-D}(p)=1}
\left( 1+{1\over p-1}\right)
\prod_{\scriptstyle p|Z_0\atop\scriptstyle \chi_{-D}(p)=-1}
\left( 1-{1\over p+1}\right)
\cdot D\log {M^2\over N},\eqno(0.15)$$
where ``RoutineFactor'' does {\it not} distinguish
 between the primes
$p|Z_0$ with $\chi_{-D}(p)=1$, $0$ or $-1$
(the power $D^4$ is an accidental choice and we ignored the negligible
 error).
Comparing (0.8) and (0.15) we obtain that they are not equal---a 
contradiction that proves (0.1). They are not equal, because
(0.8) does {\it not} distinguish, but (0.15) {\it does}
distinguish between the primes
$p|Z_0$ with $\chi_{-D}(p)=1$, $0$ or $-1$
due to the sum $\sum_{j=1}^{D^4}\chi_{-D}(j)\log j$.
Here the first factor $\chi_{-D}(j)$ comes from the Gauss sum, and
the second factor $\log j$ is explained by the fact that (0.11)
resembles a ``double harmonic sum'', which is a
Riemann sum for a logarithmic integral that we can evaluate
explicitly (somewhat like (0.14)). 

Furthermore, if the class number $h(-D)$ is ``substantially smaller'' than
$\sqrt{D}$, then we have the good approximation
$$\sum_{j=1}^{D^4}\chi_{-D}(j)\log j=-{\pi\over 6}\sqrt{D}
\sum_{(a,b,c)}{1\over a}\ +\ {\rm negligible},\eqno(0.16)$$
where $\sum_{(a,b,c)}{1\over a}$ means that
we add up the reciprocals of the leading coefficients 
in the family of
reduced, primitive, inequivalent binary quadratic forms of integer
coefficients with discriminant $-D<0$.

Also, under the same condition, we have 
$$\sum_{(a,b,c)}{1\over a}=\prod_{\scriptstyle p|Z_0
\atop\scriptstyle \chi_{-D}(p)\ne -1}
{p+1\over p-1}\prod_{p|D}{p-1\over p}\ +\ {\rm negligible}.\eqno(0.17)$$
(0.16) and (0.17) explain why (0.15) {\it does}
distinguish between the primes
$p|Z_0$ with $\chi_{-D}(p)=1$, $0$ or $-1$.

In the rest of the paper we work out the details of this brief outline.
Sections 1-5 contain (almost) all basic ideas and ingredients
of the proof of Theorem 1 (case of negative discriminants).
Section 6 is an outline of the necessary changes in the case of
positive discriminants (see Theorem 2 in Section 1 below).
Part Two, including Sections 7-30, covers all the remaining details; they
 are mostly routine estimations that I could not simplify.

The basic idea is relatively simple and (almost) elementary, but
the execution of the simple plan required dozens of elementary estimations 
(like partial summations), which made the paper ridiculously long.
Unfortunately, I don't know how to cut the paper to a ``reasonable
size''. Neither the assumption of class number one, nor the Riemann
Hypothesis seem to make a difference here.

If a paper is more than 300 pages long, then it is almost
inevitable that there are mistakes. Are these mistakes substantial?
Is the basic idea (see Section 0, or Sections 1-5) still good?
I don't know the answer to these questions.

Part One includes Sections 0-6. The reader who is interested in 
Part Two should send me an email, and I will send him/her the pdf-file.
My email address is jbeck@math.rutgers.edu

\bigskip\bigskip\bigskip

\medskip
\centerline{\bf 1. Introduction }
\bigskip\bigskip\bigskip\medskip\noindent
We assume that the reader is 
familiar with both the elements of
multiplicative number theory and the elements of the theory
of binary quadratic forms (=quadratic fields); 
see e.g.  Davenport's well-known book 
 [Da] (see also the books [Bo-Sh] and  [Za]). Suppose for simplicity that
$\Delta $ is a fundamental discriminant
(i.e., either $\Delta\equiv 1$ (mod 4) squarefree, or
$\Delta/4\equiv 3$ (mod 4) squarefree, or
$\Delta/8$ odd squarefree integer; the fundamental discriminants
are precisely the discriminants of the quadratic fields).
 The class number $h(\Delta )$
(=the number of equivalence classes of binary quadratic forms of
discriminant $\Delta $) plays a key role in higher
arithmetic, and also in the distribution of the primes
in arithmetic progressions. The class number is strongly related to 
the well-known infinite series (real Dirichlet L-function at $s=1$)
$$L(1,\chi_{\Delta} )=\sum_{n=1}^{\infty }
{\chi_{\Delta }(n)\over n},\eqno(1.1)$$
where $\chi_{\Delta }$ is the  real primitive Dirichlet character
modulo $|\Delta |$, corresponding to the quadratic field of discriminant 
$\Delta$. (Note that $\chi_{\Delta }(n)=(\Delta/n)$ is also called the
Kronecker symbol: it is totally multiplicative in $n$,
$\chi_{\Delta }(p)=(\Delta/p)$ Legendre  symbol for every 
prime $p\ge 3$,  $\chi_{\Delta }(2)=1$ or $-1$ or 0 according to
$\Delta\equiv 1$ or 5 (mod 8) or $\Delta\equiv 0$  (mod 4), and 
$\chi_{\Delta }(-1)=1$ or $-1$ according to $\Delta >0$ or $<0$.)
The connection is Dirichlet's famous (analytic) class number formula,
$$L(1,\chi_{\Delta })={2\pi h(\Delta ) \over w\sqrt{|\Delta |}},\eqno(1.2)$$
for $\Delta <0$ (here $w=2$ if $\Delta <-4$, and $w=6
{\rm\ or\ }4$ if $\Delta 
=-3{\rm\ or\ }-4$), and 
$$L(1,\chi_{\Delta })={h(\Delta )\log \eta_{\Delta }
 \over \sqrt{\Delta }},\eqno(1.3)$$
for $\Delta >0$ (here $\eta_{\Delta }>1$ denotes the fundamental unit
and $\log $ is the natural logarithm).

Dirichlet also proved a remarkable finite class number formula:
for $\Delta <0$,
$$h(\Delta )=-{w/2\over |\Delta |}\sum_{n=1}^{
|\Delta |-1}n\chi_{\Delta }(n),\eqno(1.4)$$
and for $\Delta >0$,
$$h(\Delta )=-
{1\over \log \eta_{\Delta }}\sum_{n=1}^{\Delta -1}\chi_{\Delta }(n)
\log \sin (\pi n/\Delta ).\eqno(1.5)$$

To find a good lower bound for $L(1,\chi_{\Delta })$ is a famous open problem.
The best known result is Siegel's {\it ineffective} lower bound 
$$L(1,\chi_{\Delta })>{C_0(\eps )\over |\Delta |^{\eps }},\eqno(1.6)$$
which holds for any $\eps >0$ with some positive constant factor
$C_0(\eps )$ depending only on $\eps>0$.
 Unfortunately, there is no way to compute the constant
$C_0(\eps )$ whose existence is asserted in Siegel's proof (this
explains the term {\it ineffective}).

Combining (1.3) and (1.6), for $\Delta <0$ we have
$$h(\Delta )>C'_0(\eps )|\Delta |^{{1\over 2}-\eps },$$
where
$C'_0(\eps )$ is another ineffective constant
depending only on $\eps>0$.
The best known effective result, due to a combined effort of 
Goldfeld [Go2] and Gross--Zagier [Gr-Za] involving elliptic curves, 
is much weaker than (1.6):
it is basically logarithmic instead of the correct
roughly square root order of magnitude. As an illustration, we mention
the following explicit estimation, due to
 J. Oesterl\' e, who substantially improved on Goldfeld's
implicit constants: 
$$h(\Delta )>{\log |\Delta |\over 55}\prod_{p|\Delta 
}\left( 1-{\lfloor 2\sqrt{p}\rfloor
\over p+1}\right) .\eqno(1.7)$$
 (Goldfeld's method does not seem to work for
positive discriminants.)

The objective of this paper is to prove a new effective
lower bound to $L(1,\chi_{\Delta })$.
 First I discuss
the case of negative discriminant $\Delta <0$.
\medskip\noindent
{\bf Theorem 1} {\it For every negative
fundamental discriminant $\Delta<0$ with $|\Delta|>e^{10^{100}}$,
$$L(1,\chi_{\Delta})\ge {1\over (\log |\Delta|)^{146}}.\eqno(1.8)$$}
\medskip\noindent

The lower bound (1.8) also holds for positive fundamental 
discriminants $\Delta >0$, but the proof is somewhat more complicated
(this explains why we formulate the two cases in seperate theorems).
\medskip\noindent
{\bf Theorem 2} {\it For every positive
fundamental discriminant  $\Delta>e^{10^{100}}$
we have
$$L(1,\chi_{\Delta})\ge {1\over (\log \Delta)^{146}}.\eqno(1.9)$$}
\medskip\noindent
{\bf Remarks.} 
The exponent $146$ of $\log |\Delta|$ in (1.8) and (1.9) is
certainly far from the truth. In this paper I don't make a serious
effort to find the best exponent. Instead I focus on presenting the
idea of the (unfortunately very long) proof as clearly as possible.

The basic idea of the proof of Theorem 2 is the same as that of
Theorem 1, but there are some substantial differences in the details.
We explain the necessary modifications in  Sections 6 and 15.

\medskip\noindent
{\bf Proof of Theorem 1.}
Assume that there is a ``bad'' negative fundamental discriminant 
$\Delta =-D<0$ ($D$ is always positive in the whole paper)
such that
$$L(1,\chi_{-D})< {1\over (\log D)^{146}}\ \ {\rm and} \ \ 
D>e^{10^{100}}.\eqno(1.10)$$
Our goal is to derive a contradiction from (1.10).

By (1.2) and (1.10),
$$h(-D)<{\sqrt{D}\over (\log D)^{146}}.\eqno(1.11)$$

We need the following technical lemma.
\medskip\noindent
{\bf Lemma 1.1} {\it
Let $R_{\Delta }^{(0)}[N]$ denote the number of squarefree integers $n$
in $1\le n\le N$ which are represented by an arbitrary
binary quadratic form of
discriminant $\Delta <0$. 
For all $N>|\Delta |$,
$$\sum_{\scriptstyle  p\le N\ {\rm prime}\atop
\scriptstyle \chi_{\Delta }(p)=1}1\le R_{\Delta }^{(0)}[N]\le
{12h(\Delta )N\over\sqrt{|\Delta |}}.$$}
\medskip
{\bf Proof.}
Let $F_j(x,y)=a_jx^2+b_jxy+c_jy^2$, $1\le j\le h(\Delta )$ be the family of
reduced, primitive, inequivalent binary quadratic forms of integer
coefficients with
$$4a_jc_j=-\Delta +b_j^2=|\Delta |+b_j^2,\  a_j>0,\  c_j>0,$$
where reduced means
$$-a_j<b_j\le a_j\le c_j{\rm\ \ with\ \ }b_j\ge 0{\rm\ \ if\ \ }a_j=c_j.$$
These facts imply the useful inequalities $0<a\le\sqrt{|\Delta |/3}$
and $c_j\ge |\Delta |/(4a_j)$.

We rewrite $F_j$ as follows:
$$F_j(x,y)=a_jx^2+b_jxy+c_jy^2={(2a_jx+b_jy)^2+|\Delta |y^2\over 4a_j}.$$
This implies that
 $a_j$ is the first(=least) squarefree integer represented by
the binary form $F_j$, and $c_j$ is the second squarefree integer represented by
$F_j$. Let $N>c_j$; then we have
$$\sum_{\scriptstyle  p\le N\ {\rm prime}\atop
\scriptstyle F_j{\rm\ represents\ }p}1\le
\sum_{\scriptstyle  n\le N\ {\rm squarefree}\atop
\scriptstyle F_j{\rm\ represents\ }n}1\le
{1\over 2}\sum_{\scriptstyle (x,y)\in\Z^2:F_j(x,y)\le N\atop\scriptstyle
F_j(x,y)={\rm squarefree}}1\le $$
$$\le \sum_{1\le y\le\sqrt{4a_jN/|\Delta |}}
\sum_{\scriptstyle x\in\Z:\atop\scriptstyle
F_j(x,y)\le N}1\le \sum_{1\le y\le\sqrt{4a_jN/|\Delta |}}
2\left\lceil {\sqrt{4a_jN}\over 2a_j}\right\rceil \le $$
$$\le 2\left( 1+\sqrt{4a_jN\over |\Delta |}\right)\left( 1+
{\sqrt{4a_jN}\over 2a_j}\right) ={4N\over \sqrt{|\Delta |}}
+4\sqrt{a_jN\over |\Delta |}
+2\sqrt{N\over a_j}+2\le $$
$$\le {4N\over \sqrt{|\Delta |}}+4\sqrt{\sqrt{|\Delta |}N\over |\Delta |}
+2\sqrt{N}+2= $$
$$= {4N\over \sqrt{|\Delta |}}+4\sqrt{N\over\sqrt{|\Delta |}}
+2\sqrt{N}+2\le
{12N\over \sqrt{|\Delta |}},\eqno(1.12)$$
if $N>|\Delta |$ (where, as usual,
$\lceil z\rceil$ and $\lfloor z\rfloor$
denote the upper and lower integral parts of a real number $z$). 
Since $j$ runs in $1\le j\le h(\Delta )$, (1.12)
implies Lemma 1.1.\proofend\medskip

Lemma 1.1 implies that 
the  primes$\le N$
 are concentrated in the residue classes with
$\chi_{-D}(a)=-1$, explaining why I refer to them as the 
``rich'' residue classes
modulo $D$. The rest are the ``poor'' residue classes.
(Here we assume that $N$ is larger, but not much-much larger
than $D$; the precise definition comes soon.)
Since the ``poor'' half of the residue classes contain 
relatively few primes$\le N$, the ``rich'' half of 
residue classes contain {\it on average} 
twice as many primes$\le N$
as expected.  

(It is interesting to point out that, by the well-known 
Brun--Titchmarch--Selberg
theorem, no residue class can contain more than $2+o(1)$ times as many
primes as expected. So for the ``rich'' residue classes the average 
density is 
almost the same as the maximum density. Note, however, that
in this paper we don't use the Brun--Titchmarch--Selberg theorem.)

We often use the prime number theorem, 
which in its simplest form states
$$\pi(N)=\sum_{p\le N}1={N\over\log N}+O\left( N(\log N)^{-2}
\right) ,$$
but at some point of the proof
 we need the following much more precise result (see e.g. [Da] or
[Iw-Ko] or [Ka]):
$$\left|\pi(N)-\int_2^N\, {dx\over\log x}\right| =O(Ne^{-\sqrt{\log
N}}).$$
We will also apply (usually in a weaker form) 
the following well-known asymptotic
results related to the primes ($\gamma_0$ denotes Euler's constant):
$$\prod_{p\le N}\left( 1-{1\over p}\right)
=(1+o(1)){e^{-\gamma_0}\over\log N},$$
$$\prod_{p\le N}\left( 1+{1\over p}\right)=
(1+o(1)){6e^{\gamma_0}\log N\over \pi^2},$$
$$\prod_{p\le N}p=e^{(1+o(1))N}.$$

The proof of  Theorem 1 consists of several steps.
The first step is to construct a large number of almost orthogonal complex
 vectors in a high dimensional vector space. To take advantage of
almost orthogonality, we  apply
 a simple quadratic identity (see (1.17) below). 
This is a little bit similar to the basic idea of the Large Sieve,
{\it but} here we work with a single modulus: we restrict ourselves to the
reduced residue classes modulo $D$, where $-D$ is the ``bad'' discriminant.

We begin with a simple remark.
At a later stage of the proof, estimating some error term,
it would be helpful to have an estimation like
$$\sum_{p\le N:\chi_{-D}(p)=1}{1\over p}=o(1).$$
Unfortunately, this is not necessarily true in general, due 
to  the possible existence of very small primes 
$p$  with $\chi_{-D}(p)=1$.
Luckily, we can get around this problem: the  weaker statement 
$$\sum_{D<p\le N:\chi_{-D}(p)=1}{1\over p}=O(h(-D)\log
N/\sqrt{D})=o(1)$$
suffices for our purposes. Note that this weaker statement
is a consequence of  Lemma 1.1, and $o(1)$ clearly holds under the condition
(1.11), assuming
$N$ is not too large. This motivates why
we sieve out (among others)
the primes $p\le D$ with $\chi_{-D}(p)=1$ or 0.

To define the ``sieve out'' procedure precisely, we distinguish two cases.
Let $p_0$ denote the smallest prime that is not a divisor of
the discriminant $D$; clearly
$$p_0<2\log D.\eqno(1.13)$$
Indeed, (1.13) follows from the fact
$$\prod_{p\le N}p=e^{(1+o(1))N}$$
with the choice $N=\log D$.
If $\chi_{-D}(p_0)=1$ then let
$$Z_0=\prod_{\scriptstyle p\le D:\atop\scriptstyle
\chi_{-D}(p)\ne -1}p\eqno(1.14)$$
be the product of all primes $p\le D$
with $\chi_{-D}(p)=1$ or 0. Note that 
$p_0|Z_0$; moreover,
if an integer $m$ is relatively
prime to $Z_0$, then $m$ is automatically 
coprime to $D$ (since the prime factors of $D$ are listed
in $Z_0$).

If $\chi_{-D}(p_0)=-1$ then we have two possible choices
for $Z_0$: either (1.14), i.e.,
$$Z_0=Z'_0=\prod_{\scriptstyle p\le D:\atop\scriptstyle
\chi_{-D}(p)\ne -1}p$$
or we include $p_0$ as an extra prime factor:
$$Z_0=Z''_0=p_0Z'_0=
p_0\prod_{\scriptstyle p\le D:\atop\scriptstyle
\chi_{-D}(p)\ne -1}p.\eqno(1.15)$$
We discuss the two possible choices 
$Z_0=Z'_0$ and $Z_0=Z''_0$
in a unified way.

Note that
we are going to repeatedly use the following consequence of
 sieving out the primes
 $p\le D$ with $\chi_{-D}(p)=1$ or 0: if $m$ is coprime to $Z_0$
and squarefree,
then $\mu (m)=\chi_{-D}(m)$ for all $1\le m\le D$.

{\it Warning about the notation.}  It is important to remember that 
in the rest of the paper
$p,p_1,p_2,\ldots $ and $q,q_1,q_2,\ldots $ always denote primes,
and $p_0$ is reserved for the smallest prime that is not a divisor of
the discriminant $D$ (see (1.13)). 
Also
$c_1,c_2,\ldots $ denote effectively computable constants; in fact,
from now on every constant will be effective.  I use
$\gamma_0=0.5772\ldots $ for the Euler's constant
(instead of the more common $\gamma$), 
$\log $ denotes the natural logarithm (instead of $\ln$).
The rest of the notation is more or less standard; 
e.g., $\ph (n)$ denotes the
Euler's function, $\tau (n)$ is the number of divisors of $n$,
$\mu (n)$ is the M\" obius function, $A|B$ means that the integer
$B$ is divisible by the integer $A$, 
and $A\dagger B$ means that $B$ is not divisible by $A$.

Let $m\ge 1$ be a squarefree integer relatively prime to
$Z_0$, and let
 $\l $ be an integer in $1\le\l <m$ relatively
prime to $m$.
Given $0<A_1<A_2$, let
${\bf v}_{m,\l ;a}(A_1,A_2)$ denote the  vector
$${\bf v}_{m,\l ;a}(A_1,A_2)=$$
$$=(e^{2\pi {\rm i}(j_1D+a)\l /m},
e^{2\pi {\rm i}((j_1+1)D+a)\l /m},
e^{2\pi {\rm i}((j_1+2)D+a)\l /m},
\ldots , e^{2\pi {\rm i}(j_2D+a)\l /m})\eqno(1.16)$$
that plays a key role in the proof;
here $j=j_1$ is the smallest integer with $jD+a\ge A_1$, 
 $j=j_2$ is the largest integer with $jD+a\le A_2$, and
of course ${\rm i}=\sqrt{-1}$.
Note that the dimension of
${\bf v}_{m,\l ;a}(A_1,A_2)$ is 
$j_2-j_1+1=(A_2-A_1)/D+O(1)$.

We need the usual inner product
of complex vectors: if $\psi =(a_1,\ldots ,a_n)$ and
$\xi =(b_1,\ldots ,b_n)$ with $a_j,b_j$  complex numbers, then
$$\langle \psi ,\xi\rangle =\sum_{j=1}^na_j\overline{b_j}$$
where $\overline{b_j}$ is the complex conjugate of $b_j$.
We write, as usual, 
$$\| \psi\| =\sqrt{\langle\psi ,\psi\rangle },\ \ {\rm so\ \ }
\| \psi\|^2 =\langle\psi ,\psi\rangle $$
is the square of the norm. For arbitrary complex numbers
$C_j$ we have 
$$\left\|\sum_{j=1}^LC_j\psi_j\right\|^2=
{\rm Diagonal}+{\rm OffDiagonal},\eqno(1.17)$$
where
$${\rm Diagonal}=\sum_{j=1}^L|C_j|^2\|\psi_j\|^2 $$
and
$${\rm OffDiagonal}=\sum_{j=1}^L\sum_{\scriptstyle 1\le k\le L:\atop
\scriptstyle k\ne j}
C_j\overline{C_k}\langle\psi_j,\psi_k\rangle .$$

We apply the quadratic identity (1.17) with the choice
$$\psi_j={\bf v}_{m,\l ;a}(A_1,A_2){\rm\ \ and\ \ }
C_j=C_{m,\l ;a}={\mu (m)\over\ph (m)}\eqno(1.18)$$
for a fixed residue class $a$ modulo $D$.
(Note that the alternative choice
$$C_j=C_{m,\l ;a}={\mu (m)\over m}\eqno(1.18')$$
would be equally good.)

Later we are going to  add up these equalities
(i.e., (1.17) with the choice of (1.18)) 
for all $1\le a\le D$ with weight
$\chi_{-D}(a)$.
The reason behind combining the M\"obius
function with the real character $\chi_{-D}$ is the ``similarity''
between the two functions (restricted to square-free integers)
under the condition that $h(-D)$ is ``small''.
(I will return to this guiding intuition after (2.10) below.)
We will specify  the parameters $0<A_1<A_2$ later in Section 2;
 see (2.31).

\bigskip\bigskip\bigskip
\medskip
\centerline{\bf 2. More on our high-dimensional almost orthogonal
 vectors,}
\centerline{\bf and some routine smoothing }
\bigskip\bigskip\bigskip\medskip\noindent
Write
$${\bf V}_a(M;A_1,A_2)
=\sum_{\scriptstyle 1\le m\le M,1\le \l\le m:\atop
\scriptstyle \gcd (m,\l Z_0)=1}
{\mu (m)\over\ph (m)}{\bf v}_{m,\l ;a}(A_1,A_2).$$
Thus by (1.16),
$$\left\| {\bf V}_a(M;A_1,A_2)\right\|^2=
\sum_{j:\ A_1\le jD+a\le A_2}\left|
\sum_{\scriptstyle 1\le m\le M:
\atop\scriptstyle\gcd (Z_0, m)=1}{\mu (m)\over\ph (m)}
\sum_{\scriptstyle 1\le \l\le m:\atop
\scriptstyle \gcd (\l ,m)=1}e^{2\pi {\rm i}(a+jD)\l /m}\right|^2
.\eqno(2.1)$$
I recall the so-called Ramanujan's sum 
(see Theorem 272 in Hardy--Wright [Ha-Wr]):
$$\sum_{\scriptstyle 1\le \l<m\atop\scriptstyle
\gcd (\l ,m)=1}
e^{2\pi {\rm i}n\l /m}={\mu (m')\ph (m)\over\ph (m')}
,\eqno(2.2)$$
where $m'=m/\gcd (n,m)$. 
As usual,  $\gcd $ denotes the greatest common divisor; 
$$\ph (m)=m\prod_{p|m}\left( 1-{1\over p}\right) $$
denotes the  Euler's function;
$a|b$ means that $b$ is
divisible by $a$; and $\mu (n)$ stands for the M\"obius function:
$\mu (1)=1$, and for $n\ge 2$, $\mu (n)=(-1)^r$ if $n=p_1\cdots p_r$ with
$r$ distinct prime factors, and, finally,
 $\mu (n)=0$ if $n\ge 2$ is not squarefree.

By using (2.2),  we have
$$\sum_{\scriptstyle 1\le m\le M:\atop\scriptstyle
\gcd (Z_0, m)=1}{\mu (m)\over\ph (m)}
\sum_{\scriptstyle 1\le \l\le m:\atop
\scriptstyle \gcd (\l ,m)=1}e^{2\pi{\rm i}n\l /m}=$$
$$=\sum_{\scriptstyle 1\le m\le M:\atop
\scriptstyle\gcd (Z_0, m)=1}{\mu (m)\over\ph (m)}\cdot
{\mu (m/\gcd (m,n))\ph (m)\over\ph (m/\gcd (m,n))}=$$
$$=\sum_{\scriptstyle d\ge 1:\gcd (Z_0,d)=1\atop\scriptstyle
d|n}\mu (d)\ph (d)
\sum_{\scriptstyle 1\le m\le M:\gcd (Z_0, m)=1\atop
\scriptstyle \gcd (m,n)=d}
{\mu^2 (m)\over\ph (m)}.$$
Using this with $n=a+jD$ in (2.1), we have
$$\left\| {\bf V}_a(M; A_1,A_2)\right\|^2=$$
$$=\sum_{j:\ A_1\le jD+a\le A_2}
\left(\sum_{\scriptstyle d\ge 1:\gcd (Z_0,d)=1\atop\scriptstyle
d|a+jD}\mu (d)\ph (d)
\sum_{\scriptstyle 1\le m\le M:\gcd (Z_0,m)=1\atop
\scriptstyle \gcd (m,a+jD)=d}
{|\mu (m)|\over\ph (m)}\right)^2=$$
$$=\sum_{j:\ A_1\le jD+a\le A_2}
\left(\sum_{\scriptstyle d\ge 1:\gcd (Z_0,d)=1\atop\scriptstyle
d|a+jD}\mu (d)
\sum_{\scriptstyle 1\le k\le M/d:\atop
\scriptstyle \gcd (k,(a+jD)Z_0)=1}
{|\mu (k)|\over\ph (k)}\right)^2.\eqno(2.3)$$
{\bf Remark.} At first sight the big sum (2.3)
may seem hopelessly complicated. We will show, however, that 
(2.3) is not so bad. For illustration note that
the evaluation of the square
$$\left(
\sum_{\scriptstyle d\ge 1:\gcd (Z_0,d)=1\atop\scriptstyle
d|a+jD}\mu (d)
\sum_{\scriptstyle 1\le k\le M/d:\atop
\scriptstyle \gcd (k,(a+jD)Z_0)=1}
{|\mu (k)|\over\ph (k)}\right)^2$$
at the end of (2.3) is quite simple
if $a+jD$ is a large prime:
$a+jD=p>M$. Indeed, then we have
$$\sum_{\scriptstyle d\ge 1:\gcd (Z_0,d)=1\atop\scriptstyle
d|p}\mu (d)
\sum_{\scriptstyle 1\le k\le M/d:\atop
\scriptstyle \gcd (k,pZ_0)=1}
{|\mu (k)|\over\ph (k)}=\sum_{\scriptstyle 1\le m\le M:
\atop\scriptstyle\gcd (m,Z_0)=1}
{|\mu (m)|\over\ph (m)}.$$
We use the fact
$$\sum_{\scriptstyle 1\le m\le M:\atop
\scriptstyle\gcd (m,Z_0)=1}
{|\mu (m)|\over\ph (m)}=\prod_{p|Z_0}\left( 1-{1\over p}\right)
(\log M+c)+{\rm\ negligible},\eqno(2.4)$$
where $c$ is some explicit absolute constant.
Note that (2.4) is the special case $Q=1$ of the following 
 lemma (in later applications
 we need the most general form of the lemma).
\medskip\noindent
{\bf Lemma 2.1} {\it Assume that $Z_0$ and
$Q\ge 1$ are relatively prime integers, then
$$\left|\sum_{\scriptstyle 1\le k\le L:\atop
\scriptstyle \gcd (k,QZ_0)=1}
{\mu^2 (k)\over\ph (k)}-\prod_{q|QZ_0}\left( 1-{1\over q}\right)\left(
\log L+\gamma_0+2\gamma^{\star}
-\gamma^{\star\star}+\sum_{q|QZ_0}{\log q\over q}
\right)\right|\le $$
$$\le {10^4\tau (Q)\log D\log L\over L^{1/4}}+
{10^4(10+\log Q)\over D^5}+
{4\left( 10+\log L+2(\log D)^2+(\log Q)^2\right) \over
\max\left\{ LD^{-6\log D},1\right\} } ,$$
where 
$\gamma_0$ is the  Euler's constant, 
$\gamma^{\star}$, $\gamma^{\star\star}$ are two absolute  constants
 defined by the following convergent prime-series
$$\gamma^{\star}=\sum_p{\log p\over p^2-1}\ \ {\rm and}\ \ 
\gamma^{\star\star}=\sum_p{\log p\over p(p+1)},$$
and $\tau (Q)=\sum_{d|Q}1$ is the number of divisors of $Q$. }
\medskip
Probably the reader is wondering why  we use
the at first sight artificially complicated requirement
$\gcd (k,QZ_0)=1$ in Lemma 2.1 instead of simply writing 
$\gcd (k,Q)=1$ with $Z_0|Q$. 
The reason is that in the error term of Lemma 2.1
we have the factor $\tau (Q)$ (and it is necessary 
to have it there), and
$$Z_0=Z'_0=\prod_{\scriptstyle p\le D:\atop\scriptstyle
\chi_{-D}(p)\ne -1}p{\rm\ \ or\ \ }
Z_0=Z''_0=p_0Z'_0$$
can be exponentially large in terms of  $D$. 
So if $Z_0$ is a divisor of $Q$, then
$\tau (Q)\ge \tau (Z_0)$, where
 $\tau (Z_0)$  is possibly exponentially large.
Such an extremely large factor $\tau (Q)$ in the error term would
make Lemma 2.1 nearly useless.

We postpone the elementary 
 proof of Lemma 2.1 to Section 30.

Let's return to (2.1): 
an alternative way to evaluate it is to apply (1.17):
$$\left\| {\bf V}_a(M;A_1,A_2)\right\|^2={\rm Diagonal}_a(M;A_1,A_2)+{\rm
OffDiagonal}_a(M;A_1,A_2).\eqno(2.5)$$
Note that
$$\langle {\bf v}_{m_1,\l_1;a}(A_1,A_2),{\bf v}_{m_2,\l_2;a}(A_1,A_2)
\rangle=\sum_{\scriptstyle A_1\le n\le A_2:\atop\scriptstyle n\equiv a\ 
({\rm mod}\ D)}
e^{2\pi{\rm i}n({\l_1\over m_1}-{\l_2\over m_2})}= $$
$$=e^{2\pi{\rm i}a({\l_1\over m_1}-{\l_2\over m_2})}
\sum_{j:\ A_1\le jD+a\le A_2} 
e^{2\pi{\rm i}Dj({\l_1\over m_1}-{\l_2\over m_2})} ,\eqno(2.6) $$
and 
$$\langle {\bf v}_{m,\l;a}(A_1,A_2),{\bf v}_{m,\l;a}(A_1,A_2)\rangle =\|
{\bf v}_{m,\l;a}\|^2=
\sum_{\scriptstyle A_1\le n\le A_2:\atop \scriptstyle n\equiv a\ ({\rm
mod\ }D)}1 .\eqno(2.7)$$
By (2.7) we have
$${\rm Diagonal}_a(M;A_1,A_2)=D_a(M;A_1,A_2)=$$
$$=\sum_{\scriptstyle 1\le m\le M,1\le \l\le m:\atop
\scriptstyle \gcd (\l Z_0,m)=1}\left({\mu (m)\over\ph (m)}\right)^2
\|{\bf v}_{m,\l ;a}\|^2=$$
$$=\sum_{\scriptstyle 1\le m\le M:\atop
\scriptstyle\gcd (Z_0,m)=1}\ph (m)
\left({\mu (m)\over\ph (m)}\right)^2
\|{\bf v}_{m,\l ;a}\|^2=$$
$$=\left(\sum_{\scriptstyle A_1\le n\le A_2:\atop
\scriptstyle n\equiv a\ ({\rm mod\ }D)}1\right)
\sum_{\scriptstyle 1\le m\le M:\atop
\scriptstyle\gcd (Z_0,m)=1}{|\mu (m)|\over\ph (m)}
,\eqno(2.8)$$
and by (2.6),
$${\rm OffDiagonal}_a(M;A_1,A_2)=OD_a(M;A_1,A_2)=
\sum_{\scriptstyle d\ge 1:\atop\scriptstyle \gcd (d,Z_0)=1}
OD_a(M;A_1,A_2;d),\eqno(2.9)$$
where
$$OD_a(M;A_1,A_2;d)=\sum_{{\scriptstyle (m_1,\l_1)\ne 
(m_2,\l_2):\ 1\le m_1,m_2\le M\atop
\scriptstyle 1\le\l_h\le m_h,\gcd (\l_h,m_h)=1,h=1,2}\atop\scriptstyle
\gcd (m_1,m_2)=d,\gcd (m_1m_2,Z_0)=1}{\mu(m_1)\mu(m_2)
\over\ph (m_1)\ph (m_2)}\cdot \qquad\qquad $$
$$\qquad\qquad\qquad\cdot
e^{2\pi{\rm i}a({\l_1\over m_1}-{\l_2\over m_2})}
\cdot \sum_{j:\ A_1\le jD+a\le A_2} 
e^{2\pi{\rm i}Dj({\l_1\over m_1}-{\l_2\over m_2})} .\eqno(2.10)$$

The proof of Theorem 1 is very long, but the basic idea
is relatively simple, and it is explained in the first 5 sections.

Here is a
brief intuitive explanation of  why we combine the M\"obius function
$\mu (m)$
with exponential sums
related to the reduced residue classes modulo $D$.
The exponential sum related to the fixed
modulus $D$ leads to a Gauss sum (see
Section 3; especially Lemma 3.1), 
and the Gauss sum  contains the real character $\chi_{-D}(m)$ as a factor.
If $m$ is a squarefree
integer, relatively prime to $D$, then $\chi_{-D}(m)$ and $\mu (m)$
are ``typically equal'', assuming $m$ has no ``small'' prime divisor
$p$ with $\chi_{-D}(p)=1$ and $h(-D)$ is ``substantially'' smaller than
$\sqrt{D}$. This ``almost equality'' of
$\chi_{-D}(m)$ and $\mu (m)$ makes it possible to find a good estimation
for the off-diagonal part (2.9)-(2.10), which is the hardest part of
the proof of Theorem 1.
I refer to this argument as the ``almost equality of
$\chi_{-D}(m)$ and $\mu (m)$''.

Sections 2 and 3 are routine but important preparatory sections.
Sections 4 and 5 are the two most important sections; Lemma 5.4
and Lemma 5.5 are the two crucial lemmas. The rest are basically 
long but routine estimations.

\medskip

The parameters $A_1, A_2$ 
indicate that the underlying interval is $A_1\le n\le A_2$; more precisely,
at the end of (2.10) we have the following exponential sum
$$\sum_{j:\ A_1\le jD+a\le A_2} 
e^{2\pi{\rm i}Dj({\l_1\over m_1}-{\l_2\over m_2})} .\eqno(2.11)$$
It is very useful to take a certain average here (where $A_1$ and $A_2$
are variables):
 we employ the standard trick of ``smoothing''
in Fourier analysis. The details go as follows.

I start with the well-known Dirichlet kernel
$$\sum_{k=-n}^ne^{{\rm i}kx}={\sin (n+{1\over 2})x\over \sin (x/2)}.
\eqno (2.12)$$
Squaring (2.12), we obtain the Fej\' er kernel
$$\left({\sin (n+{1\over 2})x\over \sin (x/2)}\right)^2=
\left(\sum_{k=-n}^ne^{{\rm i}kx}\right)^2=
\sum_{\l =-2n}^{2n}(2n+1-|\l |)e^{{\rm i}\l x}=$$
$$=S_{2n}(x)+S_{2n-1}(x)+S_{2n-2}(x)+\ldots +S_1(x)+S_0(x),\eqno (2.13)$$
where
$$S_m(x)=\sum_{k=-m}^me^{{\rm i}kx}.\eqno(2.14)$$
The sequence of coefficients $2n+1-|\l |$ of $e^{{\rm i}\l x}$ in (2.13)
represents the ``roof function''
$$f_2(y)=\cases{{1\over 2}-{|y|\over 4},&if $|y|\le 2$;\cr
0,&if $|y|>2$;\cr}\eqno(2.15)$$
in the sense that
$${2n+1-|\l |\over (2n+1)^2}\cdot n={1\over 2}-{|y|\over 4}+O(1/n)
{\rm\ \ with\ \ }y=\l/n\eqno(2.16)$$
(here the factor $n$ is explained by the renorming $y=\l/n$).
Note that $f_2(y)$ is familiar from probability theory:
it is the density function of the convolution---denoted by 
$\ast$---of the uniform
distribution in $[-1,1]$ with itself; formally,
$$f_2(y)=f_1\ast f_1(y)=\int_{-\infty }^{\infty}
f_1(y-z)f_1(z)\, dz,$$
where
$$f_1(y)=\cases{{1\over 2},&if $|y|\le 1$;\cr
0,&if $|y|>1$;\cr}\eqno(2.17)$$
is the density function of the uniform
distribution in $[-1,1]$.

We rewrite (2.13) in the form
$$\left({\sin (n+{1\over 2})x\over (2n+1)\sin (x/2)}\right)^2=
\sum_{k=0}^{2n}{S_k(x)\over 2k+1}\cdot {2k+1\over (2n+1)^2}
\eqno(2.18)$$
with $\sum_{k=0}^{2n}{2k+1\over (2n+1)^2}=1$, i.e., the weights
${2k+1\over (2n+1)^2}$, $0\le k\le 2n$ in (2.18) represent
 a discrete probability distribution.

It is a standard exercise in probability theory to compute
the higher convolution powers of the uniform
distribution in $[-1,1]$. A routine calculation gives
$$f_3(y)=f_1\ast f_1\ast f_1(y)=
\cases{{3-y^2\over 8},&if $|y|\le 1$;\cr
{(3-|y|)^2\over 16},&if $1\le |y|\le 3$;\cr
0,&if $|y|>3$.\cr }\eqno(2.19)$$
Consider now the third power of (2.12):
$$\left({\sin (n+{1\over 2})x\over \sin (x/2)}\right)^3=
\left(\sum_{k=-n}^ne^{{\rm i}kx}\right)^3=
\sum_{\l =-3n}^{3n}B^{(3)}_{\l } e^{{\rm i}\l x}.
\eqno (2.20)$$
It is not too difficult to compute the coefficients
$B^{(3)}_{\l }$ of $ e^{{\rm i}\l x}$ explicitly, but we don't really
need that. What we need is the analog of
(2.15)-(2.16):
$${B^{(3)}_{\l }\over (2n+1)^3}\cdot n=f_3(\l/n)+O(1/n),\eqno(2.21)$$
where $f_3(y)$ is the function defined in (2.19). 
We also need the easy fact that the sequence
$B^{(3)}_{\l }$ is monotonic in the following sense:
$$B^{(3)}_{\l_1 }\ge B^{(3)}_{\l_2 }\ge 0
{\rm\ \ if \ \ }|\l_1|\le |\l_2|.\eqno(2.22)$$
Therefore, returning to (2.20), we obtain the following analog of (2.18): 
$$\left({\sin (n+{1\over 2})x\over (2n+1) \sin (x/2)}\right)^3=
\left({1\over 2n+1}\sum_{k=-n}^ne^{{\rm i}kx}\right)^3=$$
$$=(2n+1)^{-3}\sum_{\l =-3n}^{3n}B^{(3)}_{\l } e^{{\rm i}\l x}=
\sum_{k=0}^{3n}w^{(3,n)}_k{S_k(x)\over 2k+1},
\eqno (2.23)$$
where $w^{(3,n)}_k$ are positive weights:
$w^{(3,n)}_k\ge 0$, and $\sum_{k=0}^{3n}w^{(3,n)}_k=1$, i.e., 
the coefficients $w^{(3,n)}_k$ form a
discrete probability distribution.
So (2.23) is an {\it average} of 
${S_k(x)\over 2k+1}$, $0\le k\le 3n$
(we divide by $2k+1$ since the sum $S_k(x)$ has $2k+1$ terms; see (2.14)).

In general, consider the $\kappa$th power of (2.12) where $\kappa\ge 2$ is an
arbitrary integer:
$$\left({\sin (n+{1\over 2})x\over \sin (x/2)}\right)^{\kappa}=
\left(\sum_{k=-n}^ne^{{\rm i}kx}\right)^{\kappa}
=\sum_{\l =-\kappa n}^{\kappa n}B^{(\kappa)}_{\l } e^{{\rm i}\l x}.
\eqno (2.24)$$
For a general $\kappa \ge 2$ it is not easy to compute the coefficients
$B^{(\kappa)}_{\l }$ of $ e^{{\rm i}\l x}$ explicitly, but we don't really
need that. We need less: we just need the analog of
(2.21):
$${B^{(\kappa)}_{\l }\over (2n+1)^{\kappa}}\cdot n=f_{\kappa}
(\l/n)+O(1/n),\eqno(2.25)$$
where
$$f_{\kappa}(y)=f_1\ast f_1\ast\cdots \ast f_1 (y)\eqno(2.26)$$
is the $\kappa$th convolution power of $f_1$ (see (2.17)).
It is not very hard to prove by induction on $\kappa$
the following generalization of (2.15) and (2.19) (see e.g.
R\' enyi's book [Re]):
$$f_{\kappa}(y)={1\over 2^{\kappa}(\kappa-1)!}
\sum_{j=0}^{\lfloor (\kappa+|y|)/2\rfloor }
(-1)^j{\kappa\choose j}(\kappa+|y|-2j)^{\kappa-1}{\rm\ \ if\ \ }
|y|\le \kappa,\eqno(2.27)$$
and 0 if $|y|>\kappa$. Note that $f_{\kappa}(y)$ is $(\kappa-2)$-times 
differentiable,
and consists of a few generalized parabola arcs of degree $\kappa-2$
(due to the jumps of the lower integral part function
$\lfloor (\kappa+y)/2\rfloor $ as $y$ runs in $0\le y\le \kappa$)
that smoothly fit together at the endpoints if $\kappa\ge 3$.

We also need the easy fact that the sequence
$B^{(\kappa)}_{\l }$ is monotonic in the following sense:
$$B^{(\kappa)}_{\l_1 }\ge B^{(\kappa)}_{\l_2 }\ge 0
{\rm\ \ if \ \ }|\l_1|\le |\l_2|.\eqno(2.28)$$
Therefore, returning to (2.24) we have
$$\left({\sin (n+{1\over 2})x\over (2n+1)\sin (x/2)}\right)^{\kappa}=
\left({1\over 2n+1}\sum_{k=-n}^ne^{{\rm i}kx}\right)^{\kappa }=$$
$$=(2n+1)^{-\kappa}\sum_{\l =-\kappa n}^{\kappa n}
B^{(\kappa)}_{\l } e^{{\rm i}\l x}=
\sum_{k=0}^{\kappa n}w^{(\kappa ,n)}_k{S_k\over 2k+1}(x),
\eqno (2.29)$$
where $w^{(\kappa ,n)}_k$ are positive weights:
$w^{(\kappa ,n)}_k\ge 0$ with $\sum_{k=0}^{\kappa n}w^{(\kappa ,n)}_k=1$ 
(``discrete probability distribution'').
So (2.29) is an average of 
${S_k(x)\over 2k+1}$, $0\le k\le \kappa n$.

Note that in (2.29) there is an underlying limit
as $\kappa$ increases. Indeed,
due to the Central Limit Theorem,
$$f_{\kappa}(y)={1\over \sqrt{2\pi }\cdot \sqrt{\kappa/3}}
e^{-{3y^2\over 2\kappa}}+O(1/\kappa).\eqno(2.30)$$
Here the factor $\sqrt{\kappa/3}$ is the ``standard deviation''. Indeed,
in view of (2.26), what we are dealing with
 here is a sum of $\kappa$ independent random variables, where each
component is uniformly distributed in $[-1,1]$. This is why the variance is
$\kappa\int_0^1x^2\, dx=\kappa/3$,
and this is why (2.30) is a special case of
the strong form of the Central Limit Theorem  with explicit error term
estimation.

Now we are ready to define an efficient {\it average} over the
variables $A_1,A_2$ introduced
in Section 1. Let $\kappa $ and $N$ be fixed
positive integers (in the application later we choose $\kappa $ to be
 less than 10---note in advance that
 $\kappa =8$ is a good choice---and $N$ is 
``large''), and let $k$ be an integral 
variable running in $0\le k\le\kappa N$. Write
$$A_1=A_1(a;k)=((\kappa N-k)+a)D {\rm\ \ and\ \ }
A_2=A_2(a;k)=((\kappa N+k)+a)D .\eqno(2.31)$$
By (2.31)
$$\{ j:\ A_1(a;k)\le jD+a\le A_2(a;k)\} =\{ j:\ 
\kappa N-k\le j\le \kappa N+k\} \eqno(2.32)$$
is a set of $2k+1$ consecutive integers (this explains the division by
$2k+1$ below).
Write
$$W_a(M;\kappa ;N)=\sum_{k=0}^{\kappa N}w^{(\kappa ,N)}_k
{\left\| {\bf V}_a(M;A_1(a;k),A_2(a;k))\right\|^2
\over 2k+1}=$$
$$=\sum_{k=0}^{\kappa N}w^{(\kappa ,N )}_k
{1\over 2k+1}\sum_{j:\ \kappa N-k\le j\le \kappa N+k}
\left(
\sum_{\scriptstyle d\ge 1:\gcd (Z_0,d)=1\atop\scriptstyle
d|a+jD}\mu (d)
\sum_{\scriptstyle 1\le k\le M/d:\atop
\scriptstyle \gcd (k,(a+jD)Z_0)=1}
{|\mu (k)|\over\ph (k)}\right)^2
,\eqno(2.33)$$
where $w^{(\kappa ,N)}_k$, $0\le k\le \kappa N$ is the discrete
probability distribution defined in (2.29),
and in the last step we used (2.3) and (2.32).

By (2.5)-(2.10) and (2.31)-(2.33), we have
$$W_a(M;\kappa ;N)=
{\rm WDiag}_a(M;\kappa ;N)+{\rm
WOffDiag}_a(M;\kappa ;N),\eqno(2.34)$$
where
$${\rm WDiag}_a(M;\kappa ;N)=
\sum_{k=0}^{\kappa N}w^{(\kappa ,N)}_k
{{\rm Diagonal}_a(M;A_1(a;k),A_2(a;k))\over 2k+1}=$$
$$=\sum_{k=0}^{\kappa N}w^{(\kappa ,N)}_k\cdot {1\over 2k+1}
\left(\sum_{\scriptstyle A_1(a;k)\le n\le A_2(a;k):\atop\scriptstyle
n\equiv a({\rm mod\ }D)}1\right)
\sum_{\scriptstyle 1\le m\le M:\atop
\scriptstyle\gcd (Z_0, m)=1}{|\mu (m)|\over\ph (m)}=$$
$$=\left(\sum_{k=0}^{\kappa N}w^{(\kappa ,N)}_k\right)
\sum_{\scriptstyle 1\le m\le M:\atop
\scriptstyle\gcd (Z_0,m)=1}{|\mu (m)|\over\ph (m)}=
\sum_{\scriptstyle 1\le m\le M:\atop
\scriptstyle\gcd (Z_0,m)=1}{|\mu (m)|\over\ph (m)}, \eqno(2.35)$$
and
$${\rm WOffDiag}_a(M;\kappa ;N)=
\sum_{k=0}^{\kappa N}w^{(\kappa ,N)}_k
{{\rm OffDiagonal}_a(M;A_1(a;k),A_2(a;k))\over 2k+1}=$$
$$=\sum_{\scriptstyle d\ge 1:\atop\scriptstyle \gcd (d,Z_0)=1}
WOD_a(M;\kappa ;N;d),\eqno(2.36)$$
where
$$WOD_a(M;\kappa ;N;d)
=\sum_{{\scriptstyle (m_1,\l_1)\ne 
(m_2,\l_2):\ 1\le m_1,m_2\le M\atop
\scriptstyle 1\le\l_h\le m_h,\gcd (\l_h,m_h)=1,h=1,2}\atop\scriptstyle
\gcd (m_1,m_2)=d,\gcd (m_1m_2,Z_0)=1}{\mu(m_1)\mu(m_2)
\over\ph (m_1)\ph (m_2)}\cdot \qquad\qquad $$
$$\qquad\cdot
e^{2\pi{\rm i}a({\l_1\over m_1}-{\l_2\over m_2})}
\sum_{k=0}^{\kappa N}w^{(\kappa ,N)}_k
\cdot {1\over 2k+1}
\sum_{j:\ A_1(a;k)\le jD+a\le A_2(a;k)} 
e^{2\pi{\rm i}Dj({\l_1\over m_1}-{\l_2\over m_2})} .\eqno(2.37)$$
By (2.32),
$$\sum_{j:\ A_1(a;k)\le jD+a\le A_2(a;k)} 
e^{2\pi{\rm i}Dj({\l_1\over m_1}-{\l_2\over m_2})}
=e^{{\rm i}\kappa Nx}S_k(x){\rm\ \ with\ \ }x=2\pi 
D({\l_1\over m_1}-{\l_2\over m_2}),
\eqno(2.38)$$
where $S_k(x)$ is the exponential sum defined in (2.14).

Therefore, by using (2.29) with $x=2\pi D({\l_1\over m_1}-{\l_2\over
m_2})$, (2.37) and (2.38) imply
$$WOD_a(M;\kappa ;N;d)
=\sum_{{\scriptstyle (m_1,\l_1)\ne 
(m_2,\l_2):\ 1\le m_1,m_2\le M\atop
\scriptstyle 1\le\l_h\le m_h,\gcd (\l_h,m_h)=1,h=1,2}\atop\scriptstyle
\gcd (m_1,m_2)=d,\gcd (m_1m_2,Z_0)=1}{\mu(m_1)\mu(m_2)
\over\ph (m_1)\ph (m_2)}\cdot \qquad\qquad $$
$$\quad\cdot
e^{2\pi{\rm i}a({\l_1\over m_1}-{\l_2\over m_2})}
\cdot
e^{2\pi{\rm i}D({\l_1\over m_1}-{\l_2\over m_2})\kappa N}
\left( {\sin\left( (2N+1)\pi D({\l_1\over m_1}-{\l_2\over
m_2})\right)\over (2N+1)\sin\left( \pi D({\l_1\over m_1}-{\l_2\over
m_2})\right)}\right)^{\kappa }.\eqno(2.39)$$

With $\delta =\pm 1$, write  
$$\Omega_{\delta ;\kappa ;N}(M)=
\sum_{\scriptstyle 1\le a\le D\atop\scriptstyle \chi_{-D}(a)=\delta}
W_a(M;\kappa ;N)=$$
$$=\sum_{k=0}^{\kappa N}w^{(\kappa ,N)}_k
{1\over 2k+1}\sum_{j:\ \kappa N-k\le j\le \kappa N+k}
\sum_{\scriptstyle 1\le a\le D\atop\scriptstyle \chi_{-D}(a)=\delta}$$
$$\cdot\left(
\sum_{\scriptstyle d\ge 1:\gcd (Z_0,d)=1\atop\scriptstyle
d|a+jD}\mu (d)
\sum_{\scriptstyle 1\le k\le M/d:\atop
\scriptstyle \gcd (k,(a+jD)Z_0)=1}
{|\mu (k)|\over\ph (k)}\right)^2,\eqno(2.40)$$
where in the last step we used (2.33).
Similarly, write (see (2.31)-(2.37))
$$\Omega_{\delta ;\kappa ;N}({\rm Diag};M)=
\sum_{\scriptstyle 1\le a\le D\atop\scriptstyle \chi_{-D}(a)=\delta}
{\rm WDiag}_a(M;\kappa ;N)={\ph (D)\over 2}
\sum_{\scriptstyle 1\le m\le M:\atop
\scriptstyle\gcd (Z_0,m)=1}{|\mu (m)|\over\ph (m)}
,\eqno(2.41)$$
and
$$\Omega_{1 ;\kappa ;N}({\rm OffDiag};M)-
\Omega_{-1 ;\kappa ;N}({\rm OffDiag};M)=$$
$$=\sum_{\scriptstyle 1\le a\le D\atop\scriptstyle \chi_{-D}(a)=1}
{\rm WOffDiag}_a(M;\kappa ;N)-
\sum_{\scriptstyle 1\le a\le D\atop\scriptstyle \chi_{-D}(a)=-1}
{\rm WOffDiag}_a(M;\kappa ;N)=$$
$$=\sum_{1\le a\le D}\chi_{-D}(a)
{\rm WOffDiag}_a(M;\kappa ;N)=$$
$$=\sum_{{\scriptstyle (m_1,\l_1)\ne 
(m_2,\l_2):\ 1\le m_1,m_2\le M\atop
\scriptstyle 1\le\l_h\le m_h,\gcd (\l_h,m_h)=1,h=1,2}\atop\scriptstyle
\gcd (m_1,m_2)=d,\gcd (m_1m_2,Z_0)=1}{\mu(m_1)\mu(m_2)
\over\ph (m_1)\ph (m_2)}\cdot \qquad\qquad $$
$$\cdot\left(\sum_{1\le a\le D} \chi_{-D}(a)
e^{2\pi{\rm i}a({\l_1\over m_1}-{\l_2\over m_2})}\right)
e^{2\pi{\rm i}D({\l_1\over m_1}-{\l_2\over m_2})\kappa N}
\left( {\sin\left( (2N+1)\pi D({\l_1\over m_1}-{\l_2\over
m_2})\right)\over (2N+1)\sin\left( \pi D({\l_1\over m_1}-{\l_2\over
m_2})\right)}\right)^{\kappa }
.\eqno(2.42)$$
By (2.34),
$$\Omega_{\delta ;\kappa ;N}(M)=\Omega_{\delta ;\kappa ;N}({\rm Diag};M)
+\Omega_{\delta ;\kappa ;N}({\rm OffDiag};M).\eqno(2.43)$$

Finally, for technical reasons to be explained later, we assume that
$N$ is also a
variable, and we take the average as $N$ runs in $T\le N<2T$: let
($\delta =1$ or $-1$)
$$\overline{\Omega}_{\delta ;\kappa ;T}(M)={1\over T}\sum_{T\le N<2T}
\Omega_{\delta ;\kappa ;N}(M)
.\eqno(2.44)$$
Similarly, by (2.41),
$$\overline{\Omega}_{\delta ;\kappa ;T}({\rm Diag};M)=
{1\over T}\sum_{T\le N<2T}
\Omega_{\delta ;\kappa ;N}({\rm Diag};M)
={\ph (D)\over 2}\sum_{\scriptstyle 
1\le m\le M:\atop\scriptstyle\gcd (Z_0,m)=1}{|\mu (m)|\over\ph (m)}
,\eqno(2.45)$$
and
$$\overline{\Omega}_{\delta ;\kappa ;T}({\rm OffDiag};M)=
{1\over T}\sum_{T\le N<2T}
\Omega_{\delta ;\kappa ;N}({\rm OffDiag};M).\eqno(2.46)$$
By (2.43),
$$\overline{\Omega}_{\delta ;\kappa ;T}(M)=
\overline{\Omega}_{\delta ;\kappa ;T}({\rm Diag};M)+
\overline{\Omega}_{\delta ;\kappa ;T}({\rm OffDiag};M).\eqno(2.47)$$

Next we specify the key parameters $M$, $\kappa$,
and $T$: let 
$$\kappa =8,\ \ \log T=(\log D)^{15},\ \ M=T\exp\left(-{2\over 3}
\sqrt{\log T}\right) ,\eqno(2.48)$$
noting that $N$ is an integer variable running in the interval
$T\le N<2T$.

We conclude Section 2 with a 
\medskip\noindent
{\bf Brief summary of the proof of Theorem
1.} By (2.47), 
$$\overline{\Omega}_{1 ;\kappa ;T}(M)
-\overline{\Omega}_{-1 ;\kappa ;T}(M)=$$
$$=\overline{\Omega}_{1 ;\kappa ;T}({\rm OffDiag};M)-
\overline{\Omega}_{-1 ;\kappa ;T}({\rm OffDiag};M)+
\overline{\Omega}_{1 ;\kappa ;T}({\rm Diag};M)
-\overline{\Omega}_{-1 ;\kappa ;T}({\rm Diag};M)=$$
$$=\overline{\Omega}_{1 ;\kappa ;T}({\rm OffDiag};M)-
\overline{\Omega}_{-1 ;\kappa ;T}({\rm OffDiag};M),
\eqno(2.49)$$
i.e., the diagonal part cancels out.
It is relatively easy to 
evaluate/estimate $\overline{\Omega}_{\pm 1 ;\kappa ;T}(M)$
by using (2.40), (2.44) and Lemma 2.1; see Lemma 5.4 later.
Thus we obtain the following approximation of the left-hand side of
(2.49):
$$\overline{\Omega}_{1 ;\kappa ;T}(M)-
\overline{\Omega}_{-1 ;\kappa ;T}(M)=$$
$$=\prod_{p|Z_0}\left( 1-{1\over p}\right)^2\cdot
D\log {M^2\over T}
+{\rm negligible\ error\ term},\eqno(2.50)$$
where the ``negligible error term'' will be specified in Section 5.

On the other hand, by using (2.42) and (2.46), after a long chain of
estimations we obtain 
the following approximation of the right-hand side of
(2.49) (see Lemma 5.5 later):
$$\overline{\Omega}_{1 ;\kappa ;T}({\rm
OffDiag};M)-\overline{\Omega}_{-1 ;\kappa ;T}({\rm
OffDiag};M)=$$
$$=\prod_{\scriptstyle p|Z_0\atop\scriptstyle \chi_{-D}(p)=1}
\left( 1-{1\over p}\right)
\prod_{\scriptstyle p|Z_0\atop\scriptstyle \chi_{-D}(p)=0}
\left( 1-{1\over p}\right)^2
\prod_{\scriptstyle p|Z_0\atop\scriptstyle \chi_{-D}(p)=-1}
\left( 1-{1\over p}\right)^2\left( 1-{1\over p+1}\right)
\cdot $$
$$\qquad\qquad\qquad\qquad\cdot D\log {M^2\over T}
+{\rm negligible\ error\ term},\eqno(2.51)$$
and again the ``negligible error term'' will be specified in Section 5.
(Note that
$$\prod_{\scriptstyle p|Z_0\atop\scriptstyle \chi_{-D}(p)=-1}
\left( 1-{1\over p}\right)^2\left( 1-{1\over p+1}\right)
{\rm\ is\ either\ }
\left( 1-{1\over p_0}\right)^2\left( 1-{1\over p_0+1}\right)
{\rm\ or\ }1,$$
see (1.14)-(1.15).)

Subtracting (2.51) from (2.50), by (2.49) we have
$$\prod_{p|Z_0}\left(1-{1\over p}\right)^2
\cdot D\log {M^2\over T}
\left( 1-
\prod_{\scriptstyle p|Z_0\atop\scriptstyle \chi_{-D}(p)=1}
\left( 1+{1\over p-1}\right)
\prod_{\scriptstyle p|Z_0\atop\scriptstyle \chi_{-D}(p)=-1}
\left( 1-{1\over p+1}\right)\right) =$$
$$={\rm negligible\ error\ term}.\eqno(2.52)$$
According to the definition of $Z_0$ (see (1.14)-(1.15)), we distinguish
two cases. If $\chi_{-D}(p_0)=1$ then 
$$Z_0=\prod_{\scriptstyle p\le D:\atop\scriptstyle
\chi_{-D}(p)\ne -1}p,$$
and we can estimate (2.52) as follows:
$${\rm negligible\ error\ term}=\qquad\qquad\qquad\qquad $$
$$=\prod_{p|Z_0}\left(1-{1\over p}\right)^2
\cdot D\log {M^2\over T}
\left( 1-
\prod_{\scriptstyle p|Z_0\atop\scriptstyle \chi_{-D}(p)=1}
\left( 1+{1\over p-1}\right)\right) \le $$
$$\le \prod_{p|Z_0}\left(1-{1\over p}\right)^2
\cdot D\log {M^2\over T}
\left( 1-\left( 1+{1\over p_0-1}\right)\right) =$$
$$=- \prod_{p|Z_0}\left(1-{1\over p}\right)^2
\cdot D\log {M^2\over T}\cdot {1\over p_0-1}\le $$
$$\le - \prod_{p|Z_0}\left(1-{1\over p}\right)^2
\cdot D\log {M^2\over T}\cdot {1\over 2\log
D}.\eqno(2.53)$$

If $\chi_{-D}(p_0)=-1$ then we distinguish two subcases. If
$$\prod_{\scriptstyle p|Z_0\atop\scriptstyle \chi_{-D}(p)=1}
\left( 1+{1\over p-1}\right)\ge 1+{1\over 4\log D},$$
then we choose
$$Z_0=\prod_{\scriptstyle p\le D:\atop\scriptstyle
\chi_{-D}(p)\ne -1}p,$$
and  estimate (2.52) as follows:
$${\rm negligible\ error\ term}=\qquad\qquad\qquad\qquad $$
$$=\prod_{p|Z_0}\left(1-{1\over p}\right)^2
\cdot D\log {M^2\over T}
\left( 1-
\prod_{\scriptstyle p|Z_0\atop\scriptstyle \chi_{-D}(p)=1}
\left( 1+{1\over p-1}\right)\right) \le $$
$$\le \prod_{p|Z_0}\left(1-{1\over p}\right)^2
\cdot D\log {M^2\over T}
\left( 1-\left( 1+{1\over 4\log D}\right)\right) =$$
$$= - \prod_{p|Z_0}\left(1-{1\over p}\right)^2
\cdot D\log {M^2\over T}\cdot {1\over 4\log D}.\eqno(2.54)$$

Finally, 
if $\chi_{-D}(p_0)=-1$ and
$$\prod_{\scriptstyle p|Z_0\atop\scriptstyle \chi_{-D}(p)=1}
\left( 1+{1\over p-1}\right)< 1+{1\over 4\log D},$$
then we choose
$$Z_0=p_0\prod_{\scriptstyle p\le D:\atop\scriptstyle
\chi_{-D}(p)\ne -1}p,$$
and  estimate (2.52) as follows:
$${\rm negligible\ error\ term}=\qquad\qquad\qquad\qquad $$
$$=\prod_{p|Z_0}\left(1-{1\over p}\right)^2
\cdot D\log {M^2\over T}
\left( 1-
\prod_{\scriptstyle p|Z_0\atop\scriptstyle \chi_{-D}(p)=1}
\left( 1+{1\over p-1}\right)
\prod_{\scriptstyle p|Z_0\atop\scriptstyle \chi_{-D}(p)=-1}
\left( 1-{1\over p+1}\right)\right) \ge $$
$$\ge \prod_{p|Z_0}\left(1-{1\over p}\right)^2
\cdot D\log {M^2\over T}
\left( 1-\left( 1+{1\over 4\log D}\right)
\left( 1-{1\over p_0+1}\right)\right) \ge $$
$$\ge \prod_{p|Z_0}\left(1-{1\over p}\right)^2
\cdot D\log {M^2\over T}\cdot {1\over 4\log D},\eqno(2.55)$$
since $p_0<2\log D$.

Moreover, we use the well-known  number-theoretic fact
$$\prod_{p|Z_0}\left(1-{1\over q}\right)\ge 
\prod_{p\le D}\left(1-{1\over p}\right)
\ge {1\over 3\log D},$$
which implies
$$\prod_{p|Z_0}\left(1-{1\over q}\right)^2\ge
{1\over 10(\log D)^2}.\eqno(2.56)$$

Combining (2.52)-(2.56), we obtain
$${D\over 80(\log D)^3}\log {M^2\over T}=$$
$$={\rm negligible\ error\ term},\eqno(2.57)$$
which is a contradiction, since
${D\over 80(\log D)^3}\log (M^2/T)$
is {\it not} negligible. More precisely,
by using the choice (2.48) of the parameters, and the fact that
$D>e^{10^{100}}$ is very large, a simple calculation shows that
the left side of (2.57) is larger than the right side (=negligible
error term), which contradicts the equality.
(For more details, see the second half of Section 5.) 
This contradiction proves Theorem 1.
\medskip\noindent
{\bf Emphasizing the Simple Reason behind
the proof of Theorem 1.} 
I conclude Section 2 by emphasizing 
the crucial difference between
the left-hand side and the right-hand side of (2.49). In view of
(2.51), the right-hand side {\it distinguishes} between the primes
$p|Z_0$ with $\chi_{-D}(p)=1$, $0$ or $-1$; on the other hand, in view of
(2.50), the left-hand side does {\it not} distinguish between the primes
$p|Z_0$ with $\chi_{-D}(p)=1$, $0$ or $-1$
(here we ignore the negligible
error term). I will return to this key fact in the Concluding Remark
at the end of Section 3 and Section 4. The details of the proof of
Theorem 1 are complicated, so it is particularly
important to see the simple reason  why the
method works.

The proof of Theorem 2 is similar: we start with the analog of (2.49)
for positive discriminants, and study
the left-hand side and the right-hand side of (2.49). Again the contradiction
is based on the fact that, under the assumption that there is a
``bad'' positive fundamental discriminant $D$,
the left side and the right side cannot be equal. They are not equal,
because the right-hand side {\it distinguishes} between the primes
$p|Z_0$ with $\chi_{D}(p)=1$, $0$ or $-1$; 
more precisely, we prove an analog of (2.50). On the
 other hand, the
 left-hand side does {\it not} distinguish between the primes
$p|Z_0$ with $\chi_{D}(p)=1$, $0$ or $-1$
(here we ignore the negligible
error term).  This follows from an analog of (2.51).

\bigskip\bigskip\bigskip
\medskip
\centerline{\bf 3. Gauss sums and an application of smoothing }
\bigskip\bigskip\bigskip\medskip\noindent
Let's return to (2.49). By (2.42) we need to evaluate/estimate the sum
$$\left(\sum_{1\le a\le D} 
\chi_{-D}(a)
e^{2\pi{\rm i}a({\l_1\over m_1}-{\l_2\over m_2})}\right)
e^{2\pi{\rm i}D({\l_1\over m_1}-{\l_2\over m_2})\kappa N}
\left( {\sin\left( (2N+1)\pi D({\l_1\over m_1}-{\l_2\over
m_2})\right)\over (2N+1)\sin\left( \pi D({\l_1\over m_1}-{\l_2\over
m_2})\right)}\right)^{\kappa }$$
where  $\gcd (m_1m_2,Z_0)=1$.
If $(m_1,\l_1)\ne (m_2,\l_2)$ and $\gcd (m_1,m_2)=d$, then
$${\l_1\over m_1}-{\l_2\over m_2}={\l\over n}{\rm\ \ with\ \ some\ \ }
n\ge 1{\rm\ \ and\ \ }\l{\rm\ \ satisfying\ \ the\ \ properties}$$
$$m_1m_2d^{-2}|n{\rm\ \ and\ \ }
n|{\rm lcm}(m_1,m_2)=m_1m_2/d,{\rm\ \ and\ \ }\gcd (\l ,n)=1.\eqno(3.1)$$
That is, $m_1m_2d^{-2}$ is a divisor of  $n$, 
$n$ is a divisor of $m_1m_2/d$,
and, finally, $\l/n$ is an irreducible fraction.
Here we used the assumption $\gcd (\l_h,m_h)=1$, $h=1,2$
(see (2.42)), and, as usual,
${\rm lcm}(a,b)$ denotes the least common multiple of the positive
integers $a,b\ge 1$.

Let
$$D\l\equiv s\ ({\rm mod}\ n){\rm\ \ and\ \ }
1\le |s|\le n/2,\ \gcd (s,n)=1.\eqno(3.2)$$
Then we can estimate the ``tail factor'' in terms of $|s|$:
$$\left| e^{2\pi{\rm i}D({\l_1\over m_1}-{\l_2\over m_2})\kappa N}
\left( {\sin\left( (2N+1)\pi D({\l_1\over m_1}-{\l_2\over
m_2})\right)\over (2N+1)\sin\left( \pi D({\l_1\over m_1}-{\l_2\over
m_2})\right)}\right)^{\kappa }\right| =$$
$$=\left| 
\left( {\sin\left( (2N+1)\pi |s|/n\right)
\over (2N+1)\sin\left( \pi |s|/n\right)}\right)^{\kappa }\right| \le 
\min\left\{ 1,\left|{n\over sN}\right|^{\kappa}\right\} ,$$
since in the range $1\le |s|\le n/2$ we have the trivial upper bound 
$$\left| {\sin ( (2N+1)\pi s/n)
\over (2N+1)\sin ( \pi s/n)}\right| \le 
\min\left\{ 1,\left|{n\over sN}\right|\right\} .$$
This means that the ``tail factor'' is negligible
if $|s|$ is ``large''; so the main contribution comes from
the ``small'' values of $|s|$.

Next we focus on the critical sum
$$\sum_{1\le a\le D} \chi_{-D}(a)
e^{2\pi{\rm i}a({\l_1\over m_1}-{\l_2\over m_2})}
=\sum_{1\le a\le D} \chi_{-D}(a)
e^{2\pi{\rm i}a\l/n},$$
where $|s|$ in (3.2) is ``small''. The idea is to involve Gauss sums.
In the following  crucial lemma we use 
the notation introduced in (3.1)-(3.2).
\medskip\noindent
{\bf Lemma 3.1} {\it Let  
$n$ and $\l $ be integers such that
$\gcd (\l D,n)=1$. Then
$$\left|\sum_{1\le a\le D} \chi_{-D}(a)
e^{2\pi{\rm i}a\l/n}+{\rm i}
\chi_{-D}(s)\chi_{-D}(n)\sqrt{D}\right|\le
{2\pi D|s|\over n}$$
where $s\equiv D\l$
(mod $n$),
$1\le |s|\le n/2$, $\gcd (s,n)=1$, assuming 
$\gcd (s,D)=1$. 

If $\gcd (s,D)\ge 2$, then $\chi_{-D}(s)=0$, and we have the simpler
approximation
$$\left|\sum_{1\le a\le D} \chi_{-D}(a)
e^{2\pi{\rm i}a\l/n}\right|\le
{2\pi D|s|\over n}.$$}\medskip 
{\bf Remark.} We emphasize the fact that the error term is ``small''
if $|s|$ is ``small'' relative to $n$.
\medskip
{\bf Proof of Lemma 3.1.}
The smallest possible value of $|s|$ is 1---this case 
 deserves a special study. Let $s=1$;
since $D$ and $n$ are relatively prime, 
there exists an $\l^{\star}=\l^{\star}_n$ satisfying
$$D\l^{\star}\equiv 1\ ({\rm mod}\ n){\rm\ \ and\ \ }
1\le \l^{\star}<n,\eqno(3.3)$$
that is, $\l^{\star}$ is the multiplicative inverse of $D$ (mod $n$).
(3.3) is equivalent to $D\l^{\star}-1=nj^{\star}$ for some integer 
$1\le j^{\star}=j^{\star}_n< D$, or
equivalently, $D|nj^{\star}+1$ for some $1\le j^{\star}=j^{\star}_n< D$.
So $j^{\star}=j^{\star}_n$ and $D$ are relatively prime.
Thus we have
$$e^{2\pi{\rm i}a\l^{\star}/n}=e^{2\pi{\rm i}aD\l^{\star}/(Dn)}=
e^{2\pi{\rm i}a({D\l^{\star}-1\over Dn}+{1\over Dn})}=$$
$$=e^{2\pi{\rm i}a({nj^{\star}\over Dn}+{1\over Dn})}=
e^{2\pi{\rm i}aj^{\star}/D}\cdot e^{2\pi{\rm i}a/(Dn)}.\eqno(3.4)$$
If $n$ is ``large'' then the last factor is almost 1:
$$e^{2\pi{\rm i}a/(Dn)}=1+O(1/n).$$
This leads us to the evaluation
of the exponential sum 
$$\Theta (j^{\star})=\sum_{1\le a\le D}
\chi_{-D}(a)e^{2\pi{\rm
i}aj^{\star}/D}.\eqno(3.5)$$
Notice that (3.5) is a Gauss sum
$$G(\chi ;K)=\sum_{k=1}^K\chi (k)e^{2\pi{\rm i}k/K},\eqno(3.6)$$
where $\chi $ is a primitive Dirichlet character (mod $K$). We apply
 the following classical result due to Gauss.
\medskip\noindent
{\bf Lemma 3.2} {\it For all integers $K$ and $r$,  
$$\sum_{k=1}^K\chi (k)e^{2\pi{\rm i}rk/K}=\overline{\chi (r)}
G(\chi ;K),\eqno(3.7)$$
where
$$G(\chi ;K)={\rm i}\sqrt{K}{\rm\ \ if\ \ }\chi (-1)=-1.\eqno(3.8)$$}
\medskip
Applying Lemma 3.2 in (3.5),
$$\Theta (j^{\star})=\chi_{-D}(j^{\star}){\rm i}\sqrt{D}.
\eqno(3.9)$$
Since $D|nj^{\star}+1$, we have $\chi_{-D}(j^{\star})=
 -\chi_{-D}(n)$, and using it in
(3.9), 
$$\Theta (j^{\star})=-\chi_{-D}(n){\rm i}\sqrt{D}.
\eqno(3.10)$$

Next we switch from the special case $s=1$ to the general case of arbitrary $s$.
Let
$$D\l\equiv s\ ({\rm mod}\ n){\rm\ \ and\ \ }
1\le |s|\le n/2,\ \gcd (s,n)=1,\eqno(3.11)$$
then with $\l^{\star}=\l^{\star}_n$ defined in (3.3) we have
$$s\equiv Ds\l^{\star}\ ({\rm mod}\ n),\eqno(3.12)$$
so $\l\equiv s\l^{\star} $ (mod $n$), and
$$e^{2\pi{\rm i}a\l /n}=e^{2\pi{\rm i}as\l^{\star}/n}=
e^{2\pi{\rm i}aDs\l^{\star}/(Dn)}=$$
$$=e^{2\pi{\rm i}a({D\l^{\star}-1\over Dn}s+{s\over Dn})}=
e^{2\pi{\rm i}asj^{\star}/D}\cdot e^{2\pi{\rm i}as/(Dn)}.\eqno(3.13)$$
We need the simple estimation
$$\left| e^{2\pi{\rm i}as/(Dn)}-1\right|\le
{2\pi |s|\over n}.\eqno(3.14)$$
By using (3.14), we have ($\delta =1$ or $-1$)
$$\left|\sum_{\scriptstyle 1\le a\le
D\atop\scriptstyle \chi_{-D}(a)=\delta}e^{2\pi{\rm i}a\l/n}-
\sum_{\scriptstyle 1\le a\le
D\atop\scriptstyle \chi_{-D}(a)=\delta}e^{2\pi{\rm
i}asj^{\star}_n/D}\right|
\le {\pi D|s|\over n}.\eqno(3.15)$$
Write
$$\Theta (\l ;\delta)=\sum_{\scriptstyle 1\le a\le
D\atop\scriptstyle \chi_{-D}(a)=\delta}e^{2\pi{\rm
i}a\l/D}\eqno(3.16)$$
where $\l=sj^{\star}$ and $\delta =1$ or $-1$.
By Lemma 3.2,
$$\sum_{1\le a\le D}\chi_{-D}(a)e^{2\pi{\rm i}asj^{\star}_n/D}=
\Theta (sj^{\star};+1)-\Theta
 (sj^{\star};-1)=$$
$$=\chi_{-D}(sj^{\star} ){\rm i}\sqrt{D}
=-\chi_{-D}(s)\chi_{-D}(n){\rm i}\sqrt{D},\eqno(3.17)$$
if $\gcd (D,s)=1$, and 
$$\sum_{1\le a\le D}\chi_{-D}(a)
e^{2\pi{\rm i}asj^{\star}_n/D}=0\eqno(3.17')$$
if $\gcd (D,s)\ge 2$.

Finally notice that (3.15), (3.17) and (3.17') 
imply Lemma 3.1.\proofend
\medskip 
We are going to
apply Lemma 3.1 in (2.42). First note that
$${\rm real\ part\ of\ the\ product\ }
e^{2\pi{\rm i}\kappa Ns/n}
\left(-{\rm i}
\chi_{-D}(s)\chi_{-D}(n)\sqrt{D}\right) =$$
$$=\chi_{-D}(s)\chi_{-D}(n)\sqrt{D}
\sin (2\pi\kappa Ns/n).\eqno(3.18)$$
By (2.40)
$$\Omega_{1 ;\kappa ;N}(M)
-\Omega_{-1;\kappa ;N}(M)={\rm real},$$
so we have 
$${\rm real}=\Omega_{1 ;\kappa ;N}(M)
-\Omega_{-1;\kappa ;N}(M)=$$
$$=\Omega_{1 ;\kappa ;N}({\rm OffDiag};M)
-\Omega_{-1;\kappa ;N}({\rm OffDiag};M)=$$
$$={\rm real\ part \ of\ }
\left(\Omega_{1 ;\kappa ;N}({\rm OffDiag};M)
-\Omega_{-1 ;\kappa ;N}({\rm OffDiag};M)\right) =$$
$$=2OD_{\kappa ;N}(M)+
{\rm Error}(\kappa ;N;M),\eqno(3.19)$$
where, by applying Lemma 3.1 in (2.42), using (3.18), and also using
Lemma 3.3 (to be formulated below),
 the first term in the last line of (3.19) is
$$ OD_{\kappa ;N}(M)={\sqrt{D}\over 2}
\sum_{\scriptstyle 1\le d\le M:\ \mu^2(d)=1
\atop\scriptstyle \gcd (d,Z_0)=1}\ph(d)
\sum_{d_1|d}\cdot $$
$$\cdot\prod_{p|d_1}{p-2\over p-1} 
\sum_{\scriptstyle (m_1,m_2):\ 1\le m_1,m_2\le M\atop
\scriptstyle 
\gcd (m_1,m_2)=d,\gcd (m_1m_2,Z_0)=1}{\mu(m_1)\mu(m_2)
\over\ph (m_1)\ph (m_2)}\cdot\chi_{-D}({m_1m_2d_1\over d^2})\cdot $$
$$\cdot
\left( 2\sum_{\scriptstyle 1\le s\le {n\over 2}={m_1m_2d_1\over 2d^2}:\atop
\scriptstyle \gcd (s,n)=1}\chi_{-D}(s)\sin (2\pi\kappa Ns/n)
\left( {\sin ( (2N+1)\pi s/n)
\over (2N+1)\sin ( \pi s/n)}\right)^{\kappa }\right)
\eqno(3.20)$$
with $s,n,\l$ coming from 
$${\l_1\over m_1}-{\l_2\over
m_2}={\l \over n}{\rm\  with\ } n={m_1m_2d_1\over d^2},$$
$${\rm\ and\ }
s\equiv D\l\
({\rm mod\ } n),
1\le |s|\le n/2, \gcd (s,n)=1.\eqno(3.21)$$

In (3.20) we used the following simple but important lemma.
\medskip\noindent
{\bf Lemma 3.3} {\it Let $m_1\ge 1,m_2\ge 1$ be squarefree
integers, and let $d=\gcd (m_1,m_2)$ be their greatest common divisor.
Let $d_1\ge 1$ be an arbitrary divisor of $d$, and write 
$n=m_1m_2d^{-2}d_1$.  Let $\l $ be an  integer  with
$1\le\l< n$, $\gcd (\l ,n)=1$.
Let ${\cal A}(n,\l)$ denote the number of pairs 
$(\l_1 ,\l_2)$ of integers such that
$1\le\l_h\le m_h$, $\gcd (\l_h ,m_h)=1$, $h=1,2$, and
$${\l_1 \over m_1}-{\l_2 \over m_2}={\l \over n}
{\rm\ \ or\ \ }{\l  \over n}-1.$$
Then
$${\cal A}(n,\l)=\ph_2(d_1)\ph(d/d_1)=\ph(d)\prod_{p|d_1}{p-2\over p-1}, $$
where
$$\ph (m)=m\prod_{p|m}\left( 1-{1\over p}\right){\rm\ \ 
and\ \ }
\ph_2 (m)=m\prod_{p|m}\left( 1-{2\over p}\right) .$$}
\medskip
I postpone the proof of Lemma 3.3 to the end of Section 3.

The second term in the last line of (3.19) is
the  error:
$$\left|{\rm Error}(\kappa ;N;M)\right|\le 2
\sum_{1\le d\le M}d\sum_{d_1|d}
\sum_{\scriptstyle (m_1,m_2):\ 1\le m_1,m_2\le M\atop
\scriptstyle\gcd (m_1,m_2)=d}{1\over\ph (m_1)\ph (m_2)}\cdot $$
$$\cdot 
\left( \tau (D)\sum_{1\le s\le {n\over 2}={m_1m_2d_1\over 2d^2}}
{2\pi Ds\over n}\cdot
\min\left\{ 1,\left({n\over sN}\right)^{\kappa}\right\}
\right) ,\eqno(3.22)$$
since in the range $1\le s\le n/2$ we have the trivial upper bound 
$$\left| {\sin ( (2N+1)\pi s/n)
\over (2N+1)\sin ( \pi s/n)}\right| \le 
\min\left\{ 1,\left|{n\over sN}\right|\right\}
.\eqno(3.23)$$

 I estimate the error term (3.22) from above.
This is the first application of the efficient average (``smoothing'')
introduced in Section 2.
\medskip\noindent
{\bf Lemma 3.4} {\it If $\kappa \ge 3$ then
$$\left|{\rm Error}(\kappa ;N;M)\right| <1 .$$}
\medskip
{\bf Proof of Lemma 3.4.} 
We distinguish two cases in (3.22) depending on whether
$$ 1\le {n\over sN} {\rm\ \ or\ \ } 1> {n\over sN},$$
which implies
$$\left|{\rm Error}(\kappa ;N;M)\right|\le
{\rm E}_1+{\rm E}_2,\eqno(3.24)$$
where (let $m'_i=m_i/d$, $i=1,2$ and $n=m'_1m'_2d_1$)
$${\rm E}_1=2
\sum_{1\le d\le M}{d\over\ph^2 (d)}\sum_{d_1|d}
\sum_{\scriptstyle (m'_1,m'_2):\atop
\scriptstyle 1\le m'_1, m'_2\le M/d}{1\over\ph (m'_1)\ph (m'_2)}\cdot 
 2\pi \tau (D)D\left( \sum_{1\le s\le {m'_1m'_2d_1\over N}}
{s\over m'_1m'_2d_1}\right) ,\eqno(3.25)$$
and
$${\rm E}_2=2\sum_{1\le d\le M}{d\over\ph^2 (d)}\sum_{d_1|d}
\sum_{\scriptstyle (m'_1,m'_2):\atop
\scriptstyle 1\le m'_1,m'_2\le M/d}{1\over\ph (m'_1)\ph (m'_2)}\cdot $$
$$\cdot 2\pi \tau (D)D
\left( \sum_{{m'_1m'_2d_1\over N}<s\le m'_1m'_2d_1/2}
{1\over N}\cdot 
\left({m'_1m'_2d\over sN}\right)^{\kappa -1}\right)
.\eqno(3.26)$$
Using the trivial upper bound
$$ \sum_{1\le s\le {m'_1m'_2d_1\over N}}
{s\over m'_1m'_2d_1}\le
{1\over m'_1m'_2d_1}\left( {m'_1m'_2d_1\over N}\right)^2=
{m'_1m'_2d_1\over N^2}$$
in (3.25), we have
$${\rm E}_1\le {4\pi \tau (D)D\over N^2}
\sum_{1\le d\le M}\sum_{d_1|d}{dd_1\over\ph^2 (d)}
\sum_{\scriptstyle (m'_1,m'_2):\atop
\scriptstyle 1\le m'_1,m'_2\le M/d}{m'_1m'_2\over\ph (m'_1)
\ph (m'_2)}.\eqno(3.27)$$
Here we need (and repeatedly use later) the following well-known fact
from elementary number theory:
$${n\over\ph (n)}\le 10\log\log n\eqno(3.28)$$
for all $n\ge 100$.
By (3.28),
$$\sum_{\scriptstyle (m'_1,m'_2):\atop
\scriptstyle 1\le m'_1,m'_2\le M/d}{m'_1m'_2\over\ph (m'_1)
\ph (m'_2)}\le \left(
\sum_{1\le m\le M/d}{m\over\ph (m)}\right)^2 \le 
10^4\left( {M\over d}\log\log M\right)^2.\eqno(3.29)$$
By using (3.29) in (3.27),
$${\rm E}_1\le {4\pi \tau (D)D\over N^2}
\sum_{1\le d\le M}\sum_{d_1|d}{dd_1\over\ph^2 (d)}
\cdot {10^4M^2(\log\log M)^2\over d^2}\le $$
$$\le {2\pi \tau (D)D\over N^2}\cdot
10^4\left( {M\over d}\log\log M\right)^2
\left(\sum_{d\ge 1}{\tau (d)\over\ph^2(d)}\right)\le $$
$$\le {10^7\tau (D)DM^2(\log\log M)^2\over N^2},\eqno(3.30)$$
since the infinite series in (3.30) is clearly convergent
(note that $\tau (n)=n^{o(1)}$):
$$\sum_{d\ge 1}{\tau (d)\over\ph^2(d)}\le 100. $$

To estimate ${\rm E}_2$,
 we use a standard power-of-two decomposition in (3.26): 
$${\rm E}_2\le 4\pi \tau (D)D
\sum_{1\le d\le M}{d\over\ph^2 (d)}
\sum_{\scriptstyle (m'_1,m'_2):\atop
\scriptstyle 1\le m'_1,m'_2\le M/d}
{1\over\ph (m'_1)\ph (m'_2)}\cdot $$
$$\cdot \sum_{d_1|d}
\left( \sum_{j\ge 1}\sum_{{m'_1m'_2d_1\over N}2^{j-1}<s\le 
{m'_1m'_2d_1\over N}2^j}{1\over N}\cdot 
\left({m'_1m'_2d_1\over sN}\right)^{\kappa -1}\right)\le $$
$$\le 4\pi \tau (D)D\sum_{1\le d\le M}{d\over\ph^2 (d)}
\sum_{\scriptstyle (m'_1,m'_2):\atop
\scriptstyle 1\le m'_1,m'_2\le M/d}{1\over\ph (m'_1)\ph
(m'_2)}\cdot $$
$$\cdot
\sum_{d_1|d}\left( \sum_{j\ge 1} 
{m'_1m'_2d_1\over N}2^j\cdot{1\over N}\cdot 
2^{-(\kappa -1)(j-1)}\right)\le $$
$$\le  {4\pi \tau (D)D\over N^2}
\sum_{1\le d\le M}\sum_{d_1|d}{dd_1\over\ph^2 (d)}
\sum_{\scriptstyle (m'_1,m'_2):\atop
\scriptstyle 1\le m'_1,m'_2\le M/d}{m'_1m'_2\over
\ph (m'_1)\ph (m'_2)}\cdot \left( \sum_{j\ge 1} 
2^{-(\kappa -2)(j-1)+1}\right)\le $$
$$\le  {16\pi\tau (D) D\over N^2}
\sum_{1\le d\le M}{d^2\tau (d)\over\ph^2 (d)}
\left(\sum_{1\le n\le M/d}{n\over\ph (n)}\right)^2 \le $$
$$\le {5\cdot 10^5\tau (D)D
(M\log\log M)^2\over N^2}\sum_{1\le d\le M}{\tau (d)
\over \ph^2(d)}\le $$
$$\le {5\cdot 10^7\tau (D)DM^2(\log\log M)^2\over N^2},\eqno(3.31)$$
assuming $\kappa \ge 3$
(which guarantees the convergence of the series
$\sum_{j\ge 1} 
2^{-(\kappa -2)j}$).
In the last steps of the argument in (3.31) 
we used (3.28) the same way as we did in (3.29).

Finally notice that  (3.24), (3.30) and (3.31) imply Lemma 3.4:
$$\left|{\rm Error}(\kappa ;N;M)\right| \le 
{10^8\tau (D)DM^2\over N^2}(\log\log M)^2< 1,$$
where the last step is a trivial calculation using the values of the
parameters (see (2.48)):
$$\kappa =8,\ \log T=(\log D)^{15},\ T\le N<2T,
\ M=T\exp\left(-{2\over 3}\sqrt{\log T}\right) ,
\ \log D>10^{100}.$$
\proofend
\medskip
As we promised, we include  the proof of Lemma 3.3.

{\bf Proof of Lemma 3.3.} 
Let ${\cal B}(n,\l)$ denote the number of pairs
of residue classes
$(\l_1\ ({\rm mod} \ m_1) ,\l_2\ ({\rm mod} \ m_2))$  such that
$1\le\l_h\le m_h$, $\gcd (\l_h ,m_h)=1$, $h=1,2$, and
there exist integers $\l^{(1)}$ and $\l^{(2)}$
with the property
$$\l^{(1)}\equiv\l_1\ ({\rm mod} \ m_1) ,\ \ 
\l^{(2)}\equiv\l_2\ ({\rm mod} \ m_2),$$  and
$${\l^{(1)} \over m_1}-{\l^{(2)} \over m_2}={\l \over n}.$$
Clearly ${\cal A}(n,\l)\le {\cal B}(n,\l)$.

Next I prove the upper bound
$${\cal B}(n,\l)\le
\ph_2(d_1)\ph(d/d_1)=\ph(d)\prod_{p|d_1}{p-2\over p-1},\eqno(3.32)$$
where $n=m'_1m'_2d_1$.

If
$${\l \over n}= {\l^{(1)}
 \over m_1}-{\l^{(2)} \over m_2}=
{\l^{(3)} \over m_1}-{\l^{(4)} \over m_2}$$
holds for some integers $\l^{(1)} ,\l^{(2)},\l^{(3)},
\l^{(4)}$ with 
$1\le\l^{(h)}\le m_h$, $\gcd (\l^{(h)} ,m_h)=1$, $h=1,2$
and $1\le\l^{(h+2)}\le m_h$, $\gcd (\l^{(h+2)} ,m_h)=1$, $h=1,2$,
then $(\l^{(1)}-\l^{(3)})m'_2=(\l^{(2)}-\l^{(4)})m'_1$, 
where $m'_h=m_h/d$, $h=1,2$.
Since $m'_1$ and $m'_2$ are coprime, we have
$\l^{(1)}-\l^{(3)}=jm'_1$ and $\l^{(2)}-\l^{(4)}=jm'_2$ with 
the same integer $j$. It follows that the residue class
$j\ ({\rm mod} \ d)$ uniquely determines both residue classes
$$\l^{(1)}-\l^{(3)}\ ({\rm mod} \ m_1){\rm\ \ 
 and\ \ } \l^{(2)}-\l^{(4)}\ ({\rm mod} \ m_2).\eqno(3.33)$$

We now study what restrictions apply to the residue class  
$j\ ({\rm mod} \ d)$.

First, let $p$ be a prime divisor of $d_1$
(I recall that $n=m'_1m'_2d_1$). The equalities
$\l^{(3)}m'_2=\l^{(1)}m'_2-jm'_1m'_2$ and
$\l^{(4)}m'_1=\l^{(2)}m'_1-jm'_1m'_2$ imply
$$\l^{(3)}m'_2\equiv \l^{(1)}m'_2-jm'_1m'_2\ ({\rm mod} 
\ p){\rm \ \ and\ \ } 
\l^{(4)}m'_1\equiv \l^{(2)}m'_1-jm'_1m'_2\ ({\rm mod} \ p).$$
Note that $\l^{(3)}m'_2\not\equiv 0\ ({\rm mod} \ p)$, since otherwise 
$p|\l^{(3)}$, 
implying that both $\l^{(3)}$ and $m_1$ are divisible by $p$, which
contradicts the fact that they are coprime. Similarly,
$\l^{(4)}m'_1\not\equiv 0\ ({\rm mod} \ p)$. 

Furthermore,
$\l^{(3)}m'_2-\l^{(4)}m'_1$
is not divisible by $p$, since otherwise we could simplify the
fraction
$${\l \over n}= {\l^{(3)} \over m_1}-{\l^{(4)} \over m_2}=
{\l^{(3)}m'_2-\l^{(4)}m'_1\over m'_1m'_2d}$$
by $p$, which is a contradiction.

This means that for $j$ we have 
{\it two forbidden} residue classes $({\rm mod} \ p)$: 
$$ j\not\equiv \l^{(1)}/m'_1\ ({\rm mod} \ p){\rm\ \ and\ \ }
 j\not\equiv \l^{(2)}/m'_2\ ({\rm mod} \ p).$$

On the other hand, if $q$ is a prime divisor of $d/d_1$, then the same
argument gives that there is only {\it one forbidden} residue
 class $({\rm mod} \ q)$.

Therefore, we have 
$$\ph_2(d_1)\ph(d/d_1)=\ph (d)
\prod_{p|d_1}{p-2\over p-1}\eqno(3.34)$$
available residue classes of $j\ ({\rm mod} \ d)$.
Combining (3.33) and (3.34), (3.32) follows.

Summarizing, for $m_1\ne m_2$ we have
$$\sum_{n,\l}{\cal A}(n,\l)\le\sum_{n,\l}{\cal B}(n,\l)
\le \sum_{n=m'_1m'_2d_1}\sum_{\l}\ph_2(d_1)\ph(d/d_1)=$$
$$=\sum_{d_1|d}\ph (m'_1m'_2)\ph (d_1)\cdot \ph_2(d_1)\ph(d/d_1)=
\sum_{d_1|d}\ph (m'_1m'_2)\ph (d_1)\cdot 
\ph (d)\prod_{p|d_1}{p-2\over p-1}=$$
$$=\ph (m'_1m'_2)\ph (d) \sum_{d_1|d}\prod_{p|d_1}(p-2)=
\ph (n)\prod_{p|d}(1+p-2)=$$
$$=\ph (m'_1m'_2)\ph^2 (d)=\ph (m_1)\ph (m_2).\eqno(3.35)$$
On the other hand,
$$\sum_{n,\l}{\cal A}(n,\l)=\ph (m_1)\ph (m_2),\eqno(3.36)$$
since the right-hand side of (3.36) is the number of all possible 
differences
$${\l_1\over m_1}-{\l_2\over m_2},$$
assuming $m_1\ne m_2$.

Combining (3.35) and (3.36), we have with $n=m'_1m'_2d_1$,
$$\ph (m_1)\ph (m_2)=
\sum_{n,\l}{\cal A}(n,\l)\le\sum_{n,\l}{\cal B}(n,\l)
\le \sum_{n=m'_1m'_2d_1}\sum_{\l}\ph_2(d_1)\ph(d/d_1)
\le \ph (m_1)\ph (m_2).\eqno(3.37)$$
Since
$${\cal A}(n,\l)\le {\cal B}(n,\l)
\le \ph_2(d_1)\ph(d/d_1),$$
(3.37) upgrades the inequality to equality:
$${\cal A}(n,\l)=\ph_2(d_1)\ph(d/d_1).$$

Finally, in the trivial case
 $m_1=m_2$ we have $m_1=m_2=d$ and $n=d_1$.
Then we have to exclude the case $d_1=1$ (since $(m_1,\l_1) \ne
(m_2,\l_2)$) in (3.35), and proceed as follows:
$$\sum_{\scriptstyle (d_1,\l):\atop
\scriptstyle d_1|d,d_1\ne 1}{\cal A}(d_1,\l)
\le \sum_{\scriptstyle (d_1,\l):\atop
\scriptstyle d_1|d,d_1\ne 1}{\cal B}(d_1,\l)
\le \sum_{d_1|d:\ d_1\ne 1}\sum_{\l}\ph_2(d_1)\ph(d/d_1)=$$
$$=\sum_{d_1|d:\ d_1\ne 1}\ph (d_1)\cdot \ph_2(d_1)\ph(d/d_1)=
\left(\sum_{d_1|d}\ph (d_1)\cdot \ph_2(d_1)\ph(d/d_1)\right) -
\ph (d)=\ph^2 (d)-\ph (d).\eqno(3.38)$$
On the other hand,
$$\sum_{\scriptstyle (d_1,\l):\atop
\scriptstyle d_1|d,d_1\ne 1}{\cal A}(d_1,\l)=
\ph^2 (d)-\ph (d),\eqno(3.39)$$
since the right-hand side of (3.39) is the number of all possible 
differences
$${\l_1\over d}-{\l_2\over d}$$
with $\l_1\ne\l_2$.
Again combining (3.38) and (3.39), we can upgrade the inequality to
equality, and the proof of Lemma 3.3 is complete.\proofend
\medskip\noindent
{\bf Concluding Remark of Section 3.}
Let us return to (3.19)-(3.20).
The message of Lemma 3.4 is that
${\rm Error}(\kappa ;N;M)$ is negligible compared to
$OD_{\kappa ;N}(M)$, so $2OD_{\kappa ;N}(M)$ is the dominating
part of 
$$\Omega_{1 ;\kappa ;N}({\rm OffDiag};M)-
\Omega_{-1 ;\kappa ;N}({\rm OffDiag};M).$$
At the end of Section 2 I made the claim that
$$\overline{\Omega}_{1 ;\kappa ;T}({\rm OffDiag};M)-
\overline{\Omega}_{-1 ;\kappa ;T}({\rm OffDiag};M)$$
(which is simply the average of 
$$\Omega_{1 ;\kappa ;N}({\rm OffDiag};M)-
\Omega_{-1 ;\kappa ;N}({\rm OffDiag};M)$$
as $N$ runs in $T\le N<2T$; see (2.46)) 
{\it distinguishes} between the primes
$p|Z_0$ with $\chi_{-D}(p)=1$, $0$ or $-1$. 
This claim is intuitively well justified
by the last line in (3.20)
$$\sum_{\scriptstyle 1\le s\le {n\over 2}={m_1m_2d_1\over 2d^2}:\atop
\scriptstyle \gcd (s,n)=1}\chi_{-D}(s)\sin (2\pi\kappa Ns/n)
\left( {\sin ( (2N+1)\pi s/n)
\over (2N+1)\sin ( \pi s/n)}\right)^{\kappa },$$
due to the appearance of the character $\chi_{-D}(s)$.

Note in advance that the evaluation of the sum (3.20) will eventually
involve, or rather lead to, a sum like
$\sum_{j=1}^{D^4}\chi_{-D}(j)\log j$, which will become a factor in
the dominating part of (3.20)
(note that the power $D^4$ here is accidental; any not too small
 power of $D$ would do).
The reason behind it is that the sum  in (3.20)
resembles a ``double harmonic sum'', and also the main term in Lemma
2.1 is the logarithm function. This explains why (3.20) resembles a
Riemann sum for a logarithmic integral that we can evaluate
explicitly. 
The explicitly evaluated logarithmic integral can be well
 illustrated with the example
$$\sum_{j=1}^{D^4}\chi_{-D}(j)\left(\log {M^2\over Nj}\right)^2=
\sum_{j=1}^{D^4}\chi_{-D}(j)\left(\log {M^2\over N}\right)^2+$$
$$+\sum_{j=1}^{D^4}\chi_{-D}(j)(\log j)^2 -
2\log {M^2\over N}\sum_{j=1}^{D^4}\chi_{-D}(j)\log j=$$
$$=-2\log {M^2\over N}\sum_{j=1}^{D^4}\chi_{-D}(j)\log j
+{\rm negligible }.\eqno(3.40)$$
The fact that the sum $\sum_{j=1}^{D^4}\chi_{-D}(j)(\log j)^2 $ 
represents a negligible contribution here follows from the P\'olya--Vinogradov
inequality and partial summation.


Note that the harmonic sum
$$\sum_{n=1}^M{1\over n}=\log M +O(1){\rm\ \  if\ \ }
 M\ge 1{\rm\ \  and\ \ } 0{\rm\ \  if\ \ } M<1$$
is actually associated with the function 
$\log_+ x$ instead of $\log x$
($\log_+ x=\log x$ if $x\ge 1$ and 0 if $0<x<1$).
On the other hand, in (3.40) we work with $\log x$
instead of $\log_+ x$. The intuitively plausible assumption that
``$\log_+ x$ and $\log x$ are basically the same''
requires a precise proof applied in our particular case. Unfortunately,
this part of the proof turns out to be
annoyingly cumbersome. This kind of technical problems explain why
the paper is so long.

Let's return to (3.40).
If the class number $h(-D)$ is ``substantially smaller'' than
$\sqrt{D}$, then we have the good approximation
$$\sum_{j=1}^{D^4}\chi_{-D}(j)\log j=-{\pi\over 6}\sqrt{D}
\sum_{(a,b,c)}{1\over a}\ +\ {\rm negligible}$$
(see Lemma 14.2),
where $\sum_{(a,b,c)}{1\over a}$ means that
we add up the reciprocals of the leading coefficients 
$a=a_j$ in the family of
reduced, primitive, inequivalent binary quadratic forms of integer
coefficients $ax^2+bxy+cy^2$ with discriminant
$b^2-4ac=-D<0$,  $a>0$, $c>0$, $1\le j\le h(-D)$.

Moreover, under the same condition, we have another good approximation
(see Lemma 14.3):
$$\sum_{j=1}^{h(-D)}{1\over a_j}=\prod_{\scriptstyle p|Z_0
\atop\scriptstyle \chi_{-D}(p)\ne -1}
{p+1\over p-1}\prod_{p|D}{p-1\over p}\ +\ {\rm negligible}.$$
Notice that the product
$$\prod_{\scriptstyle p|Z_0
\atop\scriptstyle \chi_{-D}(p)\ne -1}
{p+1\over p-1}\prod_{p|D}{p-1\over p}$$
clearly {\it distinguishes} between the primes
$p|Z_0$ with $\chi_{-D}(p)=1$, 0 or $-1$.

On the other hand, (2.50) does not distinguish between the primes
$p|Z_0$ with $\chi_{-D}(p)=1$, 0 or $-1$.
This is how we obtain a contradiction, which proves Theorem 1.

The proof of (2.50) is not trivial (see Sections 4, 7 and
8). What {\it is} trivial is the case of primes; see the Remark after (2.3):
the factor $\prod_{p|Z_0}\left( 1-{1\over p}\right)$
in (2.4)  does not distinguish between the primes
$p|Z_0$ with $\chi_{-D}(p)=1$, 0 or $-1$. Since the primes play a key
role in the proof of (2.50), it is not surprising that
(2.50) does not distinguish between the primes
$p|Z_0$ with $\chi_{-D}(p)=1$, 0 or $-1$.

Note in advance that the case of positive discriminants (Theorem 2)
is somewhat different. Instead of the sum
$\sum_{j=1}^{D^4}\chi_{-D}(j)\log j$, we make use of the other sum
$\sum_{j=1}^{D^4}\chi_{D}(j)(\log j)^2$.
(For $D>0$ the sum
$\sum_{j=1}^{D^4}\chi_{D}(j)\log j$ turns out to be ``negligible'' if
if $L(1,\chi_{D})$ is ``small''.)
 We have the analog approximation
$$\sum_{j=1}^{D^4}\chi_{D}(j)(\log j)^2=-{\pi^2\over 6}\sqrt{D}
\prod_{\scriptstyle p|Z_0
\atop\scriptstyle \chi_{D}(p)\ne -1}
{p+1\over p-1}\prod_{p|D}{p-1\over p}\ +\ {\rm negligible}$$
if $L(1,\chi_{D})$ is ``small''.
Again the key observation is that the product on the right-hand side
clearly {\it distinguishes} between the primes
$p|Z_0$ with $\chi_{D}(p)=1$, 0 or $-1$.

\bigskip\bigskip\bigskip
\medskip
\centerline{\bf 4. The method of the proof (I)}
\bigskip\bigskip\bigskip\medskip\noindent
As I said at the beginning of Section 2,
the basic idea of the proof of Theorem 1 is to
 evaluate/estimate the two sides of the equality (2.49)
with the choice (2.48) of the key parameters.
First we evaluate the left side of (2.49) 
by using (2.40) and (2.44);  see Lemma 5.4. 
This is the subject of Section 4 and the beginning of
Section 5.
 
Then we evaluate the right side of (2.49) by
using (2.42) and (2.45),
and we obtain a {\it different} expression.
This is the subject of the middle part of Section 5; see
 Lemma 5.5. Finally, at the end of Section 5  
we derive a contradiction, which completes the proof of Theorem 1.

The details go as follows.
We want to find a good estimation for the difference
$$\overline{\Omega}_{1 ;\kappa ;T}(M)-
\overline{\Omega}_{-1 ;\kappa ;T}(M).$$
The main difficulty is how to estimate (2.40).

By (2.41)-(2.43) and (3.19),
$$\Omega_{1 ;\kappa ;N}(M)-\Omega_{-1 ;\kappa ;N}(M)=
\Omega_{1 ;\kappa ;N}({\rm OffDiag};M)-
\Omega_{-1 ;\kappa ;N}({\rm OffDiag};M)=$$
$$=2OD_{\kappa ;N}(M)+{\rm Error}(\kappa ;N;M).\eqno(4.1)$$

By Lemma 3.4, for $\kappa \ge 3$ we have
$$\left|{\rm Error}(\kappa ;N;M)\right|  \le
{10^7\tau (D)DM^2\over N^2}\log M (\log\log M)^4.\eqno(4.2)$$

We also need the special case $Q=1$ in Lemma 2.1:  
$$\left|\sum_{\scriptstyle 1\le m\le M:
\atop\scriptstyle\gcd (Z_0,m)=1}
{|\mu (m)|\over\ph (m)}-
\prod_{q|Z_0}\left( 1-{1\over q}\right)\left(
\log M+c'+\sum_{q|Z_0}{\log q\over q}
\right)\right|\le $$
$$\le {10^4\log D\log M\over M^{1/4}}+
{10^5\over D^5}+
{4\left( 10+\log M+2(\log D)^2\right) \over
\max\left\{ MD^{-6\log D},1\right\} } ,\eqno(4.3)$$
where $c'=\gamma_0+2\gamma^{\star}-\gamma^{\star\star}$.
Note that trivially
$$\sum_{q|Z_0}{\log q\over q}
\le \sum_{q\le D}{\log q\over q}\le (\log D)^2.
\eqno(4.4)$$

By (2.40),
$$\Omega_{1 ;\kappa ;N}(M)-\Omega_{-1 ;\kappa ;N}(M)=$$
$$=\sum_{k=0}^{\kappa N}w^{(\kappa ,N)}_k
{1\over 2k+1}\sum_{j:\ \kappa N-k\le j\le \kappa
N+k}\sum_{\scriptstyle 1\le a\le D:\atop\scriptstyle
\gcd (a,D)=1}\chi_{-D}(a)S^2(M;a+jD)=$$
$$={\sum}_{n=p^r}(N;M)+{\sum}_{n\ne p^r}(N;M)
,\eqno(4.5)$$
where
$${\sum}_{n=p^r}(N;M),\ \ {\sum}_{n\ne p^r}(N;M){\rm\ \  and\ \ }
S(M;a+jD)$$
are defined in the following way.
The first sum in the last line of (4.5) is restricted to the primepowers
(of course every prime {\it is} a primepower):
$${\sum}_{n=p^r}(N;M)=\sum_{k=0}^{\kappa N}w^{(\kappa ,N)}_k
{1\over 2k+1}\sum_{\scriptstyle (\kappa N-k)D<n< (\kappa
N+k+1)D:\atop\scriptstyle n=p^r{\rm\ for\ some\ prime\ }
p{\rm\ and\ integer\ }r\ge 1}
\chi_{-D}(p^r)S^2(M;p^r),\eqno(4.6)$$
the second sum in the last line of (4.5)
is restricted to the rest of the integers:
$${\sum}_{n\ne p^r}(N;M)=
\sum_{k=0}^{\kappa N}w^{(\kappa ,N)}_k
{1\over 2k+1}
\sum_{\scriptstyle (\kappa N-k)D<n< (\kappa
N+k+1)D:\atop\scriptstyle n{\rm\ is\ not\ a\ primepower}}
\chi_{-D}(n)S^2(M;n),\eqno(4.7)$$
and, finally, $S(M;n)$ is defined as 
$$S(M;n)=
\sum_{\scriptstyle d\ge 1:\gcd (Z_0,d)=1\atop\scriptstyle
d|n}\mu (d)
\sum_{\scriptstyle 1\le k\le M/d:\atop
\scriptstyle \gcd (k,nZ_0)=1}
{|\mu (k)|\over\ph (k)}.\eqno(4.8)$$

The evaluation of 
$S(M;p)$ is quite easy. Indeed, in the case $n=p$ prime we have
$$S(M;p)=
\sum_{\scriptstyle d\ge 1:\atop\scriptstyle
d|p}\mu (d)
\sum_{\scriptstyle 1\le k\le M/d:\atop
\scriptstyle \gcd (k,pZ_0)=1}
{|\mu (k)|\over\ph (k)}=
\sum_{\scriptstyle 1\le k\le M:
\atop\scriptstyle\gcd (Z_0, k)=1}{|\mu 
(k)|\over\ph (k)}{\rm\ \ if\ \ }p>M,\eqno(4.9)$$
and
$$S(M;p)=
\sum_{\scriptstyle d\ge 1:\atop\scriptstyle
d|p}\mu (d)
\sum_{\scriptstyle 1\le k\le M/d:\atop
\scriptstyle \gcd (k,pZ_0)=1}
{|\mu (k)|\over\ph (k)}=$$
$$=\sum_{\scriptstyle 1\le k\le M:
\atop\scriptstyle\gcd (Z_0,k)=1}{|\mu (k)|\over\ph (k)}
-\sum_{\scriptstyle 1\le k\le M/p:\atop
\scriptstyle \gcd (k,pZ_0)=1}
{|\mu (k)|\over\ph (k)}
{\rm\ \ if\ \ }1<p\le M.\eqno(4.10)$$
Using (4.9)-(4.10) and  (4.3)-(4.4)
(note that $M>D^{9\log D}$),
$$\Biggl|{\sum}_{n=p^r}(N;M)-\sum_{k=0}^{\kappa N}w^{(\kappa ,N)}_k
{1\over 2k+1}\sum_{(\kappa N-k)D<p<(\kappa N+k+1)D}\chi_{-D}(p)\cdot $$
$$\qquad\qquad\qquad \cdot\prod_{q|Z_0}\left( 1-{1\over q}\right)^2\left(
\log M+c'+\sum_{q|Z_0}{\log q\over q}
\right)^2\Biggr|\le $$
$$\le 1+
\sum_{k=\kappa N-{M\over D}}^{\kappa N}w^{(\kappa ,N)}_k
{1\over 2k+1}M(\log M)^2\le  1+{M(\log M)^2\over N-M}<2,
\eqno(4.11)$$
where we used (2.48), and the simple fact that---roughly 
speaking---the overwhelming majority of
primepowers are primes.
More precisely, we used the statement that
the number of integers of the form $a^b\le x$, where
$a\ge 2$, $b\ge 2$ are integers, is clearly less than
$$x^{1/2}+x^{1/3}+x^{1/5}+x^{1/7}+x^{1/11}
+\ldots = \sum_{2\le p\le\log x}x^{1/p}\le 2\sqrt{x}$$
if $x\ge 100$.

Since we have very good estimations for the number of primes in 
{\it long} intervals, the estimations become (basically) trivial
 if we turn $N$ into
a variable running in the long interval $T\le N<2T$ and take 
the average
as follows:
$$\Biggl|{1\over T}\sum_{T\le N<2T}
\sum_{k=0}^{\kappa N}w^{(\kappa ,N)}_k
{1\over 2k+1}\sum_{(\kappa N-k)D<p<(\kappa N+k+1)D}\chi_{-D}(p)\cdot
\qquad\qquad\qquad $$
$$\qquad\cdot
\prod_{q|Z_0}\left( 1-{1\over q}\right)^2\left(
\log M+c'+\sum_{q|Z_0}{\log q\over q}
\right)^2 +
\prod_{q|Z_0}\left( 1-{1\over q}\right)^2{D(\log M)^2\over\log
T}\Biggr|\le $$
$$\le {10D\log D(\log M)^2\over (\log T)^2}
+{10D\log D\log M\over \log T}+
100D(\log M)^2L(1,\chi_{-D}) \le $$
$$\le 30D\log D,\eqno(4.12)$$
where we used the Prime Number Theorem, Lemma 1.1 with (1.11), (4.4),
(2.48) and (1.10).
(The choice  $\kappa =8$ implies that $w^{(\kappa ,N)}_k$, $0\le
k\le\kappa N$ is a ``reasonable'' probability distribution, explicitly
described in Section 2.)

Write
$$\overline{\sum}_{n=p^r}(T;M)={1\over T}\sum_{T\le N<2T}
{\sum}_{n=p^r}(N;M).\eqno(4.13)$$
Combining (4.11)-(4.13), we obtain
\medskip\noindent
{\bf Lemma 4.1} {\it We have
$$\left|\overline{ \sum}_{n=p^r}(T;M)+
\prod_{q|Z_0}\left( 1-{1\over q}\right)^2{D(\log M)^2\over\log
T}\right|\le $$
$$\le 30D\log D+2.$$}
\proofend\medskip
The estimation of the second sum
${\sum}_{n\ne p^r}(N;M)$ in (4.7)
is based on Lemmas 4.2-3-4 below.
To estimate ${\sum}_{n\ne p^r}(N;M)$, we focus on the
end-sum $S(M;n)$; see (4.8).
Lemma 2.1 implies the trivial upper bound
$$|S(M;n)|\le 20\tau (n)\log M,\eqno(4.14)$$
which will be constantly used below.

In Sections 4 and 7-8 we employ the following unusual notation:
for every integer $n\ge 1$, let $n^*$ denote the largest squarefree
divisor of $n$. That is,  $n^*=n$ holds if $n$ is squarefree, and in general,
$$n^*=p_1\cdots p_r{\rm\ \ if\ \ }n=p_1^{a_1}\cdots p_r^{a_r}
{\rm\ \ with\ \ }a_i\ge 1,\ 1\le i\le r$$
is the prime factorization of $n$.

Write 
$$M'=MD^{-9\log D}.\eqno(4.15)$$
We can rewrite the sum (4.8) as follows.

\medskip\noindent
{\bf Lemma 4.2} {\it We have
$$S(M;n)=\mu (n^*)\Biggl( 
\prod_{q|n^*Z_0}\left( 1-{1\over q}\right)
\sum_{\scriptstyle 1\le d_2< n^*/M':\atop\scriptstyle d_2|n^*}
\mu (d_2)\cdot $$
$$\cdot\left(\log (n^*/d_2)-\log M-c'
-\sum_{q|n^*Z_0}{\log q\over q}\right) +$$
$$+\sum_{\scriptstyle n^*/M\le d_1\le n^*/M':\atop\scriptstyle
d_1|n^*}\mu (d_1)
\sum_{\scriptstyle 1\le k\le Md_1/n^*:\atop
\scriptstyle \gcd (k,n^*Z_0)=1}{|\mu (k)|\over\ph (k)}\Biggr)
+{\rm Error}(M;n),$$
where
$$|{\rm Error}(M;n)|\le {10^7\tau^2 (n)\log n\over D^2}.$$}
\medskip
{\bf Proof.}
By (4.8)
$$S(M;n)=\sum_{\scriptstyle d\ge 1:\gcd (Z_0,d)=1\atop\scriptstyle
d|n}\mu (d)
\sum_{\scriptstyle 1\le k\le M/d:\atop
\scriptstyle \gcd (k,nZ_0)=1}{|\mu (k)|\over\ph (k)}=$$
$$=\sum_{\scriptstyle 1\le d\le M':\atop\scriptstyle
d|n^*}\mu (d)
\sum_{\scriptstyle 1\le k\le M/d:\atop
\scriptstyle \gcd (k,n^*Z_0)=1}{|\mu (k)|\over\ph (k)}+
\sum_{\scriptstyle M'< d\le M:\atop\scriptstyle
d|n^*}\mu (d)
\sum_{\scriptstyle 1\le k\le M/d:\atop
\scriptstyle \gcd (k,n^*Z_0)=1}{|\mu (k)|\over\ph (k)}=$$
$$=
\prod_{q|n^*Z_0}\left( 1-{1\over q}\right)
\sum_{\scriptstyle 1\le d\le M':\atop\scriptstyle
d|n^*}\mu (d)\left( \log M-\log d+c'+\sum_{q|n^*Z_0}{\log q\over q}
\right)+$$
$$+\sum_{\scriptstyle M'< d\le M:\atop\scriptstyle
d|n^*}\mu (d)
\sum_{\scriptstyle 1\le k\le M/d:\atop
\scriptstyle \gcd (k,n^*Z_0)=1}{|\mu (k)|\over\ph (k)}+{\rm Error}(M;n),
\eqno(4.16)$$
where in the last step we applied Lemma 2.1, which implies the upper bound
$$|{\rm Error}(M;n)|\le 10^6\tau (n)
\sum_{\scriptstyle 1\le d\le M':\atop\scriptstyle d|n}
\left({1+\log
(M/d)\over (M/d)^{1/4}}+{\log n\over D^5}\right)
\le {10^7\tau^2 (n)\log n\over D^2},
\eqno(4.17)$$
since by (4.15), $M/d\ge D^{9\log D}$ if $1\le d\le M'$. 

The M\" obius function has the following well-known property:
$\sum_{d|k}\mu (d)=0$ for all $k\ge 2$---it
 will play a crucial role in the argument.
 The first application of this fact goes as follows.
Suppose that $n^*\ge 2$; we have with $d_1d_2=n^*$
$$\sum_{\scriptstyle 1\le d\le M':\atop\scriptstyle d|n^*}\mu (d)=
\left(\sum_{d|n^*}\mu (d)\right) -
\sum_{\scriptstyle d_1> M':\atop\scriptstyle d_1|n^*}\mu (d_1)=$$
$$=-\sum_{\scriptstyle 1\le d_2< n^*/M':\atop\scriptstyle d_2|n^*}
\mu (n^*/d_2)=-\mu (n^*)
\sum_{\scriptstyle 1\le d_2< n^*/M':\atop\scriptstyle d_2|n^*}
\mu (d_2),\eqno(4.18)$$
since
$$\sum_{d|n^*}\mu (d)=\prod_{p|n^*}(1-1)=0$$
if $n^*\ge 2$.

Similarly, if $n\ge 2$ is not a primepower, i.e., $n^*\ge 2$ is not a
prime, then
$$-\sum_{\scriptstyle 1\le d\le M':\atop\scriptstyle d|n^*}\mu (d)\log
d=-\left(\sum_{d|n^*}\mu (d)\log d\right) +
\sum_{\scriptstyle d_1> M':\atop\scriptstyle d_1|n^*}\mu (d)\log d_1=$$
$$=\sum_{\scriptstyle 1\le d_2< n^*/M':\atop\scriptstyle d_2|n^*}
\mu (n^*/d_2)\log (n^*/d_2)=\mu (n^*)
\sum_{\scriptstyle 1\le d_2< n^*/M':\atop\scriptstyle d_2|n^*}
\mu (d_2)\log (n^*/d_2),\eqno(4.19)$$
since
$$-\sum_{d|n^*}\mu (d)\log d=\sum_{p|n^*}\log p\sum_{d|{n^*\over
p}}\mu (d)=0$$
if $n^*\ge 2$ is not a prime.

We also have with $d_1=n^*/d$,
$$\sum_{\scriptstyle M'< d\le M:\atop\scriptstyle
d|n^*}\mu (d)
\sum_{\scriptstyle 1\le k\le M/d:\atop
\scriptstyle \gcd (k,n^*Z_0)=1}{|\mu (k)|\over\ph (k)}=$$
$$=\sum_{\scriptstyle M'< n^*/d_1\le M:\atop\scriptstyle
d_1|n^*}\mu (n^*/d_1)
\sum_{\scriptstyle 1\le k\le Md_1/n^*:\atop
\scriptstyle \gcd (k,n^*Z_0)=1}{|\mu (k)|\over\ph (k)}=$$
$$=\mu (n^*)\sum_{\scriptstyle n^*/M\le d_1\le n^*/M':\atop\scriptstyle
d_1|n^*}\mu (d_1)
\sum_{\scriptstyle 1\le k\le Md_1/n^*:\atop
\scriptstyle \gcd (k,n^*Z_0)=1}{|\mu (k)|\over\ph (k)}.\eqno(4.20)$$
Note that in (4.20) we have
$$1\le k\le {Md_1\over n^*}<{M\over n^*}\cdot {n^*\over M'}=
{M\over M'}=D^{9\log D}.\eqno(4.21)$$

By using (4.18), (4.19) and (4.20), we can rewrite (4.16) as follows:
$$S(M;n)=\mu (n^*)\Biggl( 
\prod_{q|n^*Z_0}\left( 1-{1\over q}\right)
\sum_{\scriptstyle 1\le d_2< n^*/M':\atop\scriptstyle d_2|n^*}
\mu (d_2)\cdot $$
$$\cdot\left(\log (n^*/d_2)-\log M-c'
-\sum_{q|n^*Z_0}{\log q\over q}\right) +$$
$$+\sum_{\scriptstyle n^*/M\le d_1\le n^*/M':\atop\scriptstyle
d_1|n^*}\mu (d_1)
\sum_{\scriptstyle 1\le k\le Md_1/n^*:\atop
\scriptstyle \gcd (k,n^*Z_0)=1}{|\mu (k)|\over\ph (k)}\Biggr)
+{\rm Error}(M;n).\eqno(4.22)$$
Combining (4.17) and (4.22), Lemma 4.2 follows.\proofend\medskip

Given an arbitrary positive integer $H$, we define a factoring 
$n=P^-_H(n)P^+_H(n)$ 
of every integer $n\ge 1$ into a product of
two of its  divisors $P^-_H(n)$ and $P^+_H(n)$ by splitting the prime
factors $p$ of $n$, counted with multiplicity, into two groups 
depending on whether $p\le H$ or $p>H$. 
More precisely, let $P^-_H(n)$ be the product of the
prime factors of $n$ which are $\le H$, and 
let $P^+_H(n)$ be the product of the prime factors of $n$
 which are $> H$ (the prime factors are taken with multiplicity, and
the empty product is 1).

We choose
$$H=H(N)=N^{1/\sqrt{\log N}}=e^{\sqrt{\log N}},\eqno(4.23')$$
$${\rm\ \ which\ \ implies\ \ } H>D^{1200\log D}\eqno(4.23'')$$
(see (2.48)). 
Note in advance that we will repeatedly use the ``almost equality of
$\chi_{-D}(P^+_H(n))$ and $\mu (P^+_H(n))$''
(see the argument after (2.10)).

By definition (see (2.48) and (4.15))
$$H=e^{\sqrt{\log N}}\ge {5\kappa ND\over 2M'},\eqno(4.24)$$
which implies
$$H\ge {n^*\over M'},$$
and so we can rewrite (4.22) in terms of $\l =P^-_H(n)$
as follows:
$$S(M;n)=\mu (n^*)\Biggl(
\prod_{q|\l Z_0}\left( 1-{1\over q}\right)
\prod_{p|n^*:p>H}\left( 1-{1\over p}\right)\cdot $$
$$\cdot\sum_{\scriptstyle 1\le d_2< n^*/M':\atop\scriptstyle d_2|\l}
\mu (d_2)\Biggl(\log (n^*/d_2)-\log M-c'
-\sum_{q|\l Z_0}{\log q\over q}
-\sum_{p|n^*:p>H}{\log p\over p}\Biggr) +$$
$$+\sum_{\scriptstyle n^*/M\le d_1\le n^*/M':\atop\scriptstyle
d_1|\l}\mu (d_1)
\sum_{\scriptstyle 1\le k\le Md_1/n^*:\atop
\scriptstyle \gcd (k,\l Z_0)=1}{|\mu (k)|\over\ph (k)}\Biggr) 
+{\rm Error}(M;n),\eqno(4.25)$$
where in the last step we used the fact
$$\sum_{\scriptstyle 1\le k\le Md_1/n^*:\atop
\scriptstyle \gcd (k,n^*D)=1}{|\mu (k)|\over\ph (k)}=
\sum_{\scriptstyle 1\le k\le Md_1/n^*:\atop
\scriptstyle \gcd (k,\l D)=1}{|\mu (k)|\over\ph (k)}$$
if $d_1\le n^*/M'$ (since $H>D^{9\log D}$; see (4.21) and (4.23)).

Motivated by (4.21) and (4.25), we call a real number
$(\l ,M)$-bad (where $\l =P^-_H(n)$) if either $x$ has the form
$$x=M'd={Md\over D^{9\log D}}{\rm\ \ for\ some\ divisor\
}d|\l,\eqno(4.26)$$
or  $x$ has the form
$$x={Md\over r}{\rm\ \ for\ some\ divisor\
}d|\l{\rm\ and\ some\ integer\ }1\le r\le D^{9\log D}.\eqno(4.27)$$
By definition, there are at most $\tau (\l )(1+D^{9\log D})$
$(\l ,M)$-bad numbers.

Let's return to (4.25): with $\l =P^-_H(n)$ we have
$$S(M;n)=S_1(M;n)+{\rm Error}_1(M;n),\eqno(4.28)$$
where
$$S_1(M;n)=\mu (n^*)\Biggl(
\prod_{q|\l Z_0}\left( 1-{1\over q}\right)
\sum_{\scriptstyle 1\le d_2< n^*/M':\atop\scriptstyle d_2|\l}
\mu (d_2)\cdot $$
$$\cdot\left(\log (n^*/d_2)-\log M-c'
-\sum_{q|\l Z_0}{\log q\over q}\right) +
\sum_{\scriptstyle n^*/M\le d_1\le n^*/M':\atop\scriptstyle
d_1|\l}\mu (d_1)
\sum_{\scriptstyle 1\le k\le Md_1/n^*:\atop
\scriptstyle \gcd (k,\l Z_0)=1}{|\mu (k)|\over\ph (k)}\Biggr)
\eqno(4.29)$$
and
$$\left|{\rm Error}_1(M;n)\right|\le
\sum_{\scriptstyle 1\le d_2< n^*/M':\atop\scriptstyle d_2|\l}
\left( \left( 1-\prod_{p|n^*:p>H}\left( 1-{1\over p}\right)\right)
\log N+
\sum_{p|n^*:p>H}{\log p\over p}\right) +$$
$$+|{\rm Error}(M;n)|.\eqno(4.30)$$
Since $n\le 3\kappa ND$  and $H=e^{\sqrt{\log N}}$, we
have
$$ 1-\prod_{p|n^*:p>H}\left( 1-{1\over p}\right)\le 
\sum_{p|n^*:p>H}{1\over p}<{\log N\over e^{\sqrt{\log N}}}.\eqno(4.31)$$
Using (4.17) and (4.31) in (4.30), we have
$$\left|{\rm Error}_1(M;n)\right|\le
\tau (n){2(\log N)^2\over e^{\sqrt{\log N}}}+
{10^7\tau^2 (n)\log n\over D^2}.\eqno(4.32)$$
Clearly
$$S^2(M;n)=S_1^2(M;n)+2{\rm Error}_1(M;n)S_1(M;n)+
{\rm Error}^2_1(M;n).\eqno(4.33)$$
Write
$${\sum}^{(1)}_{n\ne p^r}(N;M)=
\sum_{k=0}^{\kappa N}w^{(\kappa ,N)}_k
{1\over 2k+1}
\sum_{\scriptstyle (\kappa N-k)D<n< (\kappa
N+k+1)D:\atop\scriptstyle n{\rm\ is\ not\ a\ primepower}}
\chi_{-D}(n)S_1^2(M;n).\eqno(4.34)$$
Next we prove
the following lemma.
\medskip\noindent
{\bf Lemma 4.3} {\it 
We have
$$\left|{\sum}_{n\ne p^r}(N;M)-
{\sum}^{(1)}_{n\ne p^r}(N;M)\right|\le $$
$$\le {10^6D(\log N)^7\over e^{\sqrt{\log N}}}+{10^{20}(\log N)^9\over D}
+{10^{24}(\log N)^{17}\over D^3}.$$}
\medskip
{\bf Proof.}
Let's return to (4.7): 
combining (4.8), (4.14), (4.28)-(4.34), we have
$$\left|{\sum}_{n\ne p^r}(N;M)-
{\sum}^{(1)}_{n\ne p^r}(N;M)\right|\le $$
$$\le {80(\log N)^2\over e^{\sqrt{\log N}}}
\sum_{k=0}^{\kappa N}w^{(\kappa ,N)}_k
{1\over 2k+1}\sum_{\scriptstyle n:\atop \scriptstyle 
 (\kappa N-k)D< n< (\kappa N+k+1)D}\tau^2 (n)+$$
$$+{40\cdot 10^7\log M\over D^2}
\sum_{k=0}^{\kappa N}w^{(\kappa ,N)}_k
{1\over 2k+1}\sum_{\scriptstyle n:\atop \scriptstyle 
 (\kappa N-k)D< n< (\kappa N+k+1)D}\tau^3 (n)\log n+$$
$$+{80(\log N)^4\over e^{2\sqrt{\log N}}}
\sum_{k=0}^{\kappa N}w^{(\kappa ,N)}_k
{1\over 2k+1}\sum_{\scriptstyle n:\atop \scriptstyle 
 (\kappa N-k)D< n< (\kappa N+k+1)D}\tau^2 (n)+$$
$$+{2\cdot 10^{14}\over D^4}
\sum_{k=0}^{\kappa N}w^{(\kappa ,N)}_k
{1\over 2k+1}\sum_{\scriptstyle n:\atop \scriptstyle 
 (\kappa N-k)D< n< (\kappa N+k+1)D}\tau^4 (n)(\log n)^2.
\eqno(4.35)$$
We need the well-known asymptotic results
$${1\over L}\sum_{n=1}^L\tau^2 (n)=O\left( (\log L)^3\right) ,\ 
{1\over L}\sum_{n=1}^L\tau^3 (n)=O\left( (\log L)^7\right) ,$$
$$ {1\over L}\sum_{n=1}^L\tau^4 (n)=O\left( (\log L)^{15}
\right) ,{\rm\ and\ in\ general,\  }
{1\over L}\sum_{n=1}^L\tau^k (n)=O\left( (\log L)^{2^k-1}\right) 
\eqno(4.36)$$
for an arbitrary but fixed $k\ge 1$,
where the implicit constants are
effectively computable
(note that the surprising pattern 
of power-of-two minus one in
the exponents of $\log L$   is best possible).
Moreover, for later applications I mention
two more related inequalities
$$ \sum_{1\le d\le M}{\tau^2(d)\over\ph (d)}\le 
10^6(\log M)^4{\rm\ \ and\ \ }
\sum_{1\le k\le M^2}{\tau (k)\over k}\le 100 (\log
M)^2.\eqno(4.36')$$

Working with explicit constants in (4.36),
it is easy to estimate (4.35), and thus we obtain 
$$\left|{\sum}_{n\ne p^r}(N;M)-
{\sum}^{(1)}_{n\ne p^r}(N;M)\right|\le $$
$$\le {10^6D(\log N)^7\over e^{\sqrt{\log N}}}+{10^{20}
D(\log N)^9\over D^2}
+{10^{24}D(\log N)^{17}\over D^4},\eqno(4.37)$$
which completes the proof of Lemma 4.3.\proofend\medskip

I recall the  Liouville function 
$$\lambda (n)=\lambda (p_1^{r_1}p_2^{r_2}\cdots
p_s^{r_s})=(-1)^{r_1+r_2+\ldots +r_s}.$$

Applying the notation $n=P^-_H(n)P^+_H(n)$
and the Liouville function, we have
$${\sum}^{(1)}_{n\ne p^r}(N;M)=
{\sum}^{(1)}_{\lambda =\chi_{-D}}(N;M)+{\sum}^{(1)}_{\lambda \ne\chi_{-D}}(N;M),
\eqno(4.38)$$
where
$${\sum}^{(1)}_{\lambda =\chi_{-D}}(N;M)=
\sum_{k=0}^{\kappa N}w^{(\kappa ,N)}_k
{1\over 2k+1}\sum_{\l\ge 1:\ P^+_H(\l )=1}\chi_{-D}(\l)\cdot $$
$$\cdot
\sum_{\scriptstyle (\kappa N-k)D<n< (\kappa
N+k+1)D:\gcd (n,D)=1
\atop\scriptstyle P^-_H(n)=\l ,n{\rm\ is\ not\ a\ primepower}}
\lambda (P^+_H(n))
S_1^2(M;n),\eqno(4.39)$$
and
$${\sum}^{(1)}_{\lambda \ne\chi_{-D}}(N;M)=
\sum_{k=0}^{\kappa N}w^{(\kappa ,N)}_k
{1\over 2k+1}\sum_{\l\ge 1:\ P^+_H(\l )=1}\chi_{-D}(\l)\cdot $$
$$\cdot
\sum_{\scriptstyle (\kappa N-k)D<n< (\kappa
N+k+1)D:\gcd (n,D)=1
\atop\scriptstyle P^-_H(n)=\l ,
n{\rm\ is\ not\ a\ primepower}}\left(\chi_{-D} (P^+_H(n))
-\lambda (P^+_H(n))\right)
S_1^2(M;n).\eqno(4.40)$$
Note that, if $P^-_H(n)=\l\ge 2$ then
$n$ is certainly not a prime,
and $\chi_{-D}(\l)\ne 0$ implies that
the condition $\gcd (D,n)=1$  holds automatically,
since $P^+_H(\l)=1$ and $H>D$ (see (4.23)).

Next we handle $N$ as a
variable, and take the average as $N$ runs in an interval
$T\le N<2T$. That is, roughly speaking, from now on $T$ will play the
role of $N$. As a byproduct, we slightly modify (4.23) and (4.24): let
$$H=H(T)=T^{1/\sqrt{\log T}}=e^{\sqrt{\log T}},\eqno(4.41')$$
$${\rm which\ \ implies\ \ }H\ge {5\kappa TD\over M'}
{\rm\ \ and\ \ } H>D^{1200\log D}.\eqno(4.41'')$$

To finish the estimation of
$\overline{\Omega}_{1 ;\kappa ;T}(M)-
\overline{\Omega}_{-1 ;\kappa ;T}(M)$,
we need the following technical lemma.
\medskip\noindent
{\bf Lemma 4.4} {\it  We have
$${1\over T}\sum_{T\le N<2T}
\sum_{k=0}^{\kappa N}w^{(\kappa ,N)}_k
{1\over 2k+1}\sum_{\scriptstyle \l \ge 2:\ \gcd (D,\l)=1
\atop\scriptstyle P^+_H(\l)=1}\qquad\qquad
\qquad \qquad\qquad\qquad\qquad $$
$$\qquad\qquad\qquad\qquad\left|
\sum_{\scriptstyle  (\kappa N-k)D< n< (\kappa N+k+1)D:
\gcd (n,D)=1\atop \scriptstyle P^-_H(n)=\l ,\
 n{\rm\ is\ not\ a\ primepower} } 
\lambda (n)S_1^2(M;n)\right| + $$
$$+\Biggl| {1\over T}\sum_{T\le N<2T}
\sum_{k=0}^{\kappa N}w^{(\kappa ,N)}_k
{1\over 2k+1}
\sum_{\scriptstyle  (\kappa N-k)D< n< (\kappa N+k+1)D:
\gcd (n,D)=1\atop \scriptstyle P^-_H(n)=1 ,\
 n{\rm\ is\ not\ a\ primepower}} 
\lambda (n)S_1^2(M;n)  $$
$$\qquad\qquad
 -\prod_{q|Z_0}\left( 1-{1\over q}\right)^2
{D(\log T-\log M)^2\over\log T}\Biggr|\le 
2\cdot 10^3D\log D.$$}
\medskip
The proof of Lemma 4.4 is elementary but technical. We just mention the key
ingredient---a routine sieve lemma---at the beginning 
of the next section (see Lemma 5.1), and
postpone the details of the proof of Lemma 4.4 to Sections 7-8.
\medskip\noindent
{\bf Concluding Remark of Section 4.} Notice that Lemmas 4.1 and 4.4
do {\it not} distinguish between the primes
$q|Z_0$ with $\chi_{-D}(q)=1$, $0$ or $-1$.

\bigskip\bigskip\bigskip
\medskip
\centerline{\bf 5. The method of the proof (II)}
\bigskip\bigskip\bigskip\medskip\noindent
Beside the simple but crucial fact $\sum_{d|k}\mu (d)=0$
for all $k\ge 2$,
the proof of Lemma 4.4 is based
on the next lemma, which is a routine application of the simplest
``sieve method'' in number theory.
\medskip\noindent
{\bf Lemma 5.1} {\it If $\log x>2e^2\log H(\log\log H+10)$ then
we can estimate the following Liouville sum from above:
$$\left| \sum_{\scriptstyle 1\le n \le x
\atop\scriptstyle P^-_H(n)=1 }\lambda (n)\right|\le
x\exp\left(-{\log x\over 2\log H}+1\right)+
x\log x\exp\left( -\sqrt{\log x}/15\right)+{x\over H},$$
where $\exp (y)=e^y$.}
\medskip
{\bf Proof.} If $\lambda (n)\ne\mu (n)$ and
$P^-_H(n)=1$, then $n$ can be written in the form $n=p^2s$ with some
prime $p>H$. Thus we have 
$$\left| \sum_{\scriptstyle 1\le n \le x
\atop\scriptstyle P^-_H(n)=1 }\lambda (n)-
\sum_{\scriptstyle 1\le n \le x
\atop\scriptstyle P^-_H(n)=1 }\mu (n)\right|\le 
\sum_{p>H}{x\over p^2}<\sum_{m>H}{x\over m^2}
<{x\over H}.\eqno(5.1)$$
To study the restricted M\"obius sum in (5.1), we use the 
inclusion-exclusion principle:
$$\sum_{\scriptstyle 1\le n \le x
\atop\scriptstyle P^-_H(n)=1 }\mu (n)=
\sum_{\scriptstyle 1\le m \le x
\atop\scriptstyle P^+_H(m)=1 }\mu^2 (m)
\sum_{1\le r\le x/m}\mu (r).\eqno(5.2) $$
We need the following well-known result that I put in the form of 
another lemma.
\medskip\noindent
{\bf Lemma 5.2} {\it The M\"obius sum
$${\cal M}(L)=\sum_{1\le n\le L}\mu (n)$$
has the following upper bound for all $L\ge 2$:
$$|{\cal M}(L)|\le Le^{-\sqrt{\log L}/10}.\eqno(5.3)$$}
\medskip
{\bf Remark.} It is worth to point out
that Lemma 5.2 is a deep result in analytic number
theory (see e.g. [Ka]). It is based on the fact that the Dirichlet series
$$\sum_{n=1}^{\infty}{\mu (n)\over n^s}=\prod_p\left( 1-{1\over
p^s}\right)={1\over\zeta (s)}$$
is the reciprocal of  the Riemann's zeta function
$\zeta (s)$, and so the proof
of (5.3) can be carried out along the same lines as that of
the classical analytic proof of the Prime Number
Theorem, due to Hadamard and de la Vallee Poussin.
\medskip
Returning to the proof of Lemma 5.1, by Lemma 5.2 we have
$$\left|\sum_{1\le r\le x/m}\mu (r)\right|\le
{x\over m}e^{-\sqrt{\log\sqrt{x}}/10}
{\rm\ \ if\ \ }1\le m\le\sqrt{x}.\eqno(5.4) $$
If $\sqrt{x}<m\le x$ then we just use the trivial bound
$$\left|\sum_{1\le r\le x/m}\mu (r)\right|\le
{x\over m}.\eqno(5.5) $$
Applying (5.4)-(5.5) in (5.2), we have
$$\left|\sum_{\scriptstyle 1\le n \le x
\atop\scriptstyle P^-_H(n)=1 }\mu (n)\right|\le 
\sum_{1\le m \le \sqrt{x}}{x\over m}e^{-\sqrt{\log x}/10\sqrt{2}}+$$
$$+\sum_{\scriptstyle \sqrt{x}<m\le x
\atop\scriptstyle P^+_H(m)=1 }{x\mu^2 (m)\over m}
\le {x\log x\over e^{\sqrt{\log x}/15}}+x
\sum_{\scriptstyle \sqrt{x}<m\le x
\atop\scriptstyle P^+_H(m)=1 }{\mu^2 (m)\over m}.\eqno(5.6)$$
To estimate the last sum in (5.6), we use the obvious fact
that for every integer $k\ge 1$,
$$\sum_{\scriptstyle m>H^k
\atop\scriptstyle P^+_H(m)=1 }{\mu^2 (m)\over m}
\le {1\over (k+1)!}\left(\sum_{1<p\le H}{1\over p}\right)^{k+1}+$$
$$+{1\over (k+2)!}\left(\sum_{1<p\le H}{1\over p}\right)^{k+2}+
{1\over (k+3)!}\left(\sum_{1<p\le H}{1\over p}\right)^{k+3}+\ldots
.\eqno(5.7)$$
It is well-known that
$$\sum_{p\le H}{1\over p}\le\log\log H +10,$$
so if $k\ge e^2 (\log\log H+10)$,
then for every $j>k$ we have
$${1\over j!}\left(\sum_{1<p\le H}{1\over p}\right)^j\le
\left({e\over j}\sum_{p\le H}{1\over p}\right)^j\le
\left({e(\log\log H+10)\over j}\right)^j\le
e^{-j}.$$
Using this in (5.7), we have
$$\sum_{\scriptstyle m>H^k
\atop\scriptstyle P^+_H(m)=1 }{\mu^2 (m)\over m}
\le \sum_{j>k}e^{-j}<e^{-k}{\rm\ \ 
for\ \ } k\ge e^2 (\log\log H+10).\eqno(5.8)$$
Note that
$$\sqrt{x}\ge H^k{\rm\ \ for\ \ }k=\left\lfloor {\log x\over 2\log
H}\right\rfloor .\eqno(5.9)$$
Combining (5.1), (5.6) and (5.8)-(5.9), Lemma 5.1 follows.
\proofend\medskip
The deduction of Lemma 4.4 from Lemma 5.1 is a routine but rather
long and cumbersome estimation. We postpone it to Sections 7-8.

Let's return to (4.39): write
$$\overline{\sum}^{(1)}_{\lambda =\chi_{-D}}(T;M)=
{1\over T}\sum_{T\le N<2T}
{\sum}^{(1)}_{\lambda =\chi_{-D}}(N;M)=$$
$$={1\over T}\sum_{T\le N<2T}
\sum_{k=0}^{\kappa N}w^{(\kappa ,N)}_k
{1\over 2k+1}\sum_{\scriptstyle (\kappa N-k)D<n< (\kappa
N+k+1)D:\gcd (n,D)=1
\atop\scriptstyle P^-_H(n)=1 ,n{\rm\ is\ not\ a\ primepower}}
\lambda (n)S_1^2(M;n)+$$
$$+{1\over T}\sum_{T\le N<2T}
\sum_{k=0}^{\kappa N}w^{(\kappa ,N)}_k
{1\over 2k+1}\sum_{\l\ge 2:\ P^+_H(\l )=1}{\chi_{-D}(\l)\over\lambda
(\l)}\cdot $$
$$\cdot
\sum_{\scriptstyle (\kappa N-k)D<n< (\kappa
N+k+1)D:\gcd (n,D)=1
\atop\scriptstyle P^-_H(n)=\l ,n{\rm\ is\ not\ a\ primepower}}
\lambda (n)S_1^2(M;n).\eqno(5.10)$$
Combining (5.10) with Lemma 4.4, we obtain, 
$$\left|\overline{\sum}^{(1)}_{\lambda =\chi_{-D}}(T;M)
 -\prod_{q|Z_0}\left( 1-{1\over q}\right)^2
{D(\log T-\log M)^2\over\log T}\right|\le $$
$$\le 2\cdot 10^3D\log D.\eqno(5.11)$$
Next we go to (4.40): write
$$\overline{\sum}^{(1)}_{\lambda \ne\chi_{-D}}(T;M)=
{1\over T}\sum_{T\le N<2T}
{\sum}^{(1)}_{\lambda\ne\chi_{-D}}(N;M)=$$
$$={1\over T}\sum_{T\le N<2T}
\sum_{k=0}^{\kappa N}w^{(\kappa ,N)}_k
{1\over 2k+1}\sum_{\l\ge 1:\ P^+_H(\l )=1}\chi_{-D}(\l)\cdot $$
$$\cdot
\sum_{\scriptstyle (\kappa N-k)D<n< (\kappa
N+k+1)D:\gcd (n,D)=1
\atop\scriptstyle P^-_H(n)=\l ,
n{\rm\ is\ not\ a\ primepower}}\left(\chi_{-D} (P^+_H(n))
-\lambda (P^+_H(n))\right)
S_1^2(M;n).\eqno(5.12)$$
The following lemma involves $L(1,\chi_{-D})$, which is small by hypothesis.
\medskip\noindent
{\bf Lemma 5.3} {\it We have 
$$\left|\overline{\sum}^{(1)}_{\lambda \ne\chi_{-D}}(T;M)\right| \le 
 10^9(\log T)^5D\sum_{\scriptstyle H<q\le T^2:
\atop\scriptstyle \chi_{-D}(q)=1}{1\over q}\le $$
$$\le 10^{11}(\log T)^6D\cdot L(1,\chi_{-D}).$$}
{\bf Proof.}
Let $J(\kappa N;k)$ denote the interval $(\kappa N-k)D<x<
(\kappa N+k+1)D$. 

By (5.12),
$$\left|\overline{\sum}^{(1)}_{\lambda \ne\chi_{-D}}(T;M)\right| \le $$
$$\le {2\over T}\sum_{T\le N<2T}
\sum_{k=0}^{\kappa N}w^{(\kappa ,N)}_k
{1\over 2k+1}
\sum_{\scriptstyle n\in J(\kappa N;k):\ \gcd (D,n)=1
\atop\scriptstyle \chi_{-D}(P^+_H(n))\ne\lambda (P^+_H(n))}
S_1^2(M;n).\eqno(5.13)$$

If $\chi_{-D}(P^+_H(n))\ne\lambda (P^+_H(n))$ and
$\gcd (D,n)=1$, then there is a prime $q>H$ with
$\chi_{-D}(q)=1$ such that $q|n$. Using this observation, and 
the trivial upper bound 
$$S_1^2(M;n)\le 400\tau^2 (n)(\log M)^2\eqno(5.14)$$
in (5.13), we have
$$\left|\overline{\sum}^{(1)}_{\lambda \ne\chi_{-D}}(T;M)\right| \le $$
$$\le 400(\log M)^2
{2\over T}\sum_{T\le N<2T}
\sum_{k=0}^{\kappa N}w^{(\kappa ,N)}_k
{1\over 2k+1}
\sum_{\scriptstyle n\in J(\kappa N;k):\ \gcd (D,n)=1
\atop\scriptstyle \chi_{-D}(P^+_H(n))\ne\lambda (P^+_H(n))}
\tau^2 (n)\le $$
$$\le 400(\log M)^2\left( {2\over\kappa T}\sum_{\scriptstyle 1\le n\le 5\kappa
TD:\ \gcd (D,n)=1
\atop\scriptstyle \chi_{-D}(P^+_H(n))\ne\lambda (P^+_H(n))}\tau^2 (n)\right)
{1\over T}\sum_{T\le N<2T}
\sum_{k=0}^{\kappa N}w^{(\kappa ,N)}_k=$$
$$={800(\log M)^2\over\kappa T}\sum_{\scriptstyle 1\le n\le 5\kappa
TD:\ \gcd (D,n)=1
\atop\scriptstyle \chi_{-D}(P^+_H(n))\ne\lambda (P^+_H(n))}\tau^2 (n)\le $$
$$\le {800(\log M)^2\over\kappa T}\sum_{\scriptstyle H<
 q\le 5\kappa TD:
\atop\scriptstyle \chi_{-D}(q)=1}
\sum_{ 1\le m\le 5\kappa TD/q}\tau^2 (mq)\le $$
$$\le {800(\log M)^2\over\kappa T}\sum_{\scriptstyle q>H:
\atop\scriptstyle \chi_{-D}(q)=1}\tau^2 (q)
\sum_{ 1\le m\le 5\kappa TD/q}\tau^2 (m)\le $$
$$\le {4000(\log M)^2\over\kappa T}\sum_{\scriptstyle H<q\le 5\kappa TD:
\atop\scriptstyle \chi_{-D}(q)=1}
10^3{\kappa TD\over q}(\log (T^2))^3\le
 10^9(\log T)^5D\sum_{\scriptstyle H<q\le T^2:
\atop\scriptstyle \chi_{-D}(q)=1}{1\over q}.\eqno(5.15)$$
Applying Lemma 1.1 in (5.15), we easily have
$$\left|\overline{\sum}^{(1)}_{\lambda \ne\chi_{-D}}(T;M)\right| \le 
 10^9(\log T)^5D\sum_{\scriptstyle H<q\le T^2:
\atop\scriptstyle \chi_{-D}(q)=1}{1\over q}\le $$
$$\le 10^9(\log T)^5D\cdot 100\log T\cdot L(1,\chi_{-D})=
10^{11}(\log T)^6D\cdot L(1,\chi_{-D}),\eqno(5.16)$$
completing the proof of Lemma 5.3.\proofend\medskip

Let's return to (4.38): write
$${\sum}^{(1)}_{n\ne p^r}(T;M)={1\over T}\sum_{T\le N <2T}
{\sum}^{(1)}_{n\ne p^r}(N;M).$$
Combining (4.38)-(4.40), (5.10)-(5.12) and Lemma 5.3, we have
$$\left|{\sum}^{(1)}_{n\ne p^r}(T;M)
 -\prod_{q|Z_0}\left( 1-{1\over q}\right)^2
{D(\log T-\log M)^2\over\log T}\right|\le $$
$$\le 2\cdot 10^3D\log D+
10^{11}(\log T)^6D\cdot L(1,\chi_{-D}).\eqno(5.17)$$
Next write
$${\sum}_{n\ne p^r}(T;M)={1\over T}\sum_{T\le N <2T}
{\sum}_{n\ne p^r}(N;M).$$
By (5.17) and Lemma 4.3,
$$\left|{\sum}_{n\ne p^r}(T;M)
 -\prod_{q|Z_0}\left( 1-{1\over q}\right)^2
{D(\log T-\log M)^2\over\log T}\right|\le $$
$$\le 3\cdot 10^3D\log D.\eqno(5.18)$$
Now we are ready to estimate the difference
$$\overline{\Omega}_{1 ;\kappa ;T}(M)-
\overline{\Omega}_{-1 ;\kappa ;T}(M)$$
(see (2.44)) as we promised at the beginning of Section 4.
Combining (5.18) with (4.5) and Lemma 4.1, we obtain
$$\left|\overline{\Omega}_{1 ;\kappa ;T}(M)-
\overline{\Omega}_{-1 ;\kappa ;T}(M)-
\prod_{q|Z_0}\left( 1-{1\over q}\right)^2
D(\log T-2\log M)\right|\le $$
$$\le 3\cdot 10^3D\log D+30D\log D+2\le $$
$$\le 4\cdot 10^3D\log D,\eqno(5.19)$$
where  the main term comes from
$${(\log T-\log M)^2-(\log M)^2\over\log T}=\log T-2\log M.$$

To emphasize its importance, we rewrite (5.19) as a lemma.
\medskip\noindent
{\bf Lemma 5.4} {\it We have
$$\left|\overline{\Omega}_{1 ;\kappa ;T}(M)-
\overline{\Omega}_{-1 ;\kappa ;T}(M)-
\prod_{q|Z_0}\left( 1-{1\over q}\right)^2
D(\log T-2\log M)\right|\le $$
$$\le 4\cdot 10^3D\log D.$$}
\proofend\medskip

The next lemma will  lead to a {\it different} approximation of 
$$\left|\overline{\Omega}_{1 ;\kappa ;T}(M)-.
\overline{\Omega}_{-1 ;\kappa ;T}(M)\right| .$$
\medskip\noindent
{\bf Lemma 5.5} {\it We have (see (3.20))
$$\Biggl| OD_{\kappa ;N}(M)
+{D\over 2}\log {M^2\over N}\cdot $$
$$\cdot
\prod_{\scriptstyle q|Z_0\atop\scriptstyle \chi_{-D}(q)=1}
\left(1-{1\over q}\right)
\prod_{\scriptstyle q|Z_0\atop\scriptstyle \chi_{-D}(q)=0}
\left(1-{1\over q}\right)^2
\prod_{\scriptstyle q|Z_0\atop\scriptstyle \chi_{-D}(q)=-1}
\left(1-{1\over q}\right)^2\left(1-{1\over q+1}\right)
\Biggr|\le $$
$$\le 10^{22}D(\log D)^{11}+\log {M^2\over N}\cdot
{10^5D\over (\log D)^3}.$$}
\medskip
Unfortunately the proof of Lemma 5.5 is very long
(starting from Sections 9, it is the subject of 
 more than a dozen sections).

Note that Lemma 5.4 represents the evaluation of the left side of
(2.49) via (2.40) and (2.44), 
and Lemma 5.5 represents the evaluation of the dominating part of
the right side of
(2.49) via (2.42) and (2.45). 
\medskip\noindent
{\bf Completing the proof of Theorem 1 via contradiction.}
We recall (3.19):
$$\Omega_{1 ;\kappa ;N}(M)
-\Omega_{-1;\kappa ;N}(M)=$$
$$=\Omega_{1 ;\kappa ;N}({\rm OffDiag};M)
-\Omega_{-1;\kappa ;N}({\rm OffDiag};M)=$$
$$=2OD_{\kappa ;N}(M)+
{\rm Error}(\kappa ;N;M).$$
By Lemma 3.4,
$$\left|{\rm Error}(\kappa ;N;M)\right| <1 .$$
Thus we have
$$\Biggl|\Omega_{1 ;\kappa ;N}({\rm OffDiag};M)-
\Omega_{-1 ;\kappa ;N}({\rm OffDiag};M)-
2OD_{\kappa ;N}(M)\Biggr| <1.\eqno(5.20)$$

We want to estimate 
$$\overline{\Omega}_{1 ;\kappa ;T}({\rm OffDiag};M)-
\overline{\Omega}_{-1 ;\kappa ;T}({\rm OffDiag};M)$$
where $\overline{\Omega}_{\pm 1 ;\kappa ;T}({\rm OffDiag};M)$ is
the average 
$$\overline{\Omega}_{\pm 1 ;\kappa ;T}({\rm OffDiag};M)=
{1\over T}\sum_{T\le N<2T}
{\Omega}_{\pm 1 ;\kappa ;N}({\rm OffDiag};M).$$
Since
$$\log {M^2\over N}=2\log M-\log N,$$
by Lemma 5.5 and (5.20)  we have
$$\Biggl|\overline{\Omega}_{1 ;\kappa ;T}(M)
-\overline{\Omega}_{-1 ;\kappa ;T}(M)-D(\log T-2\log M)\cdot $$
$$\cdot
\prod_{\scriptstyle q|Z_0\atop\scriptstyle \chi_{-D}(q)=1}
\left(1-{1\over q}\right)
\prod_{\scriptstyle q|Z_0\atop\scriptstyle \chi_{-D}(q)=0}
\left(1-{1\over q}\right)^2
\prod_{\scriptstyle q|Z_0\atop\scriptstyle \chi_{-D}(q)=-1}
\left(1-{1\over q}\right)^2\left(1-{1\over q+1}\right) \Biggr|\le $$
$$\le 2\cdot 10^{22}D(\log D)^{11}+\log {M^2\over N}\cdot
{2\cdot 10^5D\over (\log D)^3}+2D+1.\eqno(5.21)$$
Subtracting (5.21) from Lemma 5.4,
$\overline{\Omega}_{\pm 1 ;\kappa ;T}(M)$ cancels out, and we have
$$\left|\prod_{p|Z_0}\left(1-{1\over p}\right)^2
\cdot D\log {M^2\over T}
\left( 1-
\prod_{\scriptstyle p|Z_0\atop\scriptstyle \chi_{-D}(p)=1}
\left( 1+{1\over p-1}\right)
\prod_{\scriptstyle p|Z_0\atop\scriptstyle \chi_{-D}(p)=-1}
\left( 1-{1\over p+1}\right)\right) \right|\le $$
$$\le  10^4D(\log D)^2
+2\left( 10^{22}D(\log D)^{11}+\log {M^2\over N}\cdot
{10^5D\over (\log D)^3}\right) +2D+1.\eqno(5.22)$$

According to the definition of $Z_0$ (see (1.14)-(1.15)), 
we now distinguish
two cases ($\chi_{-D}(p_0)=1$ or $-1$).
If $\chi_{-D}(p_0)=1$, then 
$$Z_0=\prod_{\scriptstyle p\le D:\atop\scriptstyle
\chi_{-D}(p)\ne -1}p,$$
and thus we have the  estimation
$$\prod_{p|Z_0}\left(1-{1\over p}\right)^2
\left( 1-
\prod_{\scriptstyle p|Z_0\atop\scriptstyle \chi_{-D}(p)=1}
\left( 1+{1\over p-1}\right)
\prod_{\scriptstyle p|Z_0\atop\scriptstyle \chi_{-D}(p)=-1}
\left( 1-{1\over p+1}\right)\right)=$$
$$=\prod_{p|Z_0}\left(1-{1\over p}\right)^2
\left( 1-
\prod_{\scriptstyle p|Z_0\atop\scriptstyle \chi_{-D}(p)=1}
\left( 1+{1\over p-1}\right)\right) =$$
$$= \prod_{p|Z_0}\left(1-{1\over p}\right)^2
\left( 1-\left( 1+{1\over p_0-1}\right)
\cdot
\prod_{\scriptstyle p|{Z_0\over p_0}\atop\scriptstyle \chi_{-D}(p)=1}
\left( 1+{1\over p-1}\right)\right) \le $$
$$\le - \prod_{p|Z_0}\left(1-{1\over p}\right)^2
{1\over p_0-1}
\prod_{\scriptstyle p|{Z_0\over p_0}\atop\scriptstyle \chi_{-D}(p)=1}
\left( 1+{1\over p-1}\right) =$$
$$=-{1\over p_0-1}\left(1-{1\over p_0}\right)^2
\prod_{\scriptstyle p|{Z_0\over p_0D}\atop\scriptstyle \chi_{-D}(p)=1}
\left( 1-{1\over p}\right) 
 \prod_{p|D}\left(1-{1\over p}\right)^2\le $$
$$\le -{1\over 2+2\log D}
\prod_{\scriptstyle p|{Z_0\over p_0D}\atop\scriptstyle \chi_{-D}(p)=1}
\left( 1-{1\over p}\right) \cdot
{\ph^2(D)\over D^2}\le $$
$$\le -{1\over 2+2\log D}\cdot {1\over 400\log D(\log\log D)^2},
\eqno(5.23)$$
since
$$\left(1-{1\over p}\right)\left( 1+{1\over p-1}\right)=1,$$
and
$$\prod_{q|Z_0}\left(1-{1\over q}\right)
\ge \prod_{p\le D}\left(1-{1\over q}\right)
\ge {1\over 4\log D},$$
where we used the well-known number-theoretic facts
${n\over\ph (n)}\le 10\log\log n$ (see (3.28))
and
$$\prod_{p\le n}\left( 1-{1\over p}\right)
=(1+o(1)){e^{-\gamma_0}\over\log n},$$
and also used that $D>e^{10^{100}}$ is very large.

If $\chi_{-D}(p_0)=-1$ then we distinguish three subcases. If
$$\prod_{\scriptstyle p|Z_0\atop\scriptstyle \chi_{-D}(p)=1}
\left( 1+{1\over p-1}\right)\ge 2,$$
then again we choose
$$Z_0=\prod_{\scriptstyle p\le D:\atop\scriptstyle
\chi_{-D}(p)\ne -1}p,$$
and  have the estimation
$$\prod_{p|Z_0}\left(1-{1\over p}\right)^2\left( 1-
\prod_{\scriptstyle p|Z_0\atop\scriptstyle \chi_{-D}(p)=1}
\left( 1+{1\over p-1}\right)
\prod_{\scriptstyle p|Z_0\atop\scriptstyle \chi_{-D}(p)=-1}
\left( 1-{1\over p+1}\right)\right) =$$
$$=\prod_{p|Z_0}\left(1-{1\over p}\right)^2\left( 1-
\prod_{\scriptstyle p|Z_0\atop\scriptstyle \chi_{-D}(p)=1}
\left( 1+{1\over p-1}\right)\right) \le $$
$$\le -\prod_{p|Z_0}\left(1-{1\over p}\right)^2\cdot {1\over 2}
\prod_{\scriptstyle p|Z_0\atop\scriptstyle \chi_{-D}(p)=1}
\left( 1+{1\over p-1}\right) =$$
$$=-{1\over 2}\prod_{\scriptstyle 
p|{Z_0\over D}\atop\scriptstyle \chi_{-D}(p)=1}
\left( 1-{1\over p}\right) 
 \prod_{p|D}\left(1-{1\over p}\right)^2=$$
$$= -{1\over 2}
\prod_{\scriptstyle p|{Z_0\over D}\atop\scriptstyle \chi_{-D}(p)=1}
\left( 1-{1\over p}\right) \cdot
{\ph^2(D)\over D^2}\le $$
$$\le -{1\over 2}\cdot {1\over 400\log D(\log\log D)^2},
\eqno(5.24)$$
where in the last steps we repeated the estimations in (5.23).

If
$$2>\prod_{\scriptstyle p|Z_0\atop\scriptstyle \chi_{-D}(p)=1}
\left( 1+{1\over p-1}\right)\ge 1+{1\over 4\log D},$$
then again we choose
$$Z_0=\prod_{\scriptstyle p\le D:\atop\scriptstyle
\chi_{-D}(p)\ne -1}p,$$
and  have the estimation
$$\prod_{p|Z_0}\left(1-{1\over p}\right)^2\left( 1-
\prod_{\scriptstyle p|Z_0\atop\scriptstyle \chi_{-D}(p)=1}
\left( 1+{1\over p-1}\right)
\prod_{\scriptstyle p|Z_0\atop\scriptstyle \chi_{-D}(p)=-1}
\left( 1-{1\over p+1}\right)\right) =$$
$$=\prod_{p|Z_0}\left(1-{1\over p}\right)^2\left( 1-
\prod_{\scriptstyle p|Z_0\atop\scriptstyle \chi_{-D}(p)=1}
\left( 1+{1\over p-1}\right) \right)\le $$
$$\le \prod_{p|Z_0}\left(1-{1\over p}\right)^2
\left( 1-\left( 1+{1\over 4\log D}\right)\right) =$$
$$= - \prod_{p|Z_0}\left(1-{1\over p}\right)^2
\cdot {1\over 4\log D}=$$
$$=-\prod_{\scriptstyle p|{Z_0\over D}\atop\scriptstyle \chi_{-D}(p)=1}
\left( 1-{1\over p}\right)^2\cdot
{\ph^2(D)\over D^2}\cdot {1\over 4\log D}\le $$
$$\le -{1\over 2^2}\cdot {1\over 100(\log\log D)^2}
\cdot {1\over 4\log D},\eqno(5.25)$$
since in this subcase
$$\prod_{\scriptstyle p|Z_0\atop\scriptstyle \chi_{-D}(p)=1}
\left( 1-{1\over p}\right) =
\prod_{\scriptstyle p|Z_0\atop\scriptstyle \chi_{-D}(p)=1}
\left( 1+{1\over p-1}\right)^{-1}\ge {1\over 2}.$$

Finally, 
if $\chi_{-D}(p_0)=-1$ and
$$\prod_{\scriptstyle p|Z_0\atop\scriptstyle \chi_{-D}(p)=1}
\left( 1+{1\over p-1}\right)< 1+{1\over 4\log D},$$
then we choose
$$Z_0=p_0\prod_{\scriptstyle p\le D:\atop\scriptstyle
\chi_{-D}(p)\ne -1}p,$$
and  have the estimation
$$\prod_{p|Z_0}\left(1-{1\over p}\right)^2
\left( 1-
\prod_{\scriptstyle p|Z_0\atop\scriptstyle \chi_{-D}(p)=1}
\left( 1+{1\over p-1}\right)\cdot
\left( 1-{1\over p_0+1}\right)\right) \ge $$
$$\ge \prod_{p|Z_0}\left(1-{1\over p}\right)^2
\left( 1-\left( 1+{1\over 4\log D}\right)
\left( 1-{1\over p_0+1}\right)\right) \ge $$
$$\ge \prod_{p|Z_0}\left(1-{1\over p}\right)^2
\cdot {1\over 4\log D}\ge $$
$$\ge {1\over 2^2}\cdot {1\over 100(\log\log D)^2}
\cdot {1\over 4\log D},\eqno(5.26)$$
since $p_0<2\log D$ (and we also repeated the estimations in (5.25)).

Combining (5.22)-(5.26), we obtain
$${1\over 10^3(\log D)^2(\log\log D)^2}\log {M^2\over T}\le $$
$$\le  10^4(\log D)^2+
2\left( 10^{22}(\log D)^{11}+
{10^5\log M\over (\log D)^3}\right) +3
.\eqno(5.27)$$
Using $\log T=(\log D)^{15}$,  
$M=T\exp\left(-{2\over 3}\sqrt{\log T}\right) $,
and $\log D>10^{100}$, we see that the left-hand side of (5.27) is
clearly {\it larger} than
the right-hand side (since
$(\log D\log\log D)^2$ is much smaller than $(\log D)^3$),
which is a contradiction.
This contradiction proves Theorem 1.

It remains to prove Lemma 2.1 (see Section 30), 
 Lemma 4.4 (see Sections 7-8), and Lemma 5.5
(it starts from Section 9).

The proof of Lemma 2.1 is based on the well-known elementary fact
$$\left|\sum_{1\le m\le x}{1\over m}-\log x-\gamma_0\right|\le
{5\over x},$$
where $\gamma_0=0.5772\ldots $ is Euler's constant.
The proof of Lemma 4.4 is based on the ``sieve lemma'' Lemma 5.1 
(combined with Lemma 5.2) and the Prime Number Theorem.

Finally, the long proof of Lemma 5.5 uses the following key ingredients:

Fact (1): the P\' olya--Vinogradov inequality for character sums;

Fact (2): the  sum of prime-reciprocals starting from $p>D$
and ending at $x$
$$\sum_{\scriptstyle D<p\le x\atop
\scriptstyle\chi_{-D}(p)=1}{1\over p}$$
is ``small'' if the class number $h(-D)$ is substantially smaller than
$\sqrt{D}$ and $x$ is not too large; 

Fact (3): the proof technique
of the well-known average result for the divisor function
$$\sum_{1\le n\le x}\tau (n)=x\log x+(2\gamma_0-1)x+O(\sqrt{x}),$$
and other similar average type arguments/estimations based
on partial summation;

Fact (4): let
$$S(\Delta ;k)=\sum_{1\le n\le k}\chi_{\Delta}(n),$$
then for every negative fundamental 
discriminant $\Delta <0$,
$${1\over |\Delta |}\sum_{k=1}^{|\Delta |-1}
S(\Delta ;k)=h(\Delta )$$
(note that Fact (4) is a corollary of the class number formula (1.4));

Fact (5): the ``logarithmic'' character sum
$\sum_{j=1}^{D^4}\chi_{-D}(j)\log j$ can be estimated as follows:
$$\left| \sum_{j=1}^{D^4}\chi_{-D}(j)\log j+{\pi\over 6}\sqrt{D}
\sum_{(a,b,c)}{1\over a}\right| \le $$
$$\le  h(-D)(6\log D+30)+
 21(\log D)^2 D^{1/6}+ \sqrt{D}\left(
10^3\left( {h(-D)\over \sqrt{D}}\right)^{3/2}
+{10^3h(-D)\over \sqrt{D}}\right) $$
(see Lemma 14.2),
where $\sum_{(a,b,c)}{1\over a}$ means that
we add up the reciprocals of the leading coefficients 
$a=a_j$ in the family of
reduced, primitive, inequivalent binary quadratic forms of integer
coefficients $ax^2+bxy+cy^2$ with discriminant
$b^2-4ac=-D<0$,  $a>0$, $c>0$, $1\le j\le h(-D)$
(note that Fact (5) is based on 
$$\sum_{j=1}^{D^4}\chi_{-D}(j)\log j=-{\sqrt{D}\over \pi}
L'(1,\chi_{-D})+{\rm\ negligible}$$
if $h(-D)$ is small, where $L'(1,\chi_{-D})$ is
the derivative of the Dirichlet L-function $L(s,\chi_{-D})$ at $s=1$);

and to estimate the critical sum
$\sum_{(a,b,c)}{1\over a}$ in Fact (5), we need

Fact (6): 
$$0\le\prod_{\scriptstyle p|Z_0\atop\scriptstyle \chi_{-D}(p)\ne -1}
{p+1\over p-1}\prod_{p|D}{p-1\over p}
-\sum_{j=1}^{h(-D)}{1\over a_j}\le $$
$$\le {10^3h(-D)(\log D)^4\over\sqrt{D}}+{10^4\over (\log D)^3}+
{1\over D}$$
(see Lemma 14.3).

Facts (5) and (6) enter the proof of Lemma 5.5 in the following
way. The sum $OD_{\kappa ;N}(M)$ (see (3.20))
resembles a ``double harmonic sum'', and also the main term in Lemma
2.1 is the logarithm function. This explain why (3.20) resembles a
Riemann sum for a logarithmic integral that we can evaluate
explicitly. Thus we obtain---via long
estimations---the following good approximation of (3.20) 
(see (19.12) later):
$$OD_{\kappa ;N}(M)=
OD_{\kappa ;N}(M;{\rm Core};\clubsuit\clubsuit
\clubsuit ;\eps )+{\rm negligible},$$
where 
$$OD_{\kappa ;N}(M;{\rm Core};\clubsuit\clubsuit
\clubsuit ;\eps )=$$
$$=\sqrt{D}
\sum_{\scriptstyle 1\le d\le M:\atop\scriptstyle \gcd (d,Z_0)=1}
{\mu^2(d)\over\ph (d)}
\sum_{\scriptstyle 1\le r\le D^8:\atop\scriptstyle \gcd (r,dZ_0)=1}
{\mu(r)\chi_{-D}(r)\over r\ph (r)}
\sum_{\scriptstyle (r_1,r_2):\atop\scriptstyle r_1r_2=r}\cdot$$
$$\cdot\sum_{d_1|d}\chi_{-D}(d_1)\prod_{p|d_1}{p-2\over p-1}
\sum_{v|rd_1}\mu (v)\chi_{-D}(v)\sum_{j=1}^{D^4}\chi_{-D}(j)
\sum_{\l\in\Z}\eps\cdot f\left( (1+\eps )^{\l}\right)\cdot $$
$$\cdot\Biggl(
c_{13}(rd)\left( \log\left( {M^2d_1(1+\eps )^{\l }\over 2\pi Nd^2vj}
\right)\right)^3+
c_{14}(rd)
\left( \log\left( {M^2d_1(1+\eps )^{\l }\over 2\pi Nd^2vj}
\right)\right)^2+\qquad\qquad\qquad $$
$$ +c_{15}(rd) 
\left( \log\left( {M^2d_1(1+\eps )^{\l }\over 2\pi Nd^2vj}
\right)\right)
+c_{16}(rd)\Biggr) \cdot
\delta_{1,0}\left\{  {M^2d_1\over 2\pi Nd^2v}\ge C^{\star}\right\}
,\eqno(5.28)$$
where (see (9.7))
$$f(y)=
\sin (\kappa y)\left( {\sin ((1+{1\over 2N})y)\over
(2N+1)\sin (y/2N)}\right)^{\kappa }$$
if $|y|\le\pi N$, and $f(y)=0$ if $|y|>\pi N$,
and
$$\sum_{\l\in\Z}\eps\cdot f\left( (1+\eps )^{\l}\right)
={\pi\over 2}+{\rm negligible}$$
(see (12.7), (12.12), (12.15) and Lemma 14.1); moreover
$c_i(rd)$, $13\le i\le 16$ are constants depending only on 
$rd$ (see (19.4)-(19.5)), $C^{\star}=e^{(\log D)^3}$, and
$\delta_{1,0}\{\cdots \}$ is a 0-1 valued ``cutoff function''
defined as
$$\delta_{1,0}\{ {\rm true }\} =1 {\rm\ \ and\ \ }
\delta_{1,0}\{ {\rm false }\} =0.$$
Equation (5.28) is the first main step in the long proof of Lemma 5.5.

For the second main step, we
use the explicit forms of the constants
$c_i(rd)$, $13\le i\le 16$, and employ many more routine estimations
(see Sections 19-24). Thus 
we are able to rewrite (5.28) in the following form
(see (12.27)):
$$OD_{\kappa ;N}(M)=
OD_{\kappa ;N}(M;{\rm Core};\clubsuit\clubsuit
\clubsuit\heartsuit ;1;{\rm
DominatingPart})+{\rm\ negligible},$$
where
$$OD_{\kappa ;N}(M;{\rm Core};\clubsuit\clubsuit
\clubsuit\heartsuit ;1;{\rm DominatingPart})=$$
$$=\prod_{p}{p^3(p+4)\over (p+1)^4}\cdot
{(6/\pi^2)^4\over 2}
\prod_{q|Z_0}{q\over q+4}\cdot
{\pi\over 2}\cdot\left(
\sum_{j=1}^{D^4}\chi_{-D}(j)\log j\right)\cdot $$
$$\cdot\sqrt{D}
\sum_{\scriptstyle 1\le r\le D^8:\atop\scriptstyle \gcd (r,Z_0)=1}
{\mu(r)\chi_{-D}(r)\tau (r)\over r\ph (r)}\prod_{q|r}{q\over q+4}
\sum_{r_3|r}\mu (r_3)\chi_{-D}(r_3)
\cdot $$
$$\cdot \sum_{\scriptstyle 1\le d_4\le M:\atop\scriptstyle 
\gcd (d_4,rZ_0)=1}
{\mu^2 (d_4)\chi_{-D}(d_4)
\over\ph (d_4)}\prod_{p|d_4}{p(p-2)\over (p+4)(p-1)}
\sum_{\scriptstyle 1\le d_3\le M/d_4:\atop\scriptstyle \gcd 
(d_3,rd_4Z_0)=1}
{\mu (d_3)\over\ph (d_3)}
\prod_{q|d_3}{q(q-2)\over (q+4)(q-1)}\cdot $$
$$\cdot \sum_{\scriptstyle 1\le d_2\le M/(d_4d_3):
\atop\scriptstyle \gcd (d_2,rd_4d_3Z_0)=1}
{\mu^2(d_2)\over\ph (d_2)}\prod_{p|d_2}{p\over p+4}\cdot 
\Biggl( \left(2\log d_4+4\log d_3+ 4\log d_2\right)
\log {M^2\over N} $$
$$ -4\log d_4\log d_3-
8\log d_3\log d_2-4\log d_4\log d_2-4(\log d_2)^2\Biggr)
+{\rm negligible}.\eqno(5.29)$$

The evaluation/estimation of the complicated sum (5.29)
is the third main step in the proof of Lemma 5.5.
We repeatedly use the
simple fact that $\mu (r)\chi_{-D}(r)=\mu^2 (r)$ if $r\le D$ and  
$\gcd (r,Z_0)=1$, and also use the less trivial facts
$$\sum_{n=1}^{\infty}{\mu (n)\over n}=0,$$
$$\sum_{n=1}^{\infty}{\mu (n)\log n\over n}=-1,$$
see the
{\it Guiding Intuition} at the end of Section 12.
Thus we obtain---again via a long chain of estimations that
prove the error term to be  {\it negligible}---the
last approximation
 (see Lemma 13.2 and Lemma 14.1):
$$OD_{\kappa ;N}(M)=$$
$$=\sqrt{D}{\pi\over 2}\left(
\sum_{j=1}^{D^4}\chi_{-D}(j)\log j\right)\log {M^2\over N}\cdot 
{6\over \pi^2}\prod_{q|Z_0}{(q-1)^2\over q(q+1)}+{\rm\ negligible},
\eqno(5.30)$$

Finally, to evaluate the critical sum
$\sum_{j=1}^{D^4}\chi_{-D}(j)\log j$ in (5.30), we use Facts (5) and (6),
and Lemma 5.5 follows. (The proofs of Facts (5) and (6) are in Section 29.)

We may say (with some gross oversimplification)
that the rest of the proof is just routine 
elementary calculations/estimations
using the listed ingredients in the outlined way.

\bigskip\bigskip\bigskip\bigskip
\medskip
\centerline{\bf 6. Theorem 2: the case of positive discriminants (I)}
\bigskip\bigskip\bigskip\medskip\noindent
{\bf Outline of the  proof of Theorem 2.}
Let $D>0$ denote a positive fundamental discriminant
violating Theorem 2.
The proof is similar to the proof of Theorem 1, but it requires
 several  modifications.
The first change comes from the fact that 
in the proof of Lemma 1.1 
 a positive definite binary quadratic form has positive values
only. In sharp contrast, a binary quadratic form 
of positive discriminant  is indefinite, and has both positive and
negative values.
This is why we replace the elementary Lemma 1.1 with the following 
more sophisticated lemma
(its proof uses the P\' olya--Vinogradov inquality).
Let
$$R(\Delta ;n)=\sum_{i=1}^{h(\Delta )}R(\Delta ;n,F_i)\eqno(6.1)$$ 
denote the total number of primary representations 
 of a given integer $n\ge 1$
by a representative set of binary quadratic form $F(x,y)=ax^2+bxy+cy^2$ of
discriminant $b^2-4ac=\Delta>0$. By Dirichlet's theorem 
$$R(\Delta;n)=\sum_{m|n}\chi_{\Delta}(m).\eqno(6.2)$$
\medskip\noindent
{\bf Lemma 6.1} {\it For every fundamental
discriminant $\Delta>0$ and  integer $N\ge 1$ we have
$$\left|\sum_{n=1}^NR(\Delta;n)-L(1,\chi_{\Delta})N
\right|\le 4\sqrt{N}\Delta^{1/4}
\sqrt{\log \Delta}.$$ }
\medskip
{\bf Proof.} By (6.1) and (6.2),
$$\sum_{n=1}^NR(\Delta;n)=\sum_{m=1}^N\sum_{m|n}\chi_{\Delta}(m)=
\sum_{1\le m\le \sqrt{N}\Delta^{1/4}\sqrt{\log \Delta}}
\chi_{\Delta}(m)\sum_{1\le k\le N/m}1+$$
$$+\sum_{1\le k<{\sqrt{N}\over \Delta^{1/4}\sqrt{\log \Delta}}}\ \ 
\sum_{\sqrt{N}\Delta^{1/4}\sqrt{\log \Delta}<m<N/k}\chi_{\Delta}(m),$$
and combining this with the P\' olya--Vinogradov inquality (see Lemma
9.3), we obtain
$$\left|\sum_{n=1}^NR(\Delta;n)-
\sum_{1\le m\le \sqrt{N}\Delta^{1/4}
\sqrt{\log \Delta}}\chi_{\Delta}(m){N\over m}\right| \le $$
$$\le\sqrt{N}\Delta^{1/4}\sqrt{\log \Delta}+
\sum_{1\le k<{\sqrt{N}\over \Delta^{1/4}\sqrt{\log \Delta}}}
 \sqrt{\Delta}\log \Delta\le 2\sqrt{N}
\Delta^{1/4}\sqrt{\log \Delta}.\eqno(6.3)$$
Clearly
$$\sum_{1\le m\le \sqrt{N}\Delta^{1/4}\sqrt{\log \Delta}}
\chi_{\Delta}(m){N\over m}=
L(1,\chi_{\Delta})N-N\sum_{m> \sqrt{N}
\Delta^{1/4}\sqrt{\log \Delta}}
{\chi_{\Delta}(m)\over m}.\eqno(6.4)$$
Using partial summation and the P\' olya--Vinogradov inquality,
$$\left|\sum_{m>\sqrt{N}\Delta^{1/4}\sqrt{\log \Delta}}
{\chi_{\Delta}(m)\over m}\right| =$$
$$=\left|\sum_{n>\sqrt{N}\Delta^{1/4}\sqrt{\log \Delta}}
\left(\sum_{\sqrt{N}\Delta^{1/4}\sqrt{\log \Delta}<m\le n}
\chi_{\Delta}(m)\right)\left({1\over n}-{1\over n+1}\right) \right|\le $$
$$\le {2\sqrt{\Delta}\log \Delta\over\sqrt{N}
\Delta^{1/4}\sqrt{\log \Delta}}=
{2\Delta^{1/4}\sqrt{\log \Delta}\over\sqrt{N}}
.\eqno(6.5)$$
Combining (6.3)-(6.5), Lemma 6.1 follows.
(Note that we didn't  use the assumption $\Delta >0$.)\proofend\medskip
Through the proof of Theorem 2 Lemma 6.1 plays 
the same role as that of Lemma 1.1
in the proof of Theorem 1.

The second change is due to the fact
$\chi_{D}(-1)=1$ for $D>0$, which implies that
 the interval $[A_1,A_2]$ in Section 1
 is replaced by the symmetric interval $[-A,A]$ (i.e., $A_2=-A_1>0$).
Then for $D>0$ the analog of (2.39) is
$$WOD_a(M;\kappa ;N;d)
=\sum_{{\scriptstyle (m_1,\l_1)\ne 
(m_2,\l_2):\ 1\le m_1,m_2\le M\atop
\scriptstyle 1\le\l_h\le m_h,\gcd (\l_h,m_h)=1,h=1,2}\atop\scriptstyle
\gcd (m_1,m_2)=d,\gcd (m_1m_2,Z_0)=1}{\mu(m_1)\mu(m_2)
\over\ph (m_1)\ph (m_2)}\cdot \qquad\qquad $$
$$\quad\cdot
e^{2\pi{\rm i}a({\l_1\over m_1}-{\l_2\over m_2})}
\cdot
\left( {\sin\left( (2N+1)\pi D({\l_1\over m_1}-{\l_2\over
m_2})\right)\over (2N+1)\sin\left( \pi D({\l_1\over m_1}-{\l_2\over
m_2})\right)}\right)^{\kappa },\eqno(6.6)$$
and the analog of (2.40) is
($\delta =\pm 1$)  
$$\Omega_{\delta ;\kappa ;N}(M)=
\sum_{\scriptstyle 1\le a\le D\atop\scriptstyle \chi_{D}(a)=\delta}
W_a(M;\kappa ;N)=$$
$$=\sum_{k=0}^{\kappa N}w^{(\kappa ,N)}_k
{1\over 2k+1}\sum_{j:\ -k\le j\le k}
\sum_{\scriptstyle 1\le a\le D\atop\scriptstyle \chi_{D}(a)=\delta}$$
$$\cdot\left(
\sum_{\scriptstyle d\ge 1:\gcd (Z_0,d)=1\atop\scriptstyle
d|a+jD}\mu (d)
\sum_{\scriptstyle 1\le k\le M/d:\atop
\scriptstyle \gcd (k,(a+jD)Z_0)=1}
{|\mu (k)|\over\ph (k)}\right)^2.\eqno(6.7)$$
There is no change in (2.41), but the analog of (2.42)
is the following:
$$\Omega_{\delta ;\kappa ;N}({\rm OffDiag};M)=
\sum_{\scriptstyle 1\le a\le D\atop\scriptstyle \chi_{D}(a)=\delta}
{\rm WOffDiag}_a(M;\kappa ;N)=$$
$$=\sum_{{\scriptstyle (m_1,\l_1)\ne 
(m_2,\l_2):\ 1\le m_1,m_2\le M\atop
\scriptstyle 1\le\l_h\le m_h,\gcd (\l_h,m_h)=1,h=1,2}\atop\scriptstyle
\gcd (m_1,m_2)=d,\gcd (m_1m_2,Z_0)=1}{\mu(m_1)\mu(m_2)
\over\ph (m_1)\ph (m_2)}\cdot \qquad\qquad $$
$$\cdot\left(\sum_{\scriptstyle 1\le a\le D\atop\scriptstyle 
\chi_{D}(a)=\delta}
e^{2\pi{\rm i}a({\l_1\over m_1}-{\l_2\over m_2})}\right)
\left( {\sin\left( (2N+1)\pi D({\l_1\over m_1}-{\l_2\over
m_2})\right)\over (2N+1)\sin\left( \pi D({\l_1\over m_1}-{\l_2\over
m_2})\right)}\right)^{\kappa }
.\eqno(6.8)$$
Now we are ready to formulate the analog of the key equality
(2.49). Again the proof is based on the fact that, under the existence
of a ``bad'' positive discriminant $D>0$, the left-hand side and the
right-hand side of the analog of (2.49) are not equal. They are not
equal for the same
reason as in the proof of Theorem 1: 
we prove the analogs of (2.50) and (2.51),
and the underlying intuition is that the right side of the analog of 
(2.49) does, and the left side does not distinguish between the primes
$q|Z_0$ with $\chi_{D}(q)=1$, $0$ or $-1$.

The third change is that in Lemma 3.2 the Gauss sum becomes real:
$G(\chi_{D};D)=\sqrt{D}$ for $D>0$. It follows that for $D>0$ we have to
modify Lemma 3.1 by removing the factor ${\rm i}=\sqrt{-1}$.
Therefore, the analog of the critical sum (3.20)
goes as follows:
$$OD_{\kappa ;N}(M)={\sqrt{D}\over 2}
\sum_{\scriptstyle 1\le d\le M:\ \mu^2(d)=1
\atop\scriptstyle \gcd (d,Z_0)=1}\ph(d)
\sum_{d_1|d}\cdot $$
$$\cdot\prod_{p|d_1}{p-2\over p-1} 
\sum_{\scriptstyle (m_1,m_2):\ 1\le m_1,m_2\le M\atop
\scriptstyle 
\gcd (m_1,m_2)=d,\gcd (m_1m_2,Z_0)=1}{\mu(m_1)\mu(m_2)
\over\ph (m_1)\ph (m_2)}\cdot\chi_{D}({m_1m_2d_1\over d^2})\cdot $$
$$\cdot
\left( 2\sum_{\scriptstyle 1\le s\le {n\over 2}={m_1m_2d_1\over 2d^2}:\atop
\scriptstyle \gcd (s,n)=1}\chi_{D}(s)
\left( {\sin ( (2N+1)\pi s/n)
\over (2N+1)\sin ( \pi s/n)}\right)^{\kappa }\right)
\eqno(6.9)$$
with $s,n,\l$ coming from 
$${\l_1\over m_1}-{\l_2\over
m_2}={\l \over n}{\rm\  with\ } n={m_1m_2d_1\over d^2},$$
$${\rm\ and\ }
s\equiv D\l\
({\rm mod\ } n),
1\le |s|\le n/2, \gcd (s,n)=1,\eqno(6.10)$$

There will be some
change in Section 9. Lemma 9.6 states that
 for every negative fundamental discriminant $\Delta<0$,
$${1\over |\Delta|}\sum_{k=1}^{|\Delta|}S(\Delta;k)
=h(\Delta),\eqno(6.11)$$
where
$$S(\Delta;k)=\sum_{1\le j\le k}\chi_{\Delta}(j)\eqno(6.12)$$
is the character-sum.
In the proof of the case of negative discriminants
we heavily used the fact that 
the average of the  character-sum $h(-D)$
for the ``bad discriminant'' 
was $o(\sqrt{D})$ (since $h(-D)=o(\sqrt{D})$ 
was the hypothesis of Theorem 1).

The good news is that in the
 case of positive discriminants
we have the following {\it stronger} form of (6.11), which makes
the analog arguments simpler:
 for every positive fundamental discriminant $\Delta>0$,
$${1\over \Delta}\sum_{k=1}^{\Delta}S(\Delta;k)=0,\eqno(6.13)$$
where
$$S(\Delta;k)=\sum_{1\le j\le k}\chi_{\Delta}(j)$$
is the character-sum.
Indeed, since $\chi_{\Delta}(-1)=1$ for $\Delta>0$, we have
$$S(\Delta;\Delta-k)=\sum_{j=1}^{\Delta}
\chi_{\Delta}(j)-S(\Delta;k-1)=-S(\Delta;k-1),$$
and so
$$\sum_{k=1}^{\Delta}S(\Delta;k)={1\over 2}
\sum_{k=1}^{\Delta}\left( S(\Delta;k-1)+S(\Delta;\Delta-k)\right) =
0,$$
proving (6.13).

The next major change is that
$\sum_{n=1}^{ND}\chi_{D}(n)\log n$ is {\it not} about
constant times $\sqrt{D}$
for the bad positive discriminant $D>0$ with ``large'' $N$. Indeed,
we have
$$\sum_{n=1}^{ND}\chi_{D}(n)\log n=
O(L(1,\chi_{D})\sqrt{D})+O(D/N),\eqno(6.14)$$
which implies that (6.14) is negligible in the sense
$$\sum_{n=1}^{ND}\chi_{D}(n)\log n=
o(\sqrt{D})\eqno(6.15)$$
if $L(1,\chi_{D})=o(1)$ (and $N\ge D>0$).

To prove (6.14), we need
\medskip\noindent
{\bf Lemma 6.2} {\it For every fundamental
discriminant $\Delta>0$,
$$\lim_{N\to\infty }\sum_{n=1}^{N\Delta}\chi_{\Delta}(n)
\log n=-L'(0,\chi_{\Delta})
,\eqno(6.16)$$
where of course $L'(0,\chi_{\Delta})$ denotes the derivative of
the L-function $L(s;\chi_{\Delta})$ at $s=0$ (note that
 $L(s;\chi_{\Delta})$ is regular on the entire
complex plane).

Moreover, we can estimate the speed of convergence as follows:
$$\left|\sum_{n=1}^{N\Delta}\chi_{\Delta}(n)\log n
+L'(0,\chi_{\Delta})
\right|\le {2\Delta\over N}.\eqno(6.16')$$}\medskip
We postpone the proof of Lemma 6.2 to the end of Section 15.

To find a connection between $L'(0,\chi_{D})$ and $L'(1,\chi_{D})$,
we use the Functional Equation of $L(s,\chi_{D})$ for the ``bad''
discriminant $D>0$:
$$L(1-s,\chi_{D})\Gamma\left( {1-s\over 2}\right)\left(
{D\over\pi}\right)^{(1-s)/2}=L(s,\chi_{D})\Gamma\left( {s\over 2}\right)\left(
{D\over\pi}\right)^{s/2}.$$
In fact, we use the following  equivalent form:
$$L(1-s,\chi_{D})=L(s,\chi_{D})\Gamma\left( {s\over 2}\right)
{1\over \Gamma\left( (1-s)/2\right) }
\left( {D\over\pi}\right)^{s-{1\over 2}}.\eqno(6.17)$$
Note that the reciprocal of the gamma function
 $1/\Gamma (s)$ is regular on the entire complex plane,
and has a simple zero at $s=0$, so
differentiating (6.17) at $s=1$, we have
$$-L'(0,\chi_{D})=L'(1,\chi_{D})\Gamma\left( {1\over 2}\right)
{1\over\Gamma (0)}\left({D\over\pi}\right)^{1/2}+$$
$$+L(1,\chi_{D})\Gamma'\left( {1\over 2}\right)
\Gamma^{-1}(0)\left(
{D\over\pi}\right)^{1/2}-
L(1,\chi_{D})\Gamma\left( {1\over 2}\right)
{d\over ds}{1\over\Gamma (s)}\bigm|_{s=0}\left(
{D\over\pi}\right)^{1/2}+$$
$$+
L(1,\chi_{D})\Gamma\left( {1\over 2}\right)
{1\over\Gamma (0)}\left(
{D\over\pi}\right)^{1/ 2}\log (D/\pi )=$$
$$=
-L(1,\chi_{D})\Gamma\left( {1\over 2}\right)
{d\over ds}{1\over\Gamma (s)}\bigm|_{s=0}\left(
{D\over\pi}\right)^{1/2},\eqno(6.18)$$
since the rest of the terms contain the factor
$1/\Gamma (0)=0$.

Now (6.18) gives (6.14).
 But this ``idiosyncrasy'' of the case of
 positive discriminants  actually helps us; 
we will take advantage of it.

Since the first derivative of (6.17) turned out to be ``negligible'', 
we study the second derivative of (6.17) at $s=1$:
ignoring  the terms that contain the factor
$1/\Gamma (0)=0$, we have
$$L''(0,\chi_{D})=-L'(1,\chi_{D})\Gamma\left( {1\over 2}\right)
{d\over ds}{1\over\Gamma(s)}\bigm|_{s=0}
\left({D\over\pi}\right)^{1/2}$$
$$-
L(1,\chi_{D})\Gamma'\left( {1\over 2}\right)
{d\over ds}{1\over\Gamma (s)}\bigm|_{s=0}\left(
{D\over\pi}\right)^{1/2}+
L(1,\chi_{D})\Gamma\left( {1\over 2}\right)
{d^2\over ds^2}{1\over\Gamma (s)}\bigm|_{s=0}\left(
{D\over\pi}\right)^{1/2}$$
$$-
L(1,\chi_{D})\Gamma\left( {1\over 2}\right)
{d\over ds}{1\over\Gamma (s)}\bigm|_{s=0}\left(
{D\over\pi}\right)^{1/ 2}\log (D/\pi ).\eqno(6.19)$$
Since ${d\over ds}{1\over\Gamma (s)}\bigm|_{s=0}=1$
and $\Gamma (1/2) =\sqrt{\pi}$, by (6.19),
$$\left| L''(0,\chi_{D}) +L'(1,\chi_{D})\sqrt{D}\right|\le $$
$$\le 100L(1,\chi_{D})\sqrt{D}\log D.\eqno(6.20)$$
Note that for positive discriminants $D>0$
Goldfeld gave an explicit formula for
$L'(1,\chi_{D})$ (see equations (12), (13) and (14) in [Go]) 
 that goes as follows:
$$L'(1,\chi_{D})={\pi^2\over 6}\sum_{(a,b,c)}\left( {1\over
a}+Q(a,b,c)\right)\ +\ O(h(D)/\sqrt{D})+
O(L(1,\chi_{D})\log D),\eqno(6.21)$$
where $Q(a,b,c)\ge 0$ is a  constant that depends only on
the quadratic irrational
$${-b+\sqrt{D}\over 2a}={-b+\sqrt{b^2-4ac}\over 2a},\eqno(6.22)$$
where $\sum_{(a,b,c)}$ is taken over all
reduced, primitive, inequivalent binary quadratic forms
$ax^2+bxy+cy^2$ of 
discriminant $D>0$ (so there are $h(D)$ triples $(a,b,c)$).
The implicit constants in the $O$-notation in (6.21) are 
effectively computable, and also
Goldfeld gave an explicit form of $Q(a,b,c)$ in terms of  the
continued fraction expansion of the number (6.22).

In a related paper, Goldfeld and Schinzel [Go-Sch] introduced 
and studied the other sum
$$\sum_{(A,B,C)}{1\over A},\eqno(6.23)$$
where $\sum_{(A,B,C)}$ is taken over all binary quadratic forms
$Ax^2+Bxy+Cy^2$  of 
discriminant $D>0$ such that
$$-A<B\le A<{1\over 4}\sqrt{D}.\eqno(6.24)$$
(In sharp contrast with the case $D<0$, if the fundamental unit of
the real quadratic field $\Q (\sqrt{D})$ is ``large''
and the class number $h(D)$ is ``small'', then many
different triples $(A,B,C)$ satisfying (6.24) are equivalent, and
belong to the same class represented by one triple $(a,b,c)$.) 
The argument of the paper [Go-Sch] (using the theory of continued fraction)
gives that
$$\sum_{(a,b,c)}\left( {1\over
a}+Q(a,b,c)\right) =\sum_{(A,B,C)}{1\over A}\ +\ {\rm negligible}=
\prod_{\scriptstyle p|Z_0\atop\scriptstyle \chi_{D}(p)\ne -1}
{p+1\over p-1}\prod_{p|D}{p-1\over p}
\ +\ {\rm negligible},$$
if $L(1,\chi_{D})$ is small (which is true for our ``bad'' discriminant
$D>0$).

More precisely, [Go-Sch] yields the following result (an 
analog of Lemma 14.2 later).
\medskip\noindent
{\bf Lemma 6.3} {\it For our ``bad'' fundamental discriminant $D>0$
we have
$$\left|\prod_{\scriptstyle p|Z_0\atop\scriptstyle \chi_{D}(p)\ne -1}
{p+1\over p-1}\prod_{p|D}{p-1\over p}
-\sum_{(a,b,c)}\left( {1\over
a}+Q(a,b,c)\right)\right| \le $$
$$\le 10^3L(1,\chi_{D})+{10^4\over (\log D)^3}.$$}
\medskip

I recall that for $-D<0$ the sum
$$\sum_{j=1}^{D^4}\chi_{-D}(j)\log j\eqno(6.25)$$
played a key role in the proof of Theorem 1 (see Lemma 13.2 later).
For $D>0$ the sum (6.25) is ``negligible'' in the sense
$$\sum_{j=1}^{D^4}\chi_{D}(j)\log j=o(\sqrt{D})$$
(see (6.14)-(6.15)), and the other sum
$$\sum_{j=1}^{D^4}\chi_{D}(j)(\log j)^2\eqno(6.26)$$
will play  the same key role.

We need
 the following analog of Lemma 6.2 for the second derivative.
\medskip\noindent
{\bf Lemma 6.4} {\it For every fundamental
discriminant $\Delta >0$,
$$\lim_{N\to\infty }\sum_{n=1}^{N\Delta}\chi_{\Delta}
(n)(\log n)^2=L''(0,\chi_{\Delta}),\eqno(6.27)$$
where  $L''(0,\chi_{\Delta})$ denotes the second derivative of
the L-function $L(s;\chi_{\Delta})$ at $s=0$.

Moreover, we can estimate the speed of convergence:
$$\left|\sum_{n=1}^{N\Delta}\chi_{\Delta}(n)(\log n)^2-
L''(0,\chi_{\Delta})
\right|\le {8\Delta\log (N\Delta)\over N}.\eqno(6.28)$$}\medskip
We postpone the proof to the end of Section 15.

Combining (6.28), (6.21), (6.20) and Lemma 6.3,
 we conclude Section 6 with
\medskip\noindent
{\bf Lemma 6.5} {\it For our ``bad'' fundamental discriminant $D>0$
we have
$$\left|\sum_{j=1}^{D^4}\chi_{D}(j)(\log j)^2+
{\pi^2\over 6}\sqrt{D}
\prod_{\scriptstyle p|Z_0\atop\scriptstyle \chi_{D}(p)\ne -1}
{p+1\over p-1}\prod_{p|D}{p-1\over p}\right| =$$
$$=O(L(1,\chi_{D})\sqrt{D}\log D)+O(h(D))+
O(\sqrt{D}(\log D)^{-3}),$$
where the implicit constants are effectively computable.}
\medskip
We complete the outline of the proof of Theorem 2 in Section 15.

\bigskip\bigskip\bigskip
\centerline{\nagybf References}
\bigskip\bigskip
\medskip\noindent
[Ba] A. Baker,  Linear forms in the logarithms of algebraic
numbers, Mathematika 13 (1966), 204-216.
\medskip\noindent
[Bo-Sh] Z.I. Borevich and I.R. Shafarevich, {\it Number Theory},
English translation: Academic Press (New York) 1966. 
\medskip\noindent
[Br-Ch] W.E. Briggs and S. Chowla, On discriminants of binary
quadratic forms with a single class in each genus, Canadian Journal
of Mathematics 6 (1954), 463-470.
\medskip\noindent
[Da] H. Davenport, {\it Multiplicative Number Theory}, Second edition,
 revised by H.L. Montgomery, Springer-Verlag 1980.
\medskip\noindent
[De] M. Deuring, Imagin\"ere quadratische Zahlk\"orper mit der
Klassenzahl 1, Math. Zeitschrift  37 (1933), 405-415.
\medskip\noindent
[Go1] D. Goldfeld, An asymptotic formula relating the Siegel zero
 and the class number of quadratic fields,
 Annali della Scuola Normale Superiore di Pisa,
Classe de Science, Ser. 4, 2, No. 4 (1975), 611-615. 
\medskip\noindent
[Go2] D. Goldfeld, The class number of quadratic fields and the
conjectures of Birch and Swinnerton-Dyer, 
Ann. Scuola Norm. Sup. Pisa (4) 3 (1976), 623-663.
\medskip\noindent
[Go3] D. Goldfeld, Gauss' class number problem for imaginary
quadratic fields,
 Bull. Amer. Math. Soc. (1) 13 (1985), 23-37.
\medskip\noindent
[Go-Sch] D. Goldfeld and A. Schinzel, On Siegel's zero,
Ann. Scuola Norm. Sup. Pisa (4) 2 (1975), 571-583.
\medskip\noindent
[Gr-Za] B. Gross and D. Zagier, Heegner points and derivatives of
L-series, Invent. Math. 84 (1986), 225-320.
\medskip\noindent
[Ha-Wr] G. Hardy and E. Wright, {\it An introduction to the theory of
numbers}, 5th edition, Clarendon Press, Oxford, 1979.
\medskip\noindent
[Hee] K. Heegner, Diophantische Analysis und Modulfunktionen,
Math. Zeitschrift  56 (1952), 227-253.
 \medskip\noindent
[Hei] H. Heilbronn, On the class number in imaginary quadratic
fields, Quarterly Journal of Math. Oxford 5 (1934), 150-160.
\medskip\noindent
[He-Li] H. Heilbronn and E.H. Linfoot, On the imaginary quadratic
corpora of class number one, Quarterly Journal of Math. Oxford 5
(1934), 293-301.
\medskip\noindent
[Iw-Ko] H. Iwaniec and E. Kowalski, {\it Analytic Number Theory},
Amer. Math. Soc. Colloquium Publications, vol. 53, AMS Providence 2004.
\medskip\noindent
[Ka] A. Karatsuba, {\it Basic Analytic Number Theory}, 
English translation: Springer-Verlag 1993.
\medskip\noindent
[La] E. Landau, Bemerkungen zum Heilbronnschen Satz, Acta
Arithmetica 1 (1935), 1-18.
\medskip\noindent
[Po] G. P\' olya, \" Uber die Verteilung der quadratischen
Reste und Nichtreste, G\" ottingen Nachrichten (1918), 21-29.
\medskip\noindent
[Re] A. R\' enyi, {\it Probability Theory}, American Elsevier
Publishing Comp. Inc., New York, 1970.
\medskip\noindent
[Se-Ch] A. Selberg and S. Chowla, On Epstein's zeta-function, Journal
f\"ur die reine und angewandte Mathematik, 227 (1967), 86-110.
\medskip\noindent
[Si] C. L. Siegel, \"Uber die Classenzahl quadratischer
Zahlk\"orper, Acta Arithmetica 1 (1935), 83-86.
\medskip\noindent
[St1] H.M. Stark, A complete determination of the complex quadratic
fields of class-number one, Michigan Math. Journal 14 (1967), 1-27.
\medskip\noindent
[St2] H.M. Stark, On the ``gap'' in a theorem of Heegner,
Journal of the Number Theory 1 (1969), 16-27.
\medskip\noindent
[Za] D. Zagier, {\it Zetafunktionen und quadratische K\" orper},
Springer-Verlag 1981.

\bigskip\bigskip\bigskip
jbeck@math.rutgers.edu

\bigskip
J\'ozsef Beck
\medskip
Department of Mathematics, Rutgers University
\medskip
110 Frelinghuysen Road, Hill Center
\medskip
Piscataway, New Jersey, 08854-8019 USA

\end